\documentclass{article}
\pdfoutput=1
\usepackage{amsmath}
\usepackage{amssymb}
\usepackage{latexsym}
\usepackage{bbm}
\usepackage{graphicx}
\usepackage{color}
\usepackage{url}
\usepackage[english]{babel}
\usepackage{endfloat}
\usepackage[pdftex,%
a4paper=true,%
backref=page,%
pagebackref=true,%
bookmarks=true,%
bookmarksnumbered=true,%
colorlinks=true,%
urlcolor=blue,%
pdfstartview=FitH%
]{hyperref}

\newtheorem{thm}{Theorem}

\newtheorem{lem}{Lemma}
\newtheorem{rmq}{Remark}
\newtheorem{dfn}{Definition}
\newtheorem{conj}{Conjecture}

\newcommand{\be}{\begin{equation}}
\newcommand{\ee}{\end{equation}}
\newcommand{\ba}{\begin{array}}
\newcommand{\ea}{\end{array}}
\newcommand{\baa}{\left\{ \begin{array}}
\newcommand{\eaa}{\end{array} \right.}
\newcommand{\bproof}{\paragraph{Proof. }}
\newcommand{\eproof}{$\Box$\\}

\begin{document}

%----------- TITLE --------------%

\title{Experimental Study of the HUM Control Operator for Linear Waves}

\author{Gilles Lebeau\footnote{Laboratoire J.-A. Dieudonn\'e, Parc Valrose, Nice, France, {\sf lebeau@math.unice.fr}} \footnote{Institut Universitaire de France}  \\ Universit\'e de Nice \and Ma\"elle Nodet\footnote{Laboratoire J. Kuntzmann, Domaine Universitaire, Grenoble, France, {\sf maelle.nodet@inria.fr}}  \\ Universit\'e de Grenoble, INRIA}

%-----------------------------------%

\maketitle

\abstract{We consider the problem of the numerical approximation of the linear controllability of waves. All our experiments are done in a bounded domain $\Omega$ of the plane, with Dirichlet boundary conditions and internal control. We use a  Galerkin approximation of the optimal control operator of the continuous model, based on the spectral theory of the Laplace operator in $\Omega$. This allows us to obtain surprisingly good illustrations of the main theoretical results available on the controllability of waves, and
to formulate some questions for the future analysis of optimal control theory of waves.}

\newpage
\tableofcontents
\newpage

%%%%%%%%%%%%%%%%%%%%%%%%%%%%
\section{Introduction}\label{sec1}
This paper is devoted to the experimental study of the exact controllability of waves. All our experiments will be done in a bounded domain $\Omega$ of the plane, 
with Dirichlet boundary conditions and
with  internal control. We use the most natural  approach for the numerical computation:
a Galerkin approximation of the optimal control operator 
of the continuous model based on the spectral theory of the Laplace operator in $\Omega$. This will allow us to obtain surprisingly good
illustrations of the main theoretical results available on the controllability of waves, and
to formulate some questions for the future analysis of the optimal control theory of waves.\\

The problem of controllability for linear evolution equations and systems has a long story
for which we refer to the review of D.-L. Russel in \cite{Ru} and to the book of J. L. Lions
\cite{Lions88}. Concerning controllability of linear waves, the main theoretical result 
is the so called ``Geometric Control Condition" of C. Bardos, G. Lebeau and J. Rauch
 \cite{BLR92}, {\bf GCC} in short,  which gives a (almost) necessary and sufficient condition for exact controllability. This is a ``geometrical optics" condition on the behavior of
 optical rays inside $\overline\Omega$. Here, optical rays are just straight lines inside $\Omega$, 
 reflected at the boundary according to the Snell-Descartes law of reflection. The precise definition of
 optical rays near points of tangency with the boundary is given in the works of R. Melrose
 and J. Sj\" ostrand  \cite{MS1} and \cite{MS2}. For internal control, {\bf GCC} asserts that waves in 
  a regular bounded domain $\Omega$ are exactly controllable by control functions
supported  in the closure of an open sub-domain $U$ and acting during a time $T$,  if (and only if, if one allows arbitrary small perturbations of the time control $T$ and of the control domain $U$)\\

{\bf GCC}: every optical ray
of length $T$ in $\Omega$  enters the 
sub-domain $U$. \\

Even when {\bf GCC} is satisfied, the numerical computation of the control is not an easy task. 
The original approach consists first in discretizing the continuous model, and then in computing
the control of the discrete system to use it as a numerical approximation of the continuous one. This 
method has been developed by R. Glowinski et al. (see \cite{GloLi2} and \cite{GlowLionsHe08}), and
used for numerical experiments in \cite{AL98}. However, as observed in the first works of
R. Glowinski, interaction of waves with a numerical mesh produces spurious high frequency oscillations. In fact, the discrete model is not uniformly exactly controllable when the mesh size goes to zero, since the group velocity  converges to zero when the solutions wavelength is comparable to the mesh size. In other words, the processes of
numerical discretization and observation or control do not commute. A precise analysis of this
lack of commutation and its impact on the computation of the control has been done by E. Zuazua in
\cite{Zuazua02} and \cite{Zuazua05}. \\

In this paper, we shall use another approach, namely
we will discretize the optimal control of the continuous model using a projection of the
wave equation onto the finite dimensional space spanned by the eigenfunctions $e_j$ of the Laplace operator in $\Omega$  
with Dirichlet boundary conditions, $-\triangle e_j=\omega_j^2 e_j, e_j\vert_{\partial\Omega=0}$
, with $\omega_j\leq \omega$. Here, $\omega$ will be a cutoff frequency at our disposal. We prove 
(see lemma \ref{lem2} in section \ref{sec2.3}) 
that when {\bf GCC} is satisfied, our numerical control converges when 
$\omega\rightarrow\infty$ to the optimal control of the continuous model . Moreover, 
when {\bf GCC} is not satisfied, we will do
experiments  and we will see an exponential blow-up in the cutoff frequency $\omega$ of the norm of the discretized optimal control. These blow-up rates
will be compared to theoretical results in section \ref{sec2.3}.\\

The paper is organized as follows.\\

Section \ref{sec2} is devoted to the analysis of the optimal control operator for waves
in a bounded regular domain of $\mathbb R^d$.
In section \ref{sec2.1}, we recall the definition  of the optimal control operator $\Lambda$.
In section \ref{sec2.2}, we recall some known theoretical results on $\Lambda$: existence, regularity
properties and the fact that it preserves the frequency localization. We also state our first conjecture,
namely that the optimal control operator $\Lambda$ is a microlocal operator. 
In section \ref{sec2.3}, we introduce the spectral Galerkin approximation $M_{T,\omega}^{-1}$ of $\Lambda$,
where $\omega$ is a cutoff frequency. We prove the convergence of $M_{T,\omega}^{-1}$ toward $\Lambda$
when $\omega\rightarrow \infty$, and we analyze the rate of this convergence. We also state our 
second conjecture on the blow-up rate of $M_{T,\omega}^{-1}$ when {\bf GCC} is not satisfied.
Finally, in section \ref{sec2.4}, we introduce the basis of the energy space in which we compute the matrix of
the operator $M_{T,\omega}$.\\

Section \ref{sec3} is devoted to the experimental validation of our Galerkin approximation.
In section \ref{sec3.1} we introduce the $3$ different domains of the plane for our experiments: square, disc and
trapezoid. In the first two cases, the geodesic flow is totally integrable and the exact eigenfunctions
and eigenvalues of the Laplace operator with Dirichlet boundary condition are known. This is not the case
for the trapezoid. In section \ref{sec3.2smoothing} we introduce the two different choices of the control operator we use
in our experiments. In the first case (non smooth case), we use  $\chi(t,x)={\bf 1}_{[0,T]}{\bf 1}_{U}$.
In the second case (smooth case), we use a suitable regularization of the first case (see formulas
(\ref{mneq:8}) and (\ref{mneq:8bis})). Perhaps the main contribution of this paper is to give an
experimental evidence that the choice of a smooth enough control operator is the right way to get
accuracy on control computing. In section \ref{sec3.3}, in the two cases of the square and  the disc,
we compare the exact eigenvalues with the eigenvalues computed using 
the $5$-points finite difference approximation of the Laplace operator. In section \ref{sec3.4recerror}, formula
(\ref{mneq:10}), we define the
reconstruction error of our method. Finally, in section \ref{sec3.5}, in the case of the square geometry,
we compare the control function (we choose to reconstruct a single eigenvalue) when the eigenvalues and eigenvectors are computed either with finite differences or with exact formulas, and we study the experimental convergence of our numerical  optimal control to the exact optimal control of the continuous model
when the cutoff frequency goes to infinity.\\

The last section, \ref{sec4}, presents various numerical experiments which illustrate the theoretical results
of section \ref{sec2}, and are in support of our two conjectures. In section \ref{sec4.1}, our experiments
illuminate the fact that the optimal control operator preserves the frequency localization, and that
this property is far much stronger when the control function is smooth, as predicted by theoretical results.
In section \ref{sec4.2}, we use Dirac and box experiments to illustrate the fact that the optimal control operator
shows a behavior very close to the behavior of a pseudo-differential operator: this is in support
of our first conjecture. In section \ref{sec4.3}, we plot the reconstruction error as a function of the cutoff
frequency: this illuminates how the rate of convergence of the Galerkin approximation depends on the
regularity of the control function. In section \ref{sec4.4}, we present various results on the energy of the control function. In section \ref{sec4.5}, we compute the condition number of the matrix $M_{T,\omega}$ for a given
control domain $U$, as a function of the control time $T$ and the cutoff frequency $\omega$. In particular,
figures \ref{fig:cond-eigenvalues-constant} and \ref{fig:cond-time-constant} are in support of our second conjecture on the blow-up rate of
$M_{T,\omega}^{-1}$ when {\bf GCC} is not satisfied. In section \ref{sec4.6}, we perform experiments in the disc
when {\bf GCC} is not satisfied, for two different data: in the first case, every optical rays of length 
$T$ starting at a point where the data is not small enters the control domain $U$, and we observe
a rather good reconstruction error if the cutoff frequency is not too high. In the second case, there
exists an optical ray starting at a point where the data is not small and which never enters the control domain,
and we observe a very poor reconstruction at any cutoff frequency. This is a fascinating phenomena which
has not been previously studied in theoretical works. It will be of major practical interest
to get quantitative results on the best cutoff frequency which optimizes the reconstruction error
(this optimal cutoff frequency is equal to $\infty$ when {\bf GCC} is  satisfied,
the reconstruction error being equal to $0$ in that case), and to  estimate  
the reconstruction error at the optimal cutoff frequency. Clearly, our experiments indicate that
a weak Geometric Control Condition associated to the data one wants to reconstruct will enter
in such a study. \\

%%%%%%%%%%%%%%%%%%%%%%%%%%%%
\section{The  analysis of the optimal control operator}\label{sec2}

\subsection{The optimal control operator}\label{sec2.1}

Here we recall the basic facts we will need in our study 
of the optimal control operator for linear waves. For more details on the HUM
method, we refer to the book of J.-L. Lions \cite{Lions88}.\\

In the framework of the wave equation  in a bounded  open subset $\Omega$ of
$\mathbb{R}^{d}$ with boundary Dirichlet condition, and for internal control the problem of controllability is stated in the following way.
Let $T$ be a
positive time, $U$ a non void open subset of $\Omega$, and $\chi(t,x)$ as follows:
\be\label{1b}
\chi(t,x)=\psi(t)\chi_0(x)
\ee
where $\chi_0$ is a real $L^{\infty}$ function on $\overline{\Omega}$,  such that 
support$(\chi_{0})=\overline U$ and $\chi_0(x)$ is 
continuous and positive for $x\in U$,
$\psi\in C^{\infty}([0,T])$ and $\psi(t)>0$ on
$]0,T[.$ 
  For a
given $f=(u_{0},u_{1})\in H_0^{1}(\Omega)\times L^{2}(\Omega),$ the problem is to find a
source $v(t,x)\in L^{2}(0,T;L^{2}(\Omega))$  such that the solution of the system
\be\label{1.16bis}
\left\{
\begin{array}
[c]{c}%
\square u=\chi v\quad\text{in }]0,+\infty\lbrack\times\Omega\\
u_{\mid\partial\Omega}=0,\quad t>0\\
(u\vert_{t=0},\partial_{t} u\vert_{t=0})=(0,0)
\end{array}
\right. %
\ee
reaches the state $f=(u_{0},u_{1})$ at time $T$. We first rewrite 
the wave operator in (\ref{1.16bis}) as a first order system.
Let $A$ be the matrix
\be\label{1.19}
iA=\begin{pmatrix}
0 & Id \\
\triangle & 0
\end{pmatrix}
\ee
Then $A$ is a unbounded self-adjoint operator on 
$H=H_0^{1}(\Omega)\times L^{2}(\Omega)$,  where the scalar product on 
$H_0^{1}(\Omega)$ is $\int_{\Omega} \nabla u \overline{\nabla v}dx$ and   $D(A)=
\{\underline u=(u_0,u_1) \in H, \ A(\underline u)\in H, \ u_0\vert_{\partial \Omega}=0 \}$. Let $\lambda= \sqrt {-\triangle_D}$
where $-\triangle_D$ is the canonical isomorphism from $H_0^{1}(\Omega)$ onto 
$H^{-1}(\Omega)$. Then $\lambda$ is an isomorphism from 
$H_0^{1}(\Omega)$ onto 
$L^{2}(\Omega)$. 
The operator $B(t)$ given by
\be\label{1.21}
B(t)=\begin{pmatrix}
0 & 0 \\
\chi(t,.)\lambda & 0
\end{pmatrix}
\ee
is bounded on $H$, and one has 
\be\label{1.22}
B^*(t)=\begin{pmatrix}
0 & \lambda^{-1} \chi(t,.) \\
0 & 0
\end{pmatrix}
\ee
The  system (\ref{1.16bis}) is  then equivalent to 
\be\label{1.1}
(\partial_t-iA)f=B(t)g, \quad f(0)=0
\ee
 with $f=(u,\partial_t u )$, 
 $g=(\lambda^{-1}v,0)$. For any $g(t)\in L^1([0,\infty[, H)$,
 the evolution equation
\be\label{1.1bis}
(\partial_t-iA)f=B(t)g, \quad f(0)=0
\ee
admits a unique solution $f=S(g)\in C^0([0,+\infty[, H)$ given by the Duhamel formula
\be\label{1.2}
f(t)=\int_0^t e^{i(t-s)A} B(s)g(s) ds
\ee
Let $T>0$ be given. Let $\mathcal R_T$ be the  reachable set at time $T$ 

\be\label{1.3}
\mathcal R_T =\{f\in H, \  \exists g\in L^2([0,T], H), \  f=S(g)(T) \}
\ee
Then  $\mathcal R_T$ is a linear subspace of $H$, and is the set of states of the system
that one can reach in time $T$, starting from rest,  with the action of an $L^2$ source $g$ filtered by the control 
operator $B$. The control problem consists in giving an accurate description of   $\mathcal R_T$, and exact
controllability is equivalent to  the equality $\mathcal R_T=H$. Let us recall some basic facts.\\

Let $\mathcal H=L^2([0,T], H)$. Let $\mathcal F$ be the closed subspace of  $\mathcal H$
spanned by solutions of the adjoint evolution equation
\be\label{1.4}
\mathcal F =\{h\in \mathcal H, \  (\partial_t-iA^*)h=0, \  h(T)=h_T\in H\}
\ee
Observe that, in our context,  $A^*=A$, and the function $h$ in (\ref{1.4}) is given by
$h(t)=e^{-i(T-t)A}h_T$. Let $\mathcal B^*$ be the adjoint of the operator $g\mapsto S(g)(T)$. Then $\mathcal B^*$ is the bounded operator from $H$ into $\mathcal H$ defined  by
\be\label{1.5}
\mathcal B^*(h_T)(t)=B^*(t)e^{-i(T-t)A}h_T
\ee
For any $g\in L^2([0,T], H)$, one has, with  $f_T=S(g)(T)$
and $h(s)=e^{-i(T-s)A}h_T$
the fundamental identity
\be\label{1.6}
(f_T\vert h_T)_H=\int_0^T(B(s)g(s)\vert h(s))ds=(g\vert \mathcal B^*(h_T))_\mathcal H
\ee
From (\ref{1.6}), one gets easily that the following holds true
\be\label{1.7}
\mathcal R_T \ \text{is a dense subspace of } \  H \  \Longleftrightarrow \  \mathcal B^* \  \text{is an injective operator}
\ee 
which shows that approximate controllability is equivalent to a uniqueness result on 
the adjoint equation. Moreover, one gets from (\ref{1.6}), using  the Riesz and closed graph theorems, that
the following holds true
\be\label{1.8}
\mathcal R_T = H \quad  \quad  \Longleftrightarrow 
\  \exists C, \quad  \Vert h\Vert \leq C\Vert\mathcal B^*h\Vert
 \quad \forall h\in H
\ee 
This is an  observability inequality, and $\mathcal B^*$ is called the
observability operator.  We rewrite 
the observability inequality (\ref{1.8})
 in a more explicit form
\be\label{1.10}
 \exists C, \quad  \Vert h\Vert^2_H \leq  C\int_0^T\Vert B^*(s)e^{-i(T-s)A}h\Vert_H^2 ds
 \quad \forall h\in H
\ee 

Assuming that (\ref{1.10}) holds true, then $\mathcal R_T = H $, $\textrm{Im}(\mathcal B^*)$ is a closed subspace of 
$\mathcal H$, and $\mathcal B^*$ is an isomorphism
of $H$ onto $\textrm{Im}(\mathcal B^*)$. For any $f\in H$, let
$\mathcal C_f$ be the set of control functions $g$ driving $0$ to $f$ in time $T$
\be\label{1.11}
\mathcal C_f=\left\{g\in L^2([0,T],H), \quad f=\int_0^T e^{i(T-s)A}B(s)g(s)ds\right\}
\ee 
From (\ref{1.6}), one gets
\be\label{1.12}
\mathcal C_f=g_0+(\textrm{Im}\mathcal B^*)^\perp, \quad  g_0\in \textrm{Im}\mathcal B^*\cap\mathcal C_f
\ee 
and $g_0=\mathcal B^*h_T$ is the optimal control in the sense that 
\be\label{1.13}
\min \{\Vert g \Vert_{L^2([0,T],H)}, \  g\in \mathcal C_f\} \quad \text{is achieved at} \  g=g_0
\ee 
Let $\Lambda: H\rightarrow H, \ \Lambda(f)=h_T$ be the control map, 
so that the optimal control $g_0$ is equal to $g_0(t)=B^*(t)e^{-i(T-t)A}\Lambda(f)$. 
Then $\Lambda$ is exactly the inverse
of the map $M_T: H\rightarrow H$ with
\be\label{1.14}
\ba{rcl}
M_T &=&\displaystyle \int_0^T m(T-t) dt =\int_0^T m(s) ds \medskip\\
m(s) &=&e^{isA}B(T-s)B^*(T-s)e^{-isA^*}
\ea
\ee 
Observe that  $m(s)=m^*(s)$ is a bounded, self-adjoint, non-negative operator on $H$. Exact controllability is thus equivalent to
\be\label{1.15bis}
\exists C>0, \quad M_T=\int_0^T e^{i(T-t)A} \begin{pmatrix}
0 & 0 \\
0 & \chi^2(t,.)
\end{pmatrix}
e^{-i(T-t)A} dt \geq C Id
\ee  
With $\Lambda=M_T^{-1}$, the optimal
control is then given by $g_0(t)=B^*(t)e^{-i(T-t)A}\Lambda(f)$ and is by (\ref{1.22})
of the form $g_0=(\lambda^{-1}\chi \partial_t w,0)$
where $w(t)=e^{-i(T-t)A}\Lambda(f)$ is the solution of 
\be\label{1.17bis}
\left\{
\begin{array}
[c]{c}%
\square w=0\quad\text{in }\mathbb R\times\Omega, \quad
w_{\mid\partial\Omega}=0\\
(w(T,.),\partial_t w(T,.))=\Lambda(f)
\end{array}
\right. %
\ee

Thus, the optimal control 
function $v$ in (\ref{1.16bis}) is equal to $v=\chi \partial_t w$,
where $w$ is  the  solution
of the dual problem (\ref{1.17bis}). The operator $\Lambda=M_T^{-1}$, with $M_T$
given by (\ref{1.15bis}) is called the optimal control operator.

%%%%%%%%%%%%%%%%%%%%%%%%%%%%%

\subsection{Theoretical results}\label{sec2.2}

\begin{figure}
\begin{center}
\includegraphics[width=\textwidth]{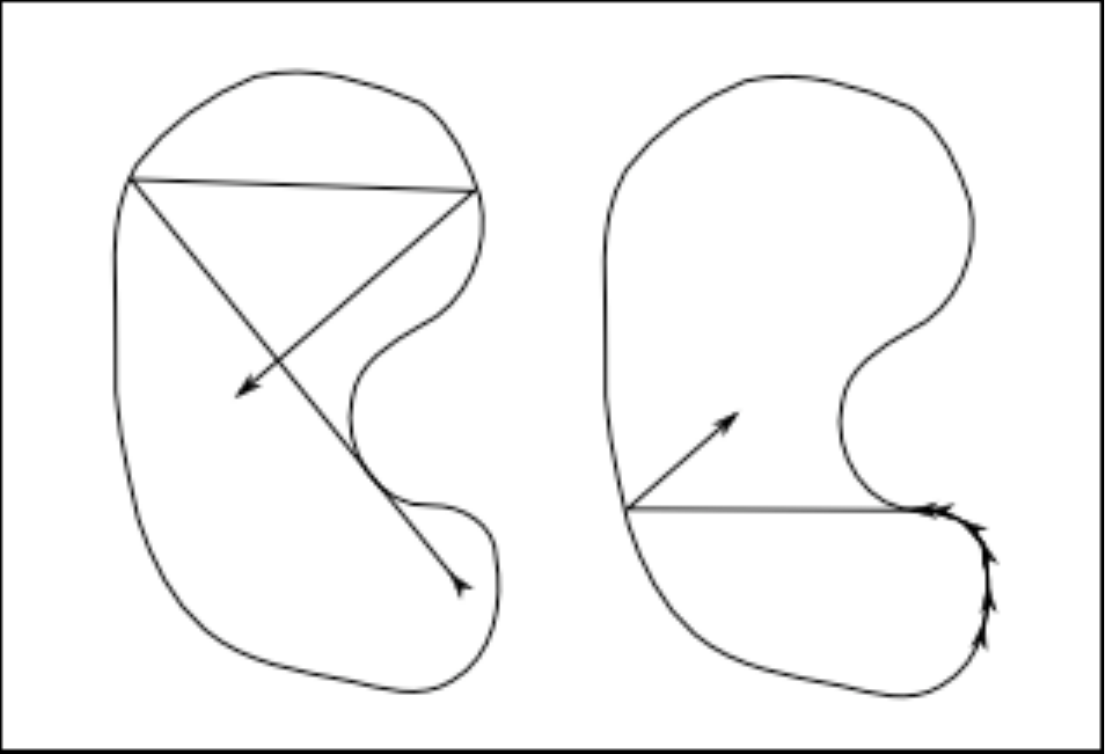}
\caption{\label{fig:rayon} 
Example of optical rays.}
\end{center}
\end{figure}

In this section we  recall some theoretical results on the analysis of the optimal 
control operator $\Lambda$. We will assume here that $\Omega$ is a bounded open subset of
$\mathbb R^d$ with smooth boundary $\partial\Omega$, and that any straight line in $\mathbb R^d$ has only finite order of contacts with the boundary. In that case, optical rays are uniquely defined. See an example of such rays in figure \ref{fig:rayon}.\\

Let $M=\Omega\times\mathbb R_t$. The phase space  
is 

$$
^bT^* M=T^*\overline M\setminus T^*_{\partial M}\simeq T^*M \cup T^*\partial M.
$$

The characteristic variety of the wave operator is the closed subset $\Sigma$ of $^bT^* M$
of points $(x,t,\xi,\tau)$
such that $\vert \tau\vert=\vert\xi\vert$ when $x\in \Omega$ and $\vert \tau\vert\geq \vert\xi\vert$ 
when $x\in \partial\Omega$. Let $^bS^*\overline\Omega$ be the set of points 

$$
^bS^*\overline\Omega = \{ (x_0,\xi_0), \textrm{ with }\vert \xi_0\vert=1 \textrm{ if } x_0\in \Omega,   \vert \xi_0\vert\leq 1 \textrm{ if } x_0\in \partial\Omega\}
$$

For $\rho_0=(x_0,\xi_0)\in {^bS}^*\overline\Omega$, and $\tau=\pm 1$, we shall denote by 
$s \rightarrow (\gamma_{\rho_0}(s), t-s\tau, \tau), s\in \mathbb R$ the generalized bicharacteristic
ray of the wave operator, issued from $(x_0,\xi_0,t,\tau)$. For the construction of the Melrose-Sj\" ostrand flow, we refer to  \cite{MS1}, \cite{MS2} and  to \cite{Hormander85}, vol 3, chapter XXIV . Then
$s\rightarrow \gamma_{\rho_0}(s)=(x(\rho_0,s),\xi(\rho_0,s))$ is the optical ray
starting at  $x_0$ in the direction $\xi_0$. When $x_0\in \partial\Omega$ and $\vert \xi_0\vert< 1$,
then the right (respectively left) derivative of $x(\rho_0,s)$ at $s=0$ is equal to the  unit vector
in $\mathbb R^d$ which projects on $\xi_0\in T^*\partial\Omega$ and which points inside (respectively outside)
$\Omega$. In all other cases, $x(\rho_0,s)$ is derivable at $s=0$ with derivative equal to $\xi_0$.\\

We first recall the theorem of \cite{BLR92}, which gives the existence of the operator $\Lambda$:
\begin{thm}\label{thmblr}
If  the geometric control condition {\bf GCC}
holds true, then $M_T$ is an isomorphism.
\end{thm} 

Next, we recall some new theoretical results obtained in \cite{DehmanLebeau07}.
For these results, the choice of the control function $\chi(t,x)$ in (\ref{1b}) will be essential.
\begin{dfn}\label{definsmooth}
The control function $\chi(t,x)=\psi(t)\chi_0(x)$ is \emph{smooth} if $\chi_0\in C^\infty(\overline{\Omega})$
and $\psi(t)$ is flat at $t=0$ and $t=T$.
\end{dfn}
 For $s\in \mathbb R$, we  denote by $H^s(\Omega,\triangle)$ the domain of the operator $(-\triangle_{\text{Dirichlet}})^{s/2}$. One has $H^0(\Omega,\triangle)=L^2(\Omega)$,
$H^1(\Omega,\triangle)=H_0^1(\Omega)$, and if $(e_j)_{j\geq 1}$ is an $L^2$ orthonormal basis of
eigenfunctions of $-\triangle$ with Dirichlet boundary conditions, $-\triangle e_j=\omega_j^2 e_j$,
$0<\omega_1\leq \omega_2\leq ...$, one has

\be\label{t1}
H^s(\Omega,\triangle)=\{f\in \mathcal D'(\Omega), f=\sum_jf_je_j, \ \sum_j \omega_j^{2s}\vert f_j\vert^2<\infty\}
\ee

The following result of \cite{DehmanLebeau07} says that, under the hypothesis that 
the control function $\chi(t,x)$ is smooth, the optimal control operator 
$\Lambda$ preserves the regularity:
\begin{thm}\label{thm1}
Assume that the geometric control condition {\bf GCC}
holds true, and that the control function $\chi(t,x)$ is smooth. Then the optimal control operator 
$\Lambda$ is an isomorphism of $H^{s+1}(\Omega,\Delta)\oplus H^{s}(\Omega,\Delta)$ for all $s\geq 0$. 
\end{thm}

Observe that theorem \ref{thmblr} is a particular case of theorem \ref{thm1} with $s=0$.
In our experimental study, we will see in section \ref{sec4.3} that the regularity of the control function 
$\chi(t,x)$ is not only a nice hypothesis to get theoretical results. It is also  very efficient
to get accuracy in the numerical computation of the control function. In other words, the
usual choice of the control function $\chi(t,x)={\bf 1}_{[0,T]}{\bf 1}_U$ is a very poor idea to compute
a control. \\

The next result 
says that the optimal control operator 
$\Lambda$ preserves the frequency localization. To state this result
 we briefly introduce the material needed for the Littlewood-Paley decomposition.
Let $\phi$ $\in$ $C^{\infty}([0,\infty[),$ with $\phi(x)=1$ for $\vert x \vert\leq 1/2$ and 
$\phi(x)=0$ for $\vert x \vert\geq 1$. Set $\psi(x)=\phi(x)-\phi(2x)$.
Then  $\psi\in
C_{0}^{\infty}(\mathbb{R}^{\ast})$,  $\psi$ vanishes outside $[1/4,1]$, and one has
\[
\phi(s)+\sum_{k=1}^{\infty}\psi(2^{-k}s)=1,\qquad\quad \forall s\in[0,\infty[%
\]
Set $\psi_{0}(s)=\phi(s)$ and 
$\psi_{k}(s)=\psi(2^{-k}s)$ for $ k\geq1$.
We then define the spectral localization operators $\psi_{k}(D),$
$k\in\mathbb{N},$ in the following way: for $u=\sum_{j}a_{j}e_{j},$ we define
\begin{equation}\label{1.33}
\psi_{k}(D)u=\sum_{j}\psi_k(\omega_{j})a_{j}e_{j} 
\end{equation}
One has $\sum_k\psi_{k}(D)=Id$ and $\psi_{i}(D)\psi_{j}(D)=0$ for $\vert i-j\vert \geq 2$. In addition,
we introduce
\begin{equation}\label{1.34}
S_{k}(D)=\sum_{j=0}^{k}\psi_j(D)=\psi_0(2^{-k}D),\qquad\quad
k\geq0 
\end{equation}
Obviously, the operators $\psi_{k}(D)$ and $S_{k}(D)$ acts as bounded operators on $H=H^1_0\times L^2$. 
The spectral localization result of \cite{DehmanLebeau07}  reads as follows.

\begin{thm}\label{thm2}Assume that the geometric control condition {\bf GCC}
holds true, and that the control function $\chi(t,x)$ is smooth.
There exists $C>0$  such that for every $k\in\mathbb{N}$,    the following inequality holds true
\be\label{1.35}\ba{rcl}
&\Vert \psi_{k}(D)\Lambda-\Lambda\psi_{k}(D)\Vert_{H}\leq C2^{-k}\\
&\Vert S_{k}(D)\Lambda-\Lambda S_{k}(D)\Vert_{H}\leq C2^{-k}\\
\ea\ee
\end{thm}

Theorem \ref{thm2} states that the optimal control operator $\Lambda$,
up to  lower order terms,  acts
individually on each frequency block of the solution.
For instance, if $e_{n}$ is the $n$-th eigenvector of the orthonormal
basis of $L^{2}(\Omega),$ if one drives the data $(0,0)$ to $(e_{n}
,0)$ in (\ref{1.16bis}) using the optimal control, both the solution $u$ and control $v$ in equation (\ref{1.16bis}) will essentially live at frequency $\omega_{n}$ for
$n$ large. We shall do experiments on this fact in section \ref{sec4.1}, and we will clearly see the impact of the regularity
of the control function  $\chi(t,x)$ on the accuracy of the frequency localization 
of the numerical control.\\

Since by the above results the optimal control operator $\Lambda$ preserves the regularity and
the frequency localization,  it is very natural to expect that $\Lambda$ is in fact a micro-local
operator, and in particular preserves the wave front set. For an introduction to micro-local
analysis and  pseudo-differential calculus, we refer to \cite{T} and \cite{Hormander85}.
In \cite{DehmanLebeau07}, it is proved that the optimal control operator $\Lambda$ for waves
on a compact Riemannian manifold without boundary is in fact an elliptic $2\times 2$ matrix
of pseudo-differential operators.
This is quite an easy result, since if $\chi(t,x)$ is smooth the Egorov theorem
implies that the operator $M_T$ given by (\ref{1.15bis}) is a $2\times 2$ matrix of pseudo-differential operators.
Moreover,  the geometric control condition {\bf GCC} implies easily that $M_T$ is elliptic.
Since $M_T$ is self adjoint, the fact that $M_T$ is an isomorphism follows then from
$Ker(M_T)=\{0\}$, which is equivalent to the injectivity of $\mathcal B^*$. This is proved
in \cite{BLR92} as a consequence of the uniqueness theorem of Calderon for the elliptic
second order operator $\triangle$. Then it follows that its inverse $\Lambda=M_T^{-1}$
is an elliptic pseudo-differential operator.\\

\begin{figure}
\begin{center}
\includegraphics[width=\textwidth]{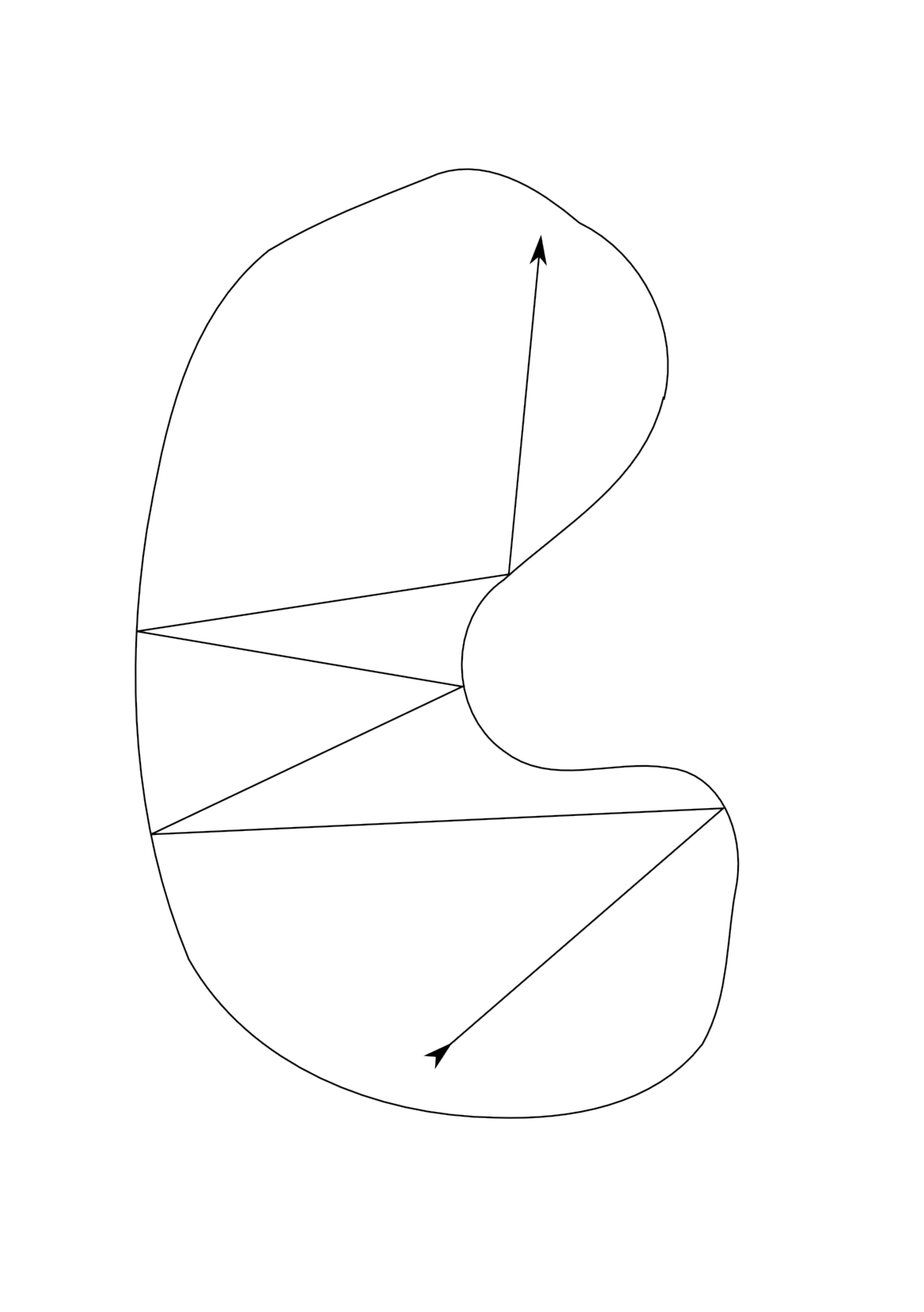}
\caption{\label{fig:rayon2} 
Example of optical ray with only transversal reflection points.}
\end{center}
\end{figure}

In our context, for waves in a 
bounded regular open subset $\Omega$ of
$\mathbb{R}^{d}$ with boundary Dirichlet condition, the situation is far much complicated, since there is
no Egorov theorem in the geometric setting of a manifold with boundary.
In fact, the Melrose-Sj\" ostrand theorem \cite{MS1}, \cite{MS2}
on propagation of singularities at the boundary 
(see also \cite{Hormander85}, vol 3, chapter XXIV for a proof) implies that the operator $M_T$
given by (\ref{1.15bis}) is a microlocal operator, but this is not sufficient to imply
that its inverse $\Lambda$ is micro-local. However let 
$\rho_0=(x_0,\xi_0)\in T^*\Omega, \vert\xi_0\vert=1$ be a point
in the cotangent space such that the two optical rays defined by the Melrose-Sj\" ostrand flow 
(see \cite{MS1}, \cite{MS2})
 $s\in [0,T]\rightarrow \gamma_{\pm\rho_0}(s)=
(x(\pm\rho_0,s),\xi(\pm\rho_0,s))$, with $\pm\rho_0=(x_0,\pm\xi_0)$, 
starting at  $x_0$ in the directions $\pm\xi_0$ have only transversal intersections with the boundary. 
Then it is not hard to show using formula (\ref{1.15bis}) and the parametrix of the wave operator, 
inside $\Omega$ and near transversal reflection points at the boundary $\partial\Omega$, as presented in figure \ref{fig:rayon2}, 
that $M_T$ is microlocally at $\rho_0$
an elliptic $2\times 2$ matrix of elliptic pseudo-differential operators. 
More precisely, let
$J$ be the isomorphism from $H^1_0\oplus L^2$ on $L^2\oplus L^2$ given by 
\be\label{t2}
J={1\over 2}\left(
\begin{array}
[c]{cc}%
\lambda & -i\\
\lambda & i%
\end{array}
\right)  
\ee 
One has $2\Vert J\underline u \Vert^2_{L^2\oplus L^2}=\Vert \underline u\Vert^2_{H^1_0\oplus L^2}$,
and if $u(t,x)$ is the solution of the wave operator $\Box u=0$ with Dirichlet boundary conditions on 
$\partial\Omega$, and Cauchy data at time $t_0$ equal to $(u_0,u_1)$, then one has
\be\label{t2bis}\ba{rcl}
\lambda u(t,.) =&\lambda\cos((t-t_0)\lambda)u_0+\sin((t-t_0)\lambda)u_1 \\
=& e^{i(t-t_0)\lambda}\Big({\lambda u_0 -iu_1\over 2}\Big)+e^{-i(t-t_0)\lambda}\Big({\lambda u_0 +iu_1\over 2}\Big)
\ea\ee
so the effect of the isomorphism $J$ is to split the solution $u(t,x)$ into a sum of two waves
with positive and negative temporal frequency. Moreover, one has
\be\label{t2ter}
Je^{itA}J^{-1}=\left(
\begin{array}
[c]{cc}%
e^{it\lambda} & 0\\
0 & e^{-it\lambda}%
\end{array}
\right)  
\ee 
Then $JM_TJ^{-1}$ acts as a non negative self-adjoint operator on $L^2\oplus L^2$, and is equal to
\be\label{t3}
\ba{rcl}
JM_TJ^{-1} =& {1\over 2}\left(
\begin{array}
[c]{cc}%
Q_+ & -\mathcal T\\
-\mathcal T^* & Q_-%
\end{array}
\right) \\  
Q_\pm =& \int_{0}^{T}e^{\pm is\lambda}\chi^{2}(T-s,.)e^{\mp is\lambda}ds\\
\mathcal T =& \int_{0}^{T}e^{is\lambda}\chi^{2}(T-s,.)e^{is\lambda}ds
\ea\ee
From (\ref{t3}), using the parametrix of the wave operator, 
inside $\Omega$ and near transversal reflection points at the boundary $\partial\Omega$,
and integration by parts to show that $\mathcal{T}$ is a smoothing operator, 
it is not difficult  to get that $JM_TJ^{-1}$ is microlocally at $\rho_0 \in T^*\Omega$ 
a pseudo-differential operator of order
zero with principal symbol
\be\label{t4}
\ba{rcl}
& \sigma_0(JM_TJ^{-1})(\rho_0)={1\over 2}\left(
\begin{array}
[c]{cc}%
q_+(x_0,\xi_0) & 0\\
0 & q_-(x_0,\xi_0)%
\end{array}
\right) \\  
&q_\pm(x_0,\xi_0)=\int_{0}^{T}\chi^{2}(T-s,x(\pm\rho_0,s))ds
\ea\ee
Obviously, condition (GCC) guarantees that $\sigma_0(JM_TJ^{-1})(\rho_0)$ is elliptic, and therefore
$J\Lambda J^{-1}$ will be at $\rho_0$
a pseudo-differential operator of order
zero with principal symbol
\be\label{t5}
\sigma_0(J\Lambda J^{-1})(\rho_0)=2\left(
\begin{array}
[c]{cc}%
q_+^{-1}(x_0,\xi_0) & 0\\
0 & q_-^{-1}(x_0,\xi_0)%
\end{array}
\right)
\ee

Therefore, the only difficulty  in order to prove that $\Lambda$ is a microlocal operator is
to get a precise analysis of the structure of the operator $M_T$ near rays which are tangent 
to the boundary. Since the set of $\rho\in ^bS^*\overline\Omega$ for which the optical ray
$\gamma_\rho(s)$ has only transversal points of intersection with the boundary is dense
in $T^*\overline\Omega\setminus T^*_{\partial\Omega}$ (see \cite{Hormander85}), it is not surprising
that our numerical experiments in section \ref{sec4.2} (where we compute the optimal control associated to a Dirac
mass $\delta_{x_0}, x_0\in \Omega$), confirms the following conjecture: 
\begin{conj} \label{conj1}
Assume that the geometric control condition {\bf GCC}
holds true, that the control function $\chi(t,x)$ is smooth, and that the optical rays have no 
infinite order of contact with the boundary. Then $\Lambda$  is a microlocal operator.
\end{conj}

Of course, part of the difficulty is to define correctly what is a microlocal operator
in our context. In the above conjecture, microlocal will implies in particular that the
optimal control operator $\Lambda$ preserves the  wave front set. A far less precise
information is to know that $\Lambda$ is a microlocal operator at the level of microlocal defect measures,
for which we refer to  \cite{G1}. But this is an easy by-product of the result of 
N. Burq and G. Lebeau in   \cite{BL01}.

 %%%%%%%%%%%%%%%%%%%%%%%%%%%%%%%
 
\subsection{The spectral Galerkin method}\label{sec2.3}
In this section we describe our numerical approximation of the optimal control operator
$\Lambda$, and we give some theoretical results on the numerical approximation $M_{T,\omega}$ 
of the operator $M_T$ given by (\ref{1.15bis}), even in the case where the 
geometric control condition {\bf GCC} is not satisfied.\\

For any cutoff frequency $\omega$, we denote by $\Pi_\omega$ the orthogonal projection, in the Hilbert space
$L^2(\Omega)$, on the finite dimensional linear subspace $L^2_\omega$ spanned by the eigenvectors $e_j$ for
$\omega_j \leq \omega$. By the Weyl formula, if $c_d$ denotes the volume of the unit ball in $\mathbb R^d$,
one has 
\be\label{t10}
N(\omega)=\text{dim}(L^2_\omega)\simeq (2\pi)^{-d}\text{Vol}(\Omega)c_d \omega^d \quad (\omega\rightarrow +\infty)
\ee

Obviously, $\Pi_\omega$ acts on $H=H^1_0\times L^2$ and commutes with $e^{itA}$, $\lambda$ and $J$. We define the Galerkin approximation $M_{T,\omega}$ of the operator $M_T$ as the operator on $L^2_\omega\times L^2_\omega$
\be\label{t11}
 M_{T,\omega}=\Pi_\omega M_T \Pi_\omega= \int_0^T e^{i(T-t)A} \Pi_\omega\begin{pmatrix}
0 & 0 \\
0 & \chi^2(t,.)
\end{pmatrix} \Pi_\omega
e^{-i(T-t)A} dt 
\ee 
Obviously, the matrix $M_{T,\omega}$ is symmetric and non negative for the Hilbert structure induced
by $H$ on $L^2_\omega\times L^2_\omega$, and by (\ref{1.14}) one has with $n_\omega(t)=B^*(t)\Pi_\omega e^{-i(T-t)A}$
\be\label{t13}
M_{T,\omega}= \int_0^T n_\omega^*(t)n_\omega(t)dt\\
\ee 
By (\ref{t3}) one has also

\be\label{t12}
\ba{rcl}
JM_{T,\omega}J^{-1} =& {1\over 2}\left(
\begin{array}
[c]{cc}%
Q_{+,\omega} & -\mathcal T_\omega\\
-\mathcal T^*_\omega & Q_{-,\omega}%
\end{array}
\right) \\  
Q_{\pm,\omega} =& \int_{0}^{T}e^{\pm is\lambda}\Pi_\omega\chi^{2}(T-s,.)\Pi_\omega e^{\mp is\lambda}ds\\
\mathcal T_\omega =& \int_{0}^{T}e^{is\lambda}\Pi_\omega\chi^{2}(T-s,.)\Pi_\omega e^{is\lambda}ds
\ea\ee

Let us first recall two easy results. For convenience, we recall here the proof of these results. The first result states that the matrix $M_{T,\omega}$
is always invertible.

\begin{lem}\label{lem1} For any (non zero) control function $\chi(t,x)$, the matrix $M_{T,\omega}$ is invertible.
\end{lem}
\bproof 
Let $u=(u_0,u_1)\in L^2_\omega\times L^2_\omega$ such that $M_{T,\omega}(u)=0$. By (\ref{t13}) one has

$$
0=(M_{T,\omega}(u)\vert u)_H=\int_0^T\Vert n_\omega(t)(u)\Vert^2_H dt.
$$
 This implies  
$n_\omega(t)(u)=0$ for almost all $t\in ]0,T[$. If $u(t,x)$ is the solution of the wave equation with
Cauchy data $(u_0,u_1)$ at time $T$, we thus get by (\ref{1.22}) and (\ref{1b}) 
$\psi(t)\chi_0(x)\partial_t u(t,x)=0$ for $t\in [0,T]$, and since $\psi(t)>0$ on $]0,T[$
and $\chi_0(x)>0$ on $U$, we get $\partial_t u(t,x)=0$ on $]0,T[\times U$. One has $u_0=\sum_{\omega_j\leq\omega}a_je_j(x)$,
$u_1=\sum_{\omega_j\leq\omega}b_je_j(x)$ and 

$$
\partial_t u(t,x)=\sum_{\omega_j\leq\omega}\omega_ja_j\sin((T-t)\omega_j)e_j(x)+\sum_{\omega_j\leq\omega}b_j\cos((T-t)\omega_j)e_j(x)
$$

Thus we get $\sum_{\omega_j\leq\omega}\omega_ja_je_j(x)=\sum_{\omega_j\leq\omega}b_je_j(x)=0$ for $x\in U$,
which implies, since the eigenfunctions $e_{j}$ are analytic in $\Omega$, 
%easily by the Calderon uniqueness theorem
 that $a_j=b_j=0$ for all $j$.
\eproof

For any $\omega_0\leq\omega$, we define $\Pi_\omega^\perp=1-\Pi_\omega$, and we set
\be\label{t17}
\ba{rcl}
&\Vert \Pi_\omega^\perp\Lambda \Pi_{\omega_0}\Vert_H=r_\Lambda(\omega,\omega_0)\\
&\Vert \Pi_\omega^\perp M_T \Pi_{\omega_0}\Vert_H=r_M(\omega,\omega_0)
\ea\ee
Since the ranges of the operators $\Lambda \Pi_{\omega_0}$ and $ M_T \Pi_{\omega_0}$ are finite dimensional vector spaces, one has
for any $\omega_0$
\be\label{t17bis}
\ba{rcl}
&\lim_{\omega\rightarrow\infty}r_\Lambda(\omega,\omega_0)=0\\
&\lim_{\omega\rightarrow\infty}r_M(\omega,\omega_0)=0
\ea\ee
The second result states that when {\bf GCC} holds true, the inverse matrix $M_{T,\omega}^{-1}$
converges in the proper sense to the optimal control operator $\Lambda$ when the cutoff frequency 
$\omega$ goes to infinity.

\begin{lem}\label{lem2} Assume that the geometric condition {\bf GCC} holds true. There exists $c>0$ such that the following holds true:
for any given $f\in H$, let $g=\Lambda(f)$, 
$f_\omega=\Pi_\omega f$ and $g^\omega=M_{T,\omega}^{-1}(f_\omega)$. Then, one has
\be\label{t14}
\Vert g-g^\omega\Vert_H \leq c\Vert f-f_\omega\Vert_H+ \Vert \Lambda(f_\omega)-M_{T,\omega}^{-1}(f_\omega)\Vert_H
\ee
with
\be\label{t15}
\lim_{\omega\rightarrow\infty}\Vert \Lambda(f_\omega)-M_{T,\omega}^{-1}(f_\omega)\Vert_H=0
\ee
\end{lem}
\bproof Since {\bf GCC} holds true, there exists $C>0$ such that one has by (\ref{1.15bis})
$(M_Tu\vert u)_H\geq C\Vert u\Vert_H^2$
for all $u\in H$, hence $(M_{T,\omega}u\vert u)_H\geq C\Vert u\Vert_H^2$ for all $u\in L^2_\omega\times L^2_\omega$. Thus, with $c=C^{-1}$, one has $\Vert\Lambda\Vert_H\leq c$ and  $\Vert M_{T,\omega}^{-1}\Vert_H\leq c$
for all $\omega$. Since $g-g^\omega=\Lambda(f-f_\omega)+\Lambda(f_\omega)-M_{T,\omega}^{-1}f_\omega$, 
(\ref{t14}) holds true. Let us prove that (\ref{t15}) holds true. 
With $\Lambda_\omega=\Pi_\omega\Lambda\Pi_\omega$, one has 
\be\label{t16}
\Lambda(f_\omega)-M_{T,\omega}^{-1}(f_\omega)=\Pi_\omega^\perp\Lambda(f_\omega)+(\Lambda_\omega-M_{T,\omega}^{-1})f_\omega\ee
Set for $\omega_0\leq\omega$, $f_{\omega_0,\omega}=(\Pi_\omega-\Pi_{\omega_0})f$. Then one has 
\be\label{t18}
\ba{rcl}
\Vert \Pi_\omega^\perp\Lambda(f_\omega)\Vert_H & = \Vert \Pi_\omega^\perp\Lambda\Pi_{\omega_0}(f)+
\Pi_\omega^\perp\Lambda (f_{\omega_0,\omega})\Vert_H \\
&\leq r_\Lambda(\omega,\omega_0)\Vert f\Vert_H+c\Vert
f_{\omega_0,\omega}\Vert_H
\ea\ee
On the other hand, one has
\be\label{t19}
\ba{rcl}
 (\Lambda_\omega-M_{T,\omega}^{-1})f_\omega 
&= M_{T,\omega}^{-1}(\Pi_\omega M_{T}\Pi_\omega^2\Lambda-\Pi_\omega)f_\omega\\
&=M_{T,\omega}^{-1}(\Pi_\omega M_{T}\Pi_\omega-\Pi_\omega M_T)\Lambda f_\omega\\
&=-M_{T,\omega}^{-1}\Pi_\omega M_T\Pi_\omega^\perp\Lambda \Pi_\omega f
\ea\ee
From (\ref{t19})  we get
\be\label{t20}
\Vert (\Lambda_\omega-M_{T,\omega}^{-1})f_\omega\Vert_H \leq \\
c \Vert M_T \Vert (\Vert\Lambda\Vert \Vert f_{\omega_0,\omega}\Vert_H + r_\Lambda(\omega,\omega_0)\Vert f\Vert_H)
\ee
Thus, for all $\omega_0\leq\omega$,  we get from (\ref{t18}), (\ref{t20}), and (\ref{t16}) 
\be\label{t21}
\Vert \Lambda(f_\omega)-M_{T,\omega}^{-1}(f_\omega)\Vert_H
\leq (1+c\Vert M_T \Vert)\Big( r_\Lambda(\omega,\omega_0)\Vert f\Vert_H+
c\Vert f_{\omega_0,\omega}\Vert_H\Big)
\ee
and (\ref{t15}) follows from (\ref{t17bis}),  (\ref{t21}) and 
$\Vert f_{\omega_0,\omega}\Vert_H\leq\Vert \Pi_{\omega_0}^\perp f\Vert_H \rightarrow 0$ when
$\omega_0\rightarrow\infty$.
\eproof

We shall now discuss two important points linked to the previous lemmas.
The first point is about the growth of the function
\be\label{t22}
\omega \rightarrow \Vert M_{T,\omega}^{-1}\Vert_H
\ee
when $\omega\rightarrow \infty$. This function is bounded when the geometric control condition 
{\bf GCC} is satisfied. Let us recall some known results in the general case.  
For simplicity, we assume that $\partial\Omega$ is an analytic hyper-surface of $\mathbb R^d$.
We know from \cite{Leb92} that, for $T>T_{u}$, where $T_{u} = 2 \sup_{x\in\Omega} \textrm{dist}_{\Omega}(x,U)$ is the uniqueness time, there exists $A>0$ such that
\be\label{t23}
\lim\text{sup}_{\omega \rightarrow \infty} {\log \Vert M_{T,\omega}^{-1}\Vert_H \over \omega}\leq A
\ee
On the other hand, when   there exists $\rho_0\in T^*\Omega$ such that the optical ray $s\in [0,T]\rightarrow \gamma_{\rho_0}(s)$
has only transversal points of intersection with the boundary and is such that 
$x(\rho_0,s)\notin \overline U$ for all $s\in [0,T]$, then {\bf GCC} is not satisfied. Moreover 
it is proven in \cite{Leb92}, using an explicit construction of a wave concentrated near this optical ray, 
that there exists $B>0$ such that

\be\label{t23bis}
\lim\text{inf}_{\omega \rightarrow \infty} {\log \Vert M_{T,\omega}^{-1}\Vert_H \over \omega}\geq B
\ee

Our experiments lead us to think that the following conjecture may be true
for a ``generic" choice of the control function $\chi(t,x)$:

\begin{conj} \label{conj2} There exists $C(T,U)$ such that
\be\label{t24}
\lim_{\omega \rightarrow \infty} {\log \Vert M_{T,\omega}^{-1}\Vert_H \over \omega}=C(T,U)
\ee
\end{conj}

In our experiments, we have studied (see section \ref{sec4.5}) the behavior of $C(T,U)$ as a function of $T$,
when the geometric control condition 
{\bf GCC} is satisfied for the control domain $U$ for $T\geq T_0$.
These experiments confirm the conjecture \ref{conj2} when $T<T_0$. We have not seen any clear change in the behavior of
the constant $C(T,U)$ when $T$ is smaller than the uniqueness time $T_{u}$.\\

The second  point we shall discuss is the rate of convergence of our Galerkin approximation.
By lemma \ref{lem2}, and formulas (\ref{t20}) and (\ref{t21}), this speed of convergence is governed
by the function $r_\Lambda(\omega,\omega_0)$ defined in (\ref{t17}). The following lemma tells us
that when the control function is smooth, the convergence in (\ref{t17bis}) is very fast.

\begin{lem}\label{lem3}
Assume that the geometric condition {\bf GCC} holds true and that the control function $\chi(t,x)$
is smooth. Then there exists a function $g$ with rapid decay such that
\be\label{t24bis}
r_\Lambda(\omega,\omega_0)\leq g\left(\frac\omega{\omega_0}\right)
\ee
\end{lem}

\bproof
By theorem \ref{thm1}, the operator $\lambda^s \Lambda \lambda^{-s}$ is bounded on $H$ for all $s\geq0$. Thus we get, for all $s\geq 0$:
$$
\|\Pi^{\perp}_{\omega} \Lambda \Pi_{\omega_{0}}\|_{H} = \|\Pi^{\perp}_{\omega}\lambda^{-s}\lambda^{s} \Lambda \lambda^{-s}\lambda^{s}\Pi_{\omega_{0}}\|_{H} \leq C_{s} \left(\frac{\omega_{0}}\omega\right)^s
$$
where we have used $\|\lambda^{s}\Pi_{\omega_{0}}\|_{H} \leq \omega_{0}^s$ and $\|\Pi^{\perp}_{\omega}\lambda^{-s}\|_{H} \leq \omega^{-s}$.
The proof of lemma \ref{lem3} is complete.
\eproof

Let us recall that $JM_{T,\omega}J^{-1}$ and the operators $Q_{\pm,\omega}$ and 
$\mathcal T_\omega$ are defined by formula (\ref{t12}). For any bounded operator $M$ on $L^2$, the matrix
coefficients of $M$ in the basis of the eigenvectors $e_n$ are
\be\label{t25}
M_{i,j}=(Me_i\vert e_j)
\ee
From (\ref{1b}) and (\ref{t12}) one has, for $\omega_{i}\leq \omega$, $\omega_{j}\leq \omega$:
$$
\ba{rcl}
Q_{\pm,\omega,i,j} &=& \displaystyle \int_{0}^T \left(e^{\pm i s \lambda} \Pi_{\omega}\psi^2 (T-s) \chi_{0}^2 (x) \Pi_{\omega}e^{\mp i s \lambda}(e_{i}) | e_{j} \right)\, ds\\
&=& \displaystyle \int_{0} ^T \psi^2(T-s) e^{\pm i s (\omega_{j}-\omega_{i})}\, ds \, \left(\chi_{0}^2 e_{i}|e_{j}\right)
\ea
$$
Since $\psi(t)\in C_{0}^{\infty}$ has support in $[0,T]$, we get that, for any $k\in \mathbb{N}$,
there exists
 a constant $C_k$ independent of the cutoff frequency $\omega$, such that one has
\be\label{t29}
\sup_{i,j}\vert (\omega_i-\omega_j)^kQ_{\pm,\omega,i,j}\vert\leq C_k
\ee
Moreover, by the results of \cite{DehmanLebeau07}, we know that the operator $\mathcal T$ defined in (\ref{t3})
is smoothing, and 
therefore we get
\be\label{t31}
\sup_{i,j}\vert (\omega_i+\omega_j)^k\vert\mathcal T_{\omega,i,j}\leq C_k
\ee

\begin{figure}
\begin{center}
\includegraphics[width=\textwidth]{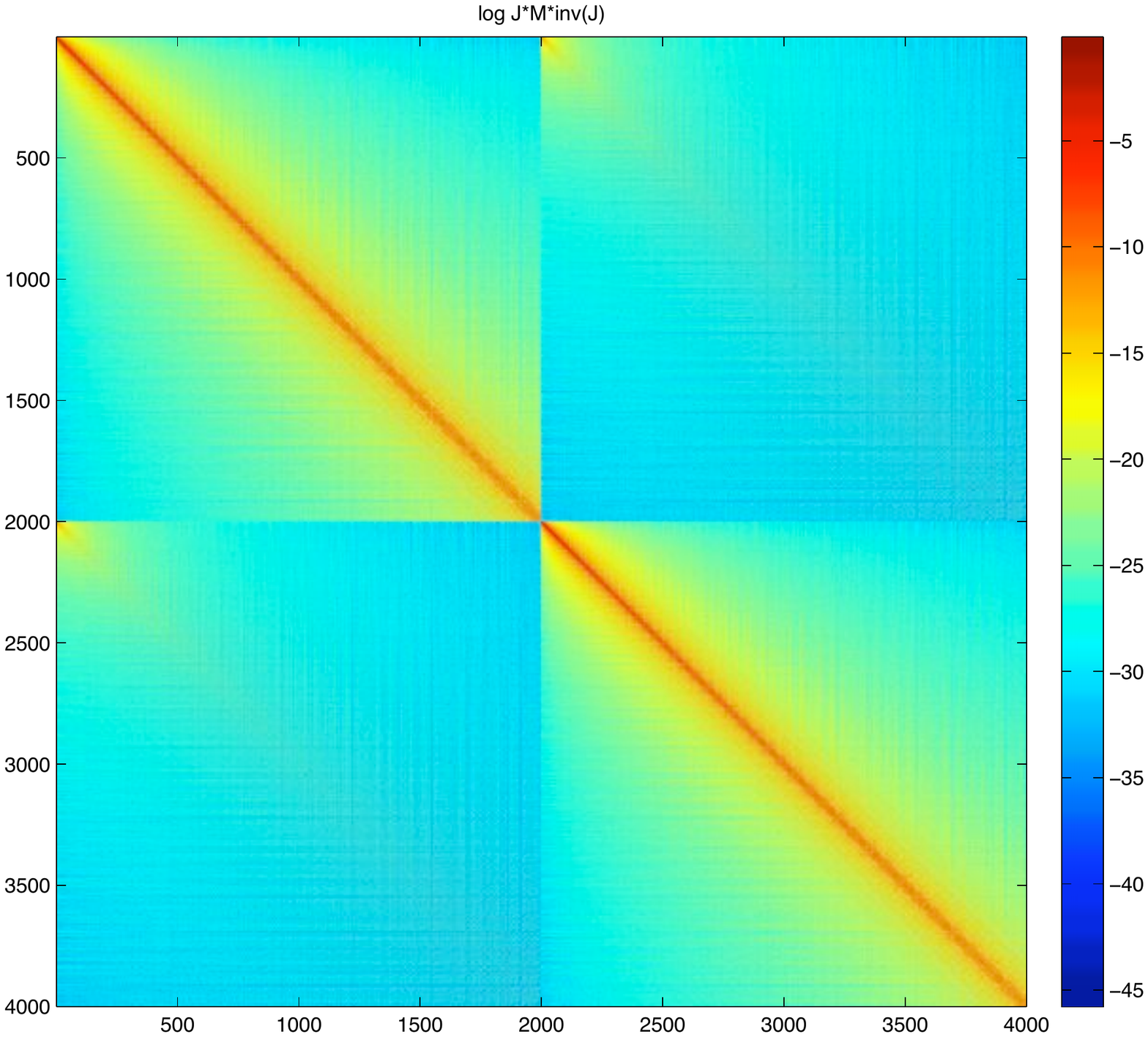}
\caption{\label{fig:MTsmooth} 
View of the logarithm of the coefficients of the matrix $JM_TJ^{-1}$, for the square geometry, with smooth control. This illustrates decay estimates (\ref{t29}) and (\ref{t31}).}
\end{center}
\end{figure}

\begin{figure}
\begin{center}
\includegraphics[width=\textwidth]{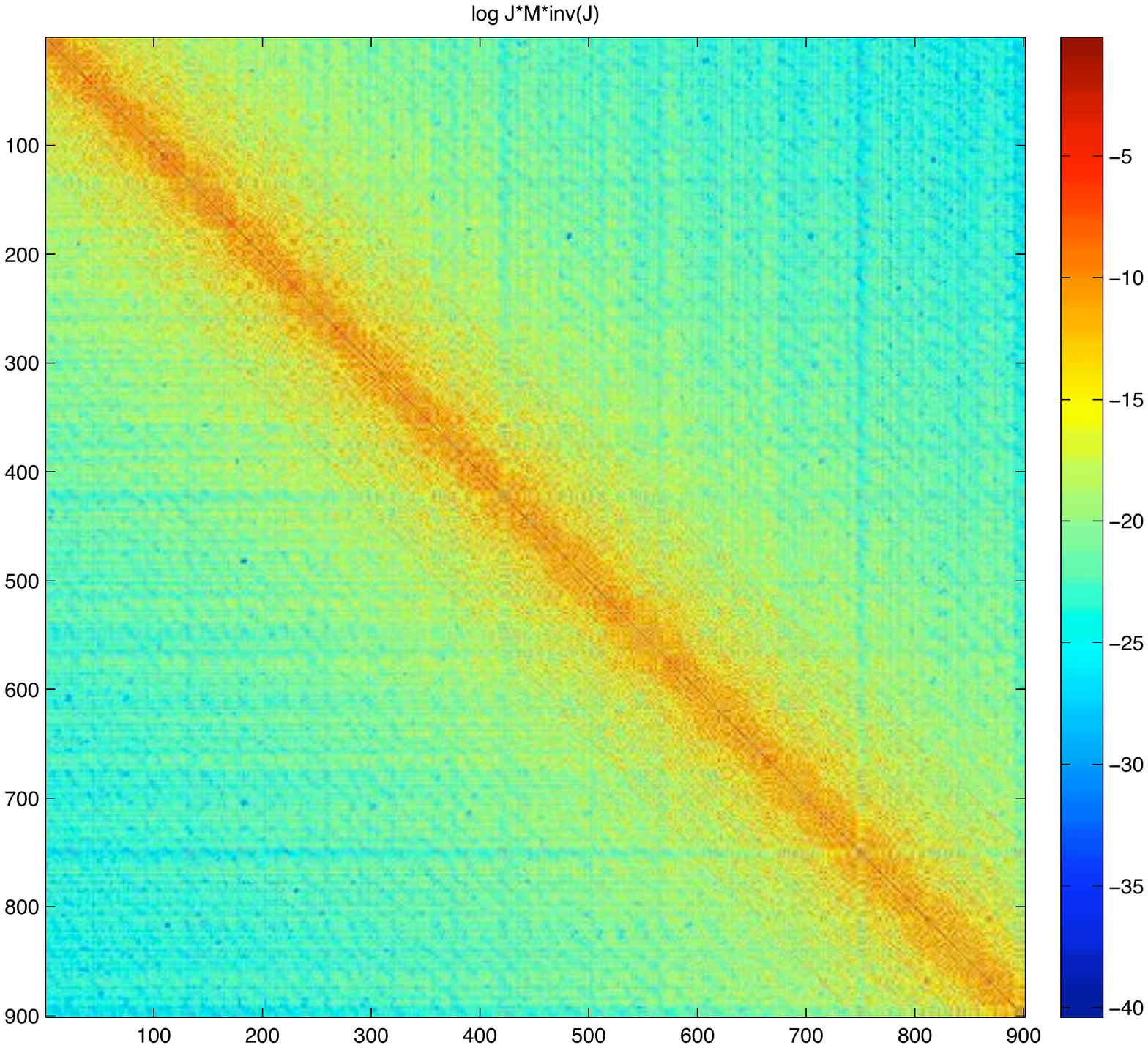}
\caption{\label{fig:MTsmoothzoom} 
View of the logarithm of the coefficients of the matrix $JM_TJ^{-1}$, for the square geometry, with smooth control (zoom). This illustrates decay estimate (\ref{t29}).}
\end{center}
\end{figure}

\begin{figure}
\begin{center}
\includegraphics[width=\textwidth]{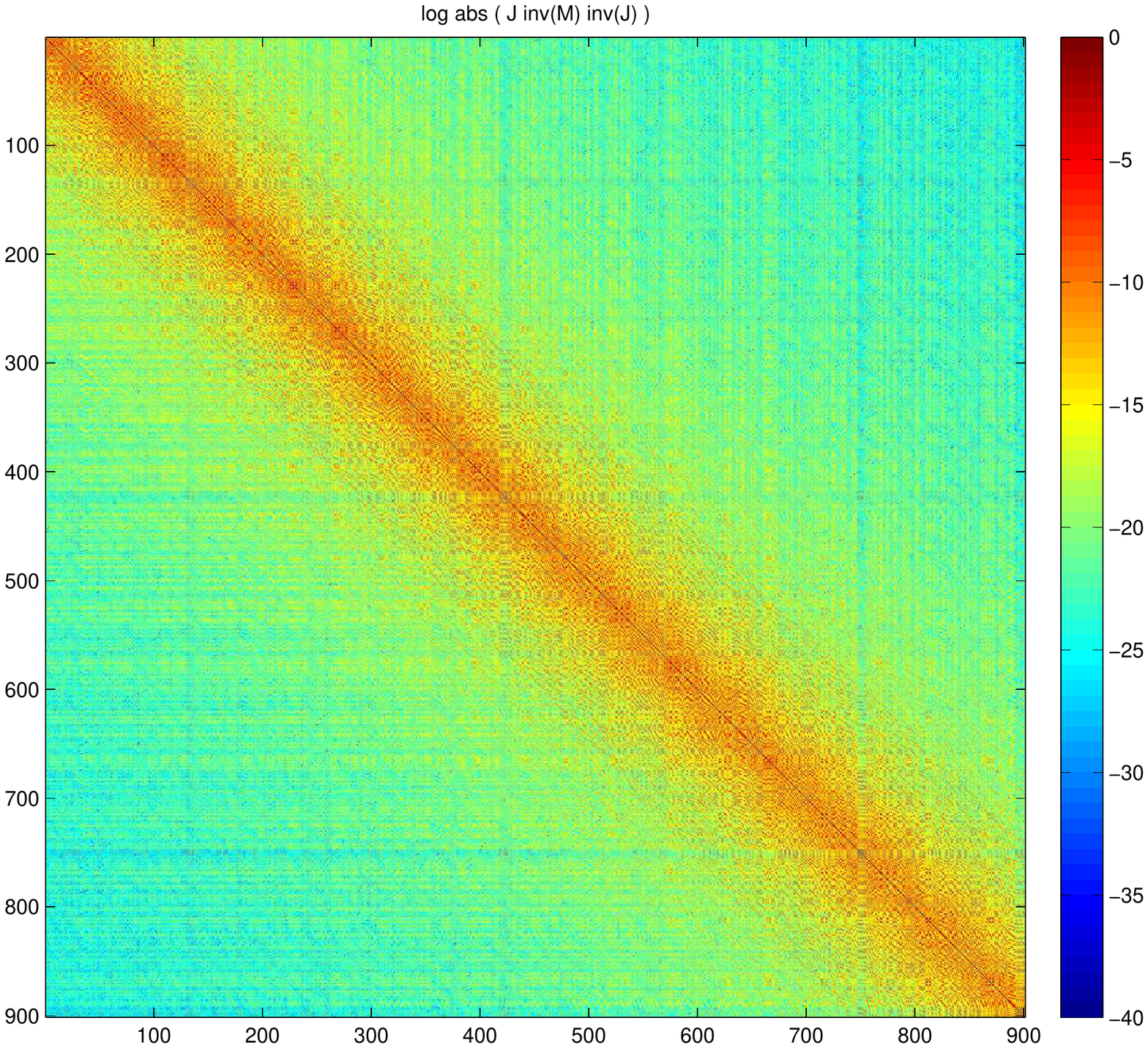}
\caption{\label{fig:invMTsmoothzoom} 
View of the logarithm of the coefficients of the matrix $JM_T^{-1}J^{-1}$, for the square geometry, with smooth control (zoom).}
\end{center}
\end{figure}

\begin{figure}
\begin{center}
\includegraphics[width=\textwidth]{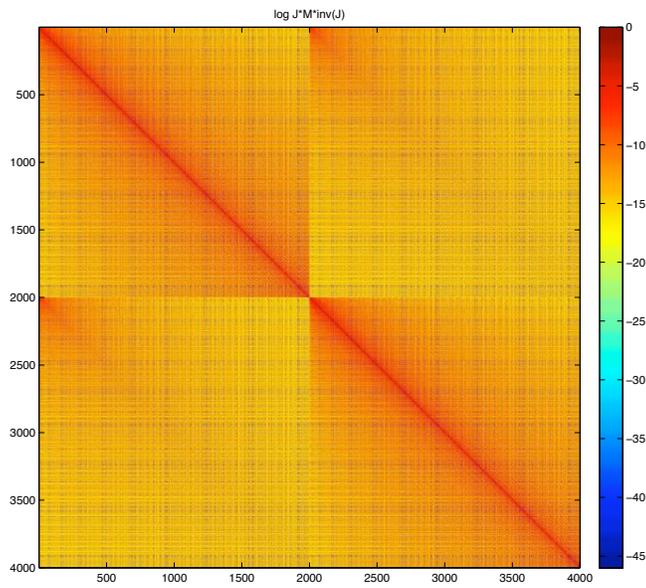}
\caption{\label{fig:MTconstant} 
View of the logarithm of the coefficients of the matrix $JM_TJ^{-1}$, for the square geometry, with non-smooth control. Note that the color scaling is the same as in Figure \ref{fig:MTsmooth}.}
\end{center}
\end{figure}

\begin{figure}
\begin{center}
\includegraphics[width=\textwidth]{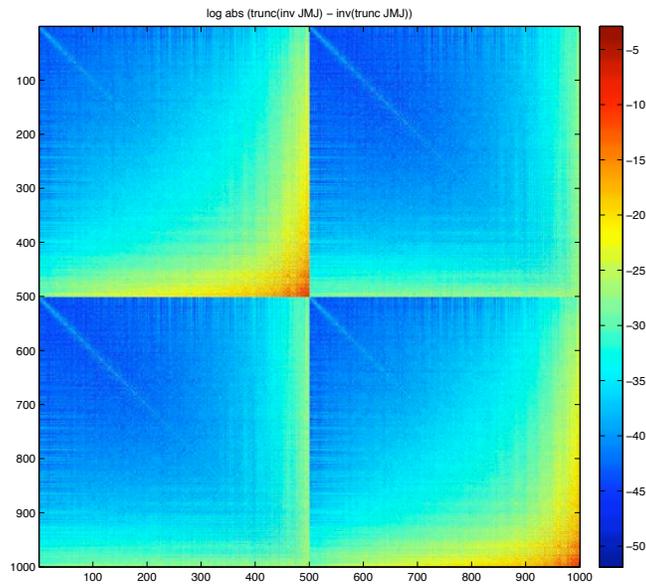}
\caption{\label{fig:MTconstant2} 
View of the logarithm of the coefficients of the matrix $J\left[((M_T)^{-1})_{\omega}-((M_T)_{\omega})^{-1}\right]J^{-1} = J\left[\Lambda_{\omega}-((M_T)_{\omega})^{-1}\right]J^{-1}$, for the square geometry, with smooth control. The $M_{T}$ matrix is computed with 2000 eigenvalues, the cutoff frequency $\omega$ being associated with the 500th eigenvalue.}
\end{center}
\end{figure}

Figure \ref{fig:MTsmooth} shows the log of $JM_{T,\omega}J^{-1}$ coefficients and illustrates the decay estimates (\ref{t29}) and (\ref{t31}). A zoom in is done in figure \ref{fig:MTsmoothzoom}, so that we can observe more precisely (\ref{t29}). In particular we can notice that the distribution of the coefficients along the diagonal of the matrix is not regular. Figure \ref{fig:invMTsmoothzoom} presents the same zoom for $JM_{T,\omega}^{-1}J^{-1}$. This gives an illustration of the matrix structure
of a microlocal operator. Figure \ref{fig:MTconstant} represents the log of $JM_{T,\omega}J^{-1}$ coefficients \emph{without smoothing}. And finally, figure \ref{fig:MTconstant2} gives a view of the convergence of our Galerkin approximation, as it presents the matrix entries of $J (\Lambda_{\omega} - M_{T,\omega}^{-1})J^{-1}$, illustrating lemma \ref{lem2}, its proof, and lemma \ref{lem3}.

%%%%%%%%%%%%%%%%%%%%%%%%%%%%%%%%%%%%%%%%%%%%%%%%%%
\subsection{Computation of the discrete control operator}\label{sec2.4}

For any real $\omega$, let $N(\omega)=\sup\{n,\omega_n\leq \omega\}$. Then the dimension of the vector space $L^2_\omega$
is equal to $N(\omega)$. Let us define the following $(\phi_j)_{1\leq j\leq 2N(\omega)}$:
\be\label{mneq:base}
\baa{rcl}
\phi_j={e_j\over \omega_j}&\textrm{for}&1\leq j\leq N(\omega)\\
\phi_j=e_{j-N(\omega)}&\textrm{for}&N(\omega)+1\leq j\leq 2N(\omega)
\eaa
\ee
Then $(\phi_j)_{1\leq j\leq 2N(\omega)}$ is an orthonormal basis  of the 
Hilbert space 
$H_\omega=\Pi_\omega(H^1_0(\Omega)\oplus L^2(\Omega))$. 

In this section we compute explicitly $\big(M_T \phi_l | \phi_k \big)_{H}$ for all $1\leq k,l \leq 2N(\omega)$. We recall 
\be
\label{mneq:4}
\ba{rcl}
e^{isA} \left[\ba{l} e_i \\0\ea\right] &=& \left[\ba{l} \cos(s\omega_i) e_i(x)\\-\omega_i \sin(s\omega_i) e_i(x)\ea\right] \medskip\\
e^{isA} \left[\ba{l} 0 \\e_i\ea\right] &=& \left[\ba{l} \sin(s\omega_i) e_i(x)/\omega_i\\ \cos(s\omega_i) e_i(x)\ea\right] 
\ea
\ee
We now compute the coefficients of the $M_T$ matrix, namely ${M_T}_{n,m} = \big(M_T \phi_n | \phi_m \big)_{H}$:
\be
\label{mneq:5}
\ba{rcl}
{M_T}_{n,m} &=& \big(M_T \phi_n | \phi_m \big)_{H}\\
&=& \int_0^T \big( e^{isA} B B^* e^{-isA} \phi_n | \phi_m \big)_{H}\, dt\\
&=& \int_0^T \big(  \left(\ba{ll} 0&0\\0&\chi^2\ea\right) e^{-isA} \phi_n | e^{-isA}\phi_m \big)_{H}\, dt\ea
\ee
We now have to distinguish four cases, depending on $m,n$ being smaller or larger than $N(\omega)$. 
For the case $(m,n)\leq N(\omega)$ we have:
\be
\label{mneq:5-1}
\ba{rcl}
{M_T}_{n,m}
 &=& \int_0^T \big(  \left(\ba{ll} 0&0\\0&\chi^2\ea\right) e^{-isA} \phi_n | e^{-isA}\phi_m \big)_{H}\, ds \medskip\\
&=& \int_0^T \big(  \left(\ba{ll} 0&0\\0&\chi^2\ea\right) \left[\ba{l} \cos(s\omega_n) f_n(x)\\\omega_n \sin(s\omega_n) f_n(x)\ea\right]  | \left[\ba{l} \cos(s\omega_m) f_m(x)\\\omega_m \sin(s\omega_m) f_m(x)\ea\right]   \big)_{H}\, ds \medskip\\
&=& \int_0^T \big(  \left[\ba{l} 0\\ \omega_n \chi^2\sin(s\omega_n) f_n(x)\ea\right]  | \left[\ba{l} \cos(s\omega_n) f_m(x)\\\omega_m \sin(s\omega_m) f_m(x)\ea\right]   \big)_{H}\, ds \medskip\\
&=& \int_0^T (  (\psi(t)\chi_0(x))^2 \omega_n \sin(s\omega_n) f_n(x)  | \omega_m \sin(s\omega_m) f_m(x)   )_{L^2(\Omega)}\, ds\medskip\\
&=& \int_0^T\psi^2 \sin(s\omega_n) \sin(s\omega_m) \,ds \int_\Omega    \chi_0^2 e_n(x)   e_m(x)   \, dx\medskip\\
&=& a_{n,m} G_{n,m}
\ea
\ee
where
\be
\label{mneq:6-1}
a_{n,m} = \int_0^T \psi^2 \sin(s\omega_m) \sin(s\omega_n) \,ds 
\ee
and
\be
\label{mneq:7}
 G_{n,m} = \int_\Omega    \chi_0^2(x)\, e_m(x) \,  e_n(x)   \, dx   
 \ee
Similarly, for the case $n>N(\omega)$, $m\leq N(\omega)$ we have:
\be
\label{mneq:5-2}
\ba{rcl}
{M_T}_{n,m}
 &=& \int_0^T \big(  \left(\ba{ll} 0&0\\0&\chi^2\ea\right) e^{-isA} \phi_n | e^{-isA}\phi_m \big)_{H}\, ds \medskip\\
&=& \int_0^T \big(  \left(\ba{ll} 0&0\\0&\chi^2\ea\right) \left[\ba{l} \sin(s\omega_n) f_n(x)/\omega_n\\ \cos(s\omega_n) f_n(x)\ea\right]  | \left[\ba{l} \cos(s\omega_m) f_m(x)\\\omega_m \sin(s\omega_m) f_m(x)\ea\right]   \big)_{H}\, ds \medskip\\
&=& \int_0^T \big(  \left[\ba{l} 0\\  \chi^2\cos(s\omega_n) f_n(x)\ea\right]  | \left[\ba{l} \cos(s\omega_n) f_m(x)\\\omega_m \sin(s\omega_m) f_m(x)\ea\right]   \big)_{H}\, ds \medskip\\
&=& \int_0^T (    \chi^2\cos(s\omega_n) f_n(x)  | \omega_m \sin(s\omega_m) f_m(x)   )_{L^2(\Omega)}\, ds\medskip\\
&=& \int_0^T \psi^2\cos(s\omega_n) \sin(s\omega_m) \,ds \int_\Omega    \chi_0^2 e_n(x)   e_m(x)   \, dx\medskip\\
&=&  b_{n,m} G_{n,m}
\ea
\ee
where
\be
\label{mneq:6-2}
b_{n,m} = \int_0^T \psi^2\cos(s\omega_n) \sin(s\omega_m) \,ds 
\ee
For $n\leq N(\omega)$ and $m>N(\omega)$ we get:
\be
\label{mneq:5-3}
\ba{rcl}
{M_T}_{n,m}
 &=&  c_{n,m} G_{n,m}
\ea
\ee
where
\be
\label{mneq:6-3}
c_{n,m} = b_{m,n} = \int_0^T \psi^2 \cos(s\omega_m) \sin(s\omega_n) \,ds 
\ee
And for $m,n > N(\omega)$:
\be
\label{mneq:5-4}
\ba{rcl}
{M_T}_{n,m}
 &=&  d_{n,m} G_{n,m}
\ea
\ee
where
\be
\label{mneq:6-4}
d_{n,m} = \int_0^T\psi^2 \cos(s\omega_m) \cos(s\omega_n) \,ds 
\ee
The above integrals have to be implemented carefully when $|\omega_n-\omega_m|$ is small, even when $\psi(t)=1$.

%%%%%%%%%%%%%%%%%%%%%%%%%%%%
\section{Numerical setup and validation}\label{sec3}

\subsection{Geometries and control domains}\label{sec3.1}

\begin{figure}
\begin{center}
\includegraphics[width=\textwidth]{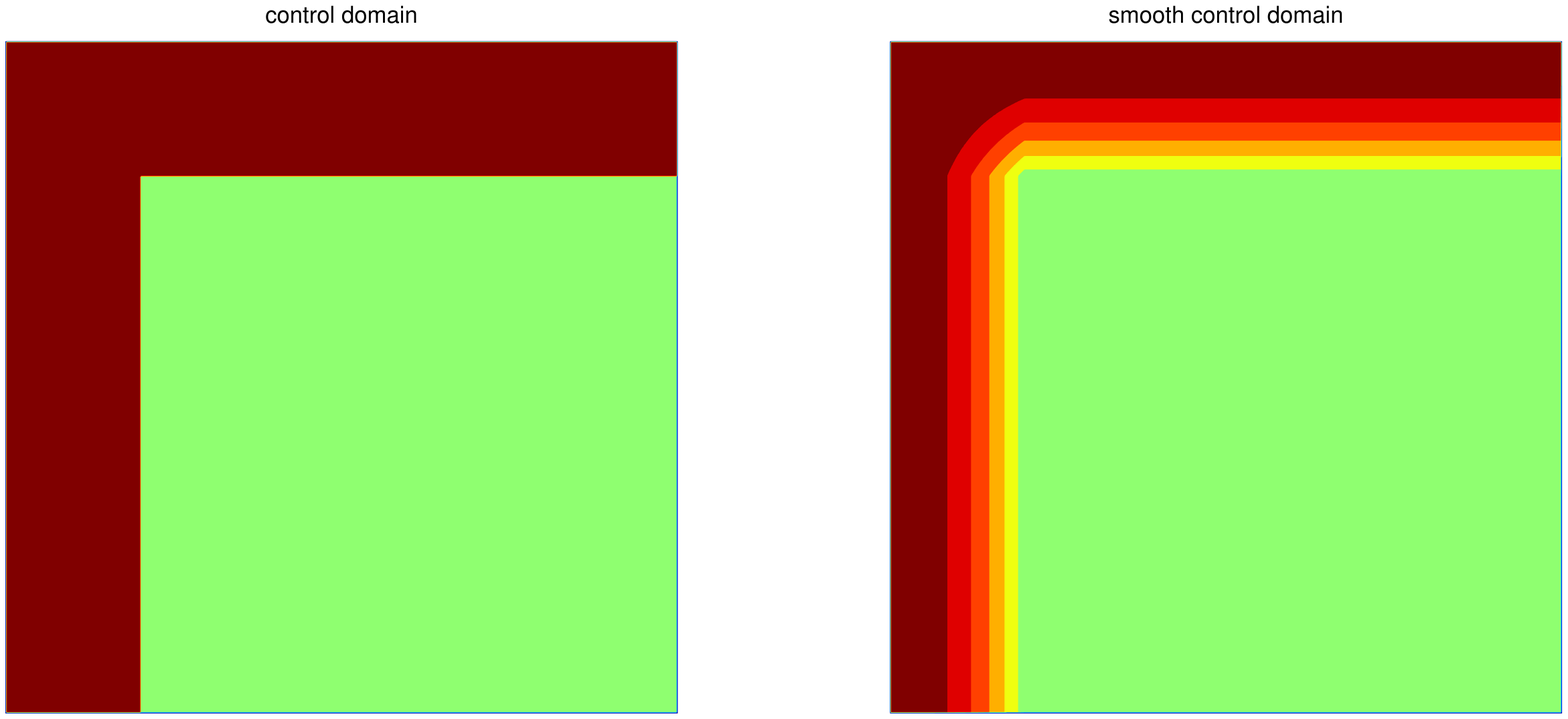}
\caption{\label{fig:domain:sq} 
Domain and example of a control domain for the square, with smoothing in space (right panel) or without (left panel).}
\end{center}
\end{figure}

\begin{figure}
\begin{center}
\includegraphics[width=\textwidth]{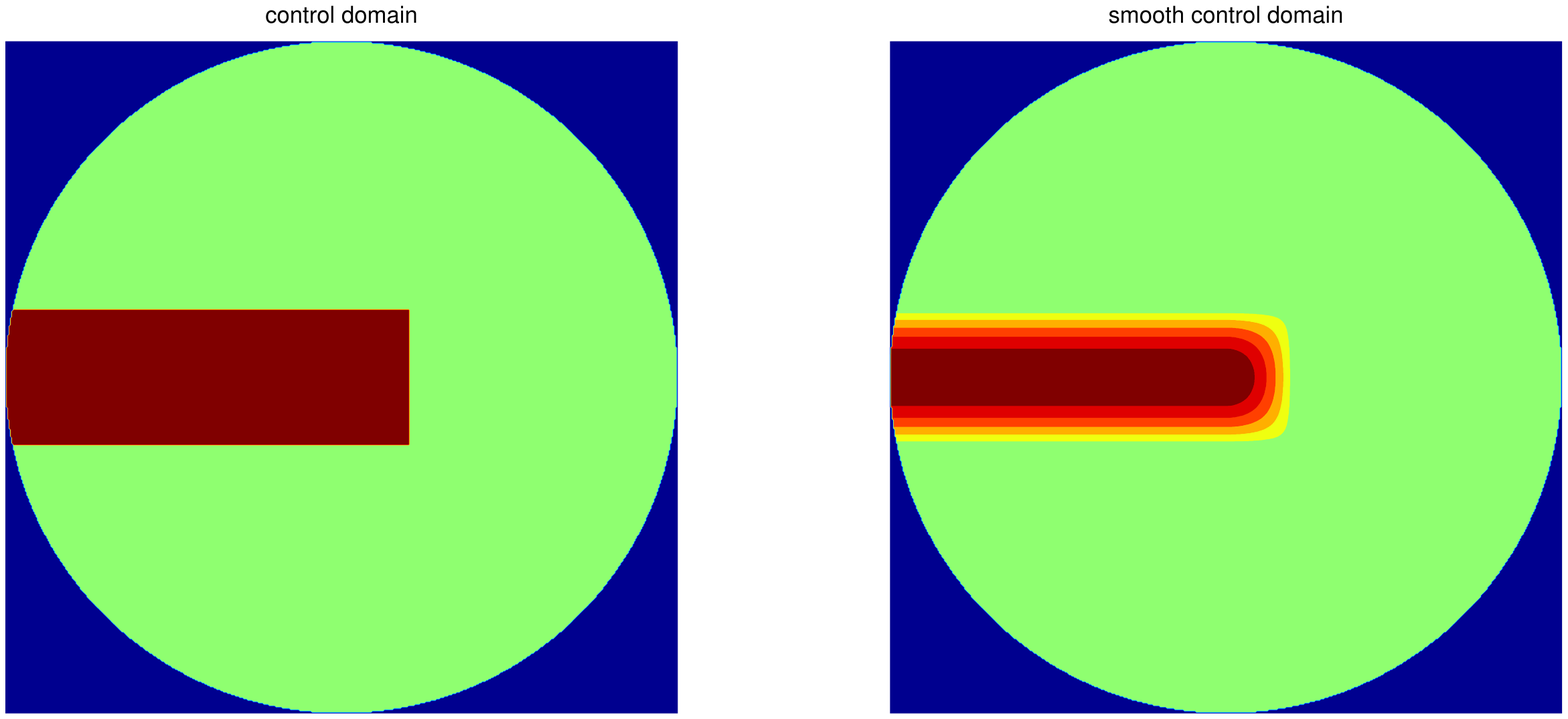}
\caption{\label{fig:domain:di} 
Domain and example of a control domain for the disc, with smoothing in space (right panel) or without (left panel).}
\end{center}
\end{figure}

\begin{figure}
\begin{center}
\includegraphics[width=\textwidth]{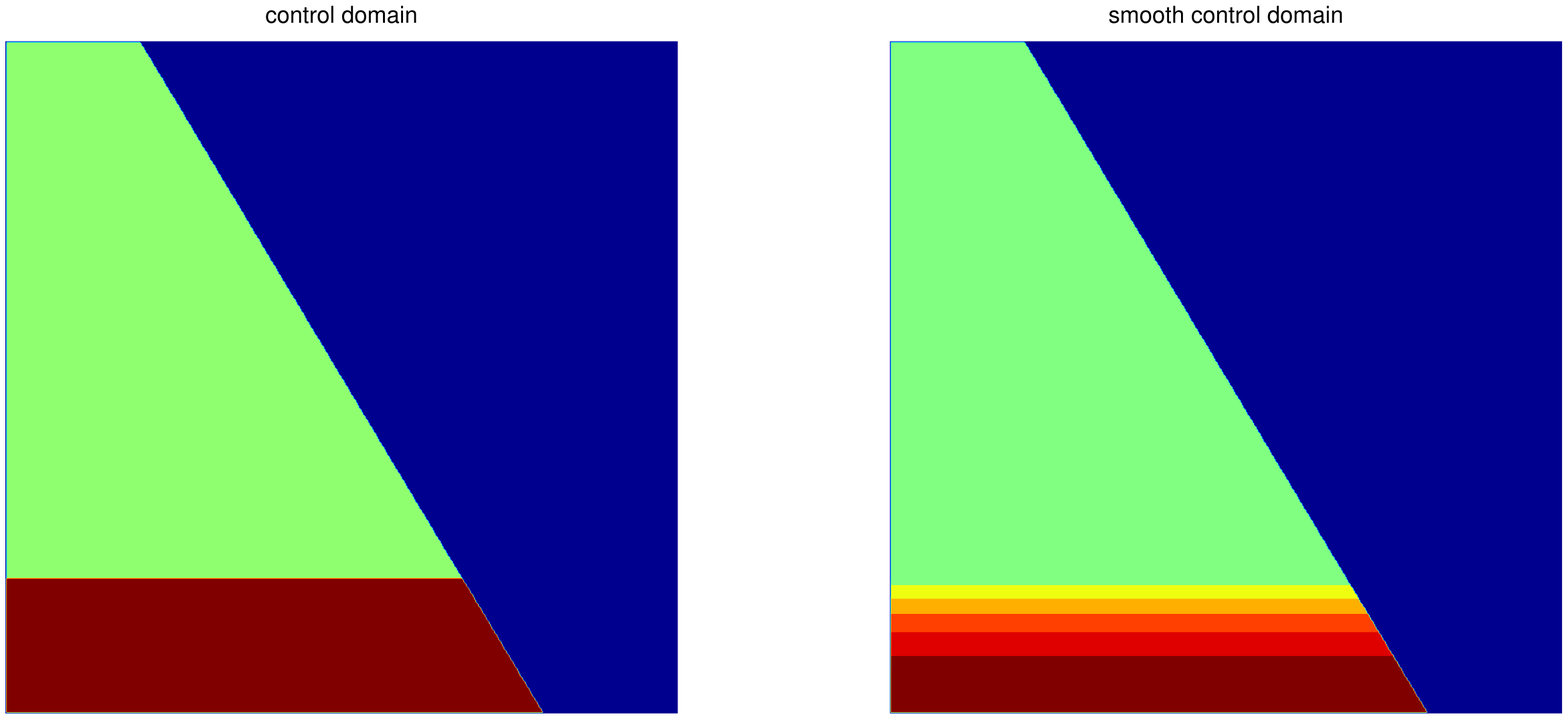}
\caption{\label{fig:domain:tr} 
Domain and example of a control domain for the trapezoid, with smoothing in space (right panel) or without (left panel).}
\end{center}
\end{figure}

The code we implemented allows us to choose the two-dimensional domain $\Omega$, as well as the control domain $ U $. In the sequel, we will present some results with three different geometries: square, disc and trapezoid. For each geometry, we have chosen a reference shape of control domain. It consists of the neighborhood of two adjacent sides of the boundary (in the square), of a radius (in the disc), of the base side (in the trapezoid). Then we adjust the width of the control domain, and also its smoothness (see next paragraph). Figures \ref{fig:domain:sq}, \ref{fig:domain:di} and \ref{fig:domain:tr} present these domains, and their respective control domains, either non-smooth (left panels) or smooth (right panels).

\subsection{Time and space smoothing}\label{sec3.2smoothing}

We will investigate the influence of the regularity of the function $\chi(t,x) = \psi(t) \chi_0(x)$. Different options have been set. 
\paragraph{Space-smoothing.} The integral (\ref{mneq:7}) defining $G_{n,m}$ features $\chi_0$. In the literature we find $\chi_0 = {\bf 1}_{U}$, so that
\be
\label{mneq:7bis}
G_{n,m} = \int_ U  e_n(x)\, e_m(x)\,dx
\ee
In \cite{DehmanLebeau07} the authors show that a smooth $\chi_0^2$ leads to a more regular control (see also theorem \ref{thm2} and lemma \ref{lem3}). Thus for each control domain $U$ we implemented both smooth and non-smooth (constant) cases.
The different implementations of $\chi_0$ are:
\begin{itemize}
\item constant case: $\chi_0(x,y) = {\bf 1}_U $,
\item ``smooth'' case: $\chi_0(x,y)$ has the same support of $U $, the width $a$ of the domain $\{ x\in \Omega, 0<\chi_0(x)<1 \}$ is adjustable, and on this domain $\chi$ is a polynomial of degree 2. For example, in the square we have:
\be
\ba{rcl}
\label{mneq:8} 
\chi_{0}(x,y) &=& {\bf 1}_{U } \Big[ 1- \Big({\bf 1}_{x\geq a} + \frac{x^2}{a^2}. {\bf 1}_{x< a}\Big) \\ 
&&\qquad \quad\Big({\bf 1}_{y\leq 1-a} + \frac{(1-y)^2}{a^2}. {\bf 1}_{y>1- a} \Big) \Big]
\ea
\ee
\end{itemize}

\paragraph{Time-smoothing.} Similarly, the time integrals (\ref{mneq:6-1},\ref{mneq:6-2},\ref{mneq:6-3},\ref{mneq:6-4}) defining $a$, $b$, $c$ and $d$ features $\psi(t)$, which is commonly chosen as ${\bf 1}_{[0,T]}$. As previously, better results are expected with a smooth $\psi(t)$. In the code, the integrals (\ref{mneq:6-1},\ref{mneq:6-2},\ref{mneq:6-3},\ref{mneq:6-4}) are computed explicitly, the different implementations of $\psi$ being:
\begin{itemize}
\item constant case $\psi = {\bf 1}_{[0,T]}$,
\item ``smooth case"
\be\label{mneq:8bis}
\psi(t) = \frac{4 t (T-t)}{T^2} {\bf 1}_{[0,T]}
\ee 
\end{itemize}

\subsection{Validation of the eigenvalues computation}\label{sec3.3}

\begin{figure}
\begin{center}
\includegraphics[width=\textwidth]{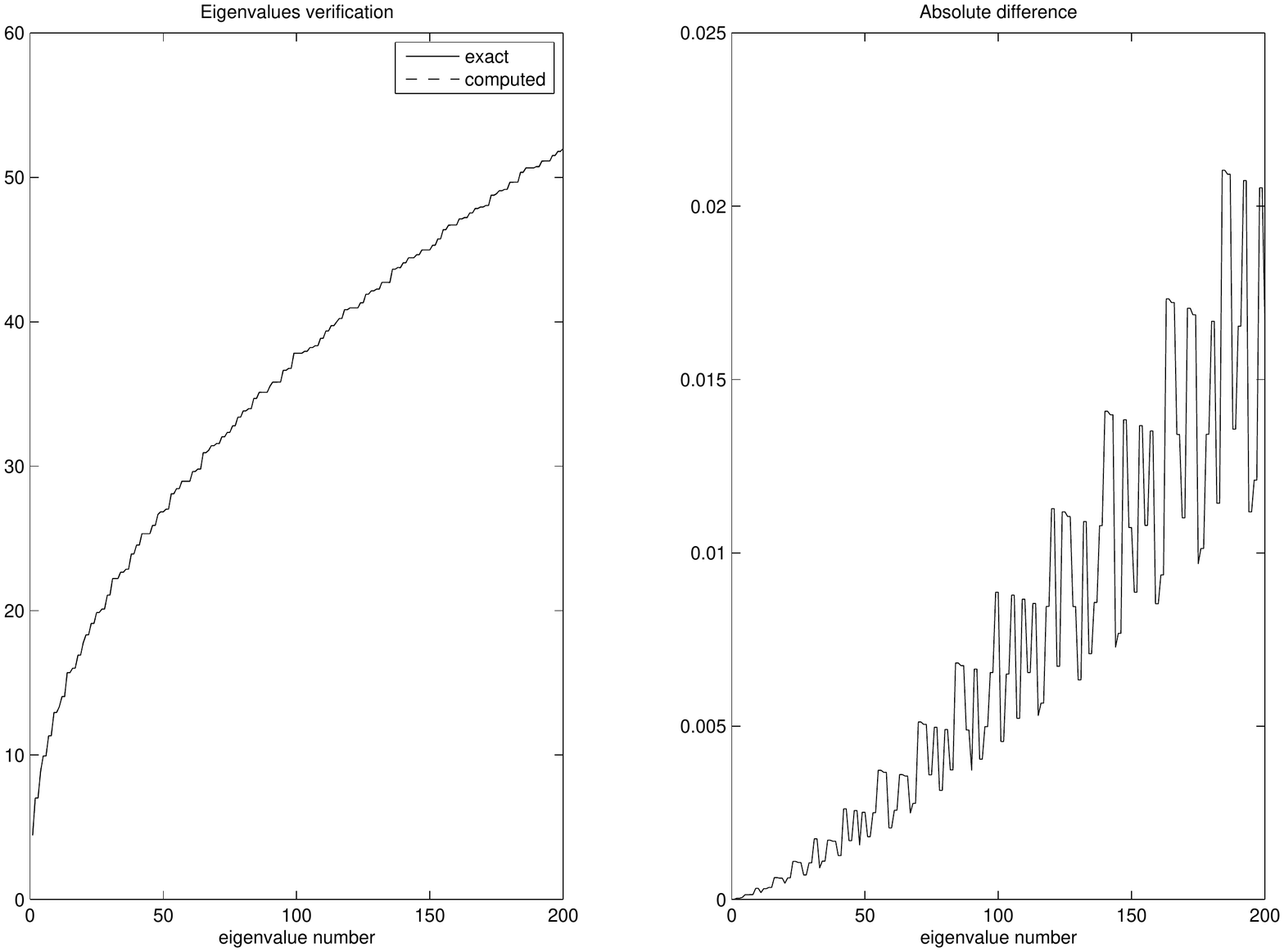}
\caption{\label{fig:eigs:sq} 
Verification of the eigenvalues computation in the square: exact and finite-differences-computed eigenvalues (left panel), and their absolute difference (right panel).}
\end{center}
\end{figure}

\begin{figure}
\begin{center}
\includegraphics[width=\textwidth]{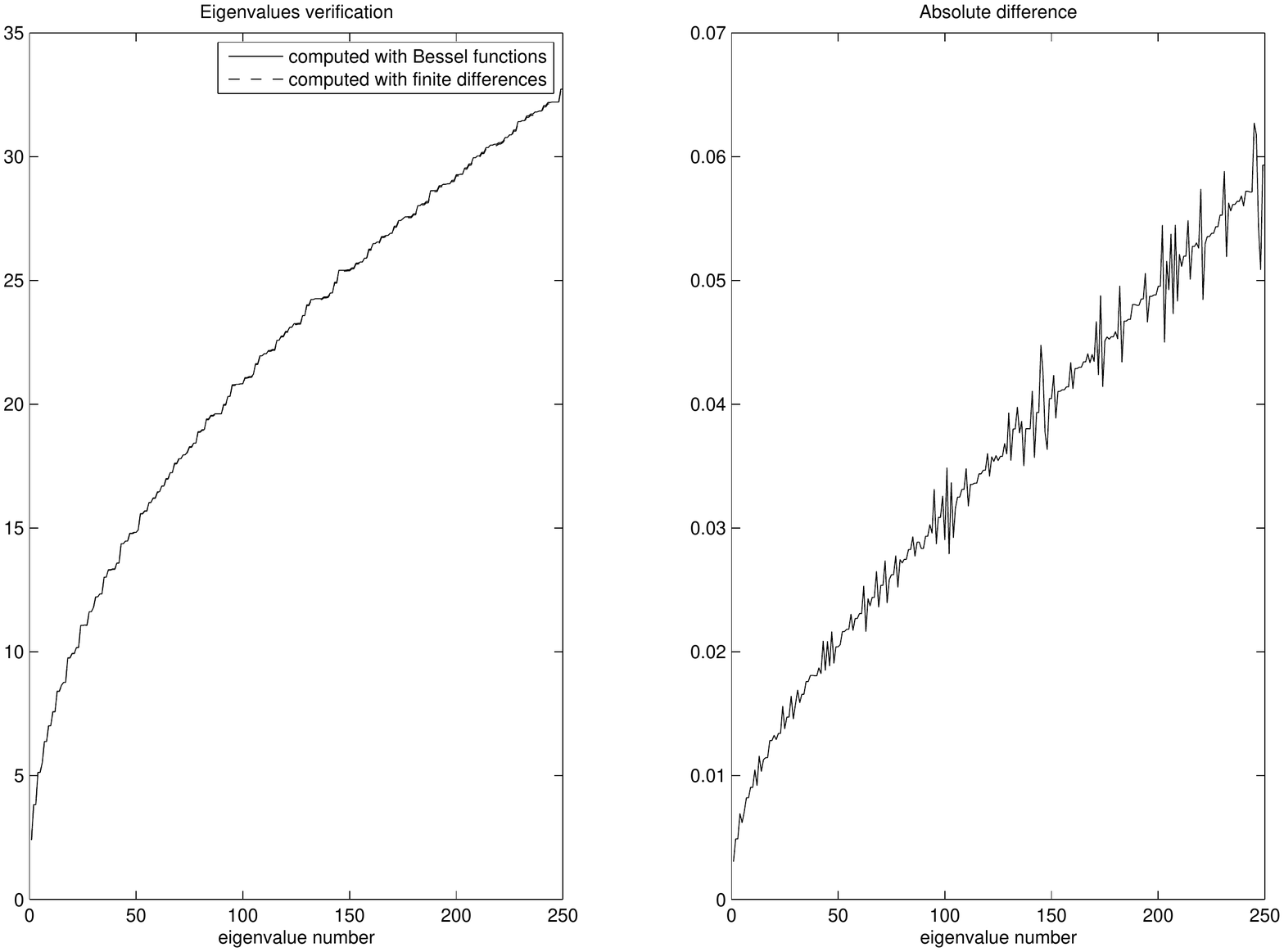}
\caption{\label{fig:eigs:di} 
Verification of the eigenvalues computation in the disc: eigenvalues computed either as zeros of Bessel functions or with finite differences (left panel), and their absolute difference (right panel).}
\end{center}
\end{figure}

The code we implemented has a wide range of geometries for $\Omega$. As it is a spectral-Galerkin method, it requires the accurate computation of eigenvalues and eigenvectors. We used Matlab \emph{eigs}\footnote{\label{matlab-eigs}\url{www.mathworks.com/access/helpdesk/help/techdoc/ref/eigs.html}} function. Figure \ref{fig:eigs:sq} shows the comparison between the first 200 exact eigenvalues in the square, and those computed by Matlab with $500\times 500$ grid-points. Figure \ref{fig:eigs:di} presents the same comparison in the disc, for 250 eigenvalues, the ``exact'' ones being computed as zeros of Bessel function.

\subsection{Reconstruction error. }\label{sec3.4recerror} 
In the sequel, we will denote the input data $u=(u_{0},u_{1})$, and its image by the control map  $w=(w_{0},w_{1}) = \Lambda(u_{0},u_{1})$, which will often be called the ``control". We recall from section \ref{sec2.1} that for a given data $u=(u_0,u_1)$ to be reconstructed at
time $T$, the optimal control $v(t)$ is given by 
\be
\label{mneq:9}
v(t)= \chi \partial_t e^{-i(T-t)A} w = \chi \partial_t e^{-i(T-t)A}\Lambda(u)
\ee  
Then, solving the wave equations (\ref{1.16bis}) forward, with null initial conditions and $\chi v$ as a forcing source, we reach $y=(y_0,y_1)$ in time $T$. Should the experiment be perfect, we would have $(y_0,y_1)=(u_0,u_1)$.  The reconstruction error is then by definition:
\be
\label{mneq:10}
E = \sqrt{\frac{\|u_0-y_0\|^2_{H^1(\Omega)} + \|u_1-y_1\|^2_{L^2(\Omega)}}{\|u_0\|^2_{H^1(\Omega)} + \|u_1\|^2_{L^2(\Omega)}}}
\ee

\subsection{Validation for the square geometry}\label{sec3.5}

\subsubsection{Finite differences versus exact eigenvalues}\label{sec3.5-1verif1}

\begin{figure}
\begin{center}
\includegraphics[width=1\textwidth]{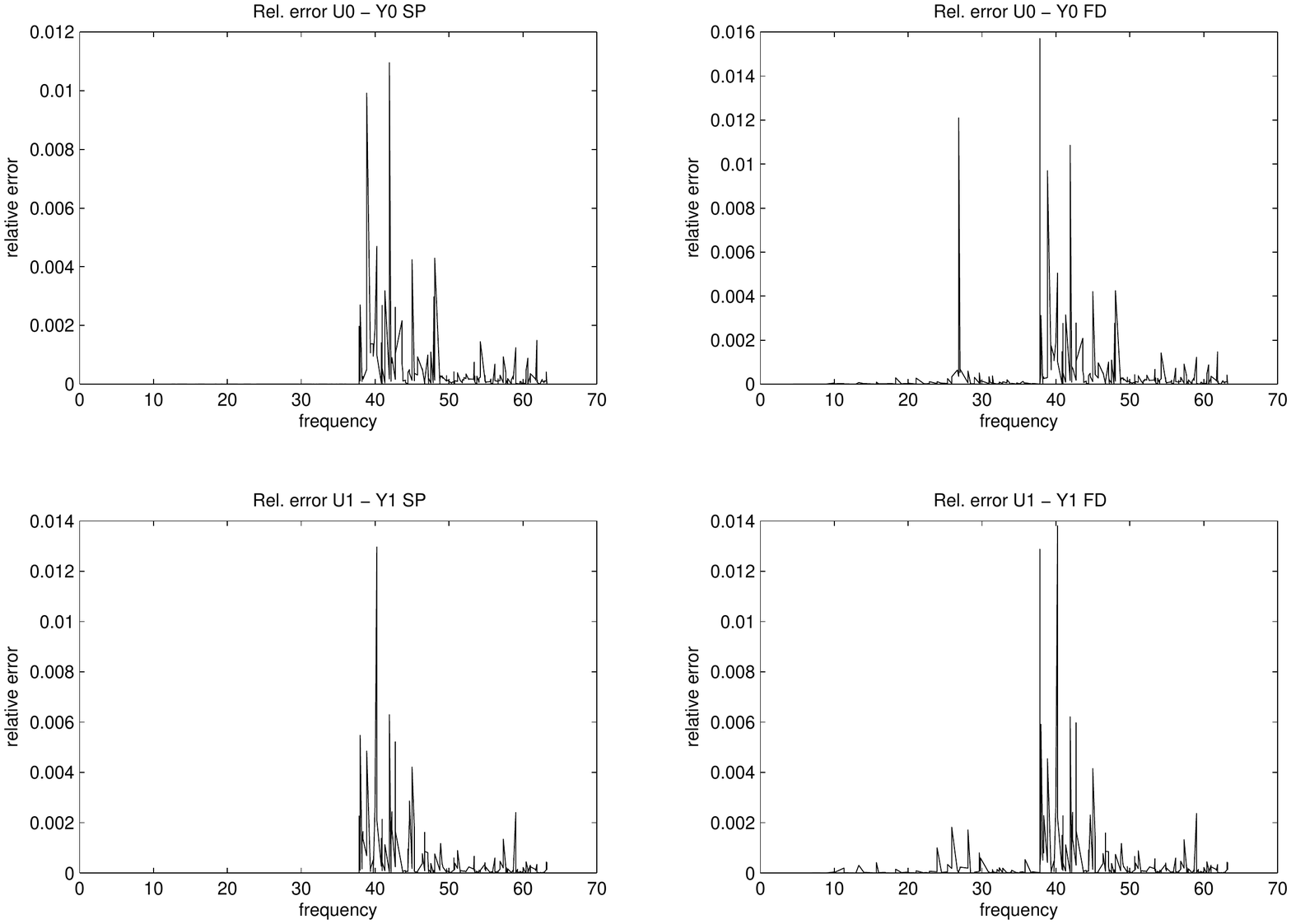}
\caption{\label{fig:valid-errorUY} 
Validation experiments in the square: Relative errors between the spectral coefficients of the original function $u$ and the reconstructed function $y$, computed with exact eigenvalues (left panels) or finite differences eigenvalues (right panels), for $u_{0}$ and $y_{0}$ (top panels) or $u_{1}$ and $y_{1}$ (bottom panels). The errors are plotted as a function of the frequency of the eigenvalues. The computation are performed with 100 eigenvalues, corresponding to a frequency of about 38, the reconstruction with 2000, corresponding to a frequency of about 160. For the readability of the figure, we plot only the major counterparts of the error, i.e. we stop the plot after frequency 63 (300th eigenvalue).}
\end{center}
\end{figure}

In this paragraph, we compare various outputs for our spectral method, when the eigenvalues and eigenvectors are computed either with finite-differences or with exact formulas. In this first experiment, we have $N\times N = 500\times 500$ grid-points and we use $N_{e}=100$ eigenvalues to compute the $G$ and $M_{T}$ matrices. The data $(u_{0},u_{1})$ is as follows:
\be
\label{mneq:11}
\baa{rcl}
u_{0} &=& e_{50}\\
u_{1} &=& 0
\eaa
\ee
Where $e_{n}$ denotes the $n$-th \emph{exact} eigenvector. The control time $T$ is equal to 3, the control domain $U $ is 0.2 wide, and we do not use any smoothing. For reconstruction we use 2000 eigenvalues and eigenvectors.\\
Table \ref{tab:valid:1} shows the condition number of the $M_{T}$ matrices, and reconstruction errors, which are very similar for both experiments.\\
\begin{table}
\begin{center}
\begin{tabular}{|l|l|l|l|}
\hline
Eigenvalues computation  &  Condition Number & Reconstruction error\\ % Exp
\hline
Finite differences & 7.4 & 1.8 \% \\
\hline
Exact & 7.5 & 1.6 \% \\
\hline
\end{tabular}
\caption{\label{tab:valid:1} 
Validation experiments in the square: Condition numbers and validation errors for a 100-eigenvalues-experiment in the square (without smoothing), where the eigenvalues are exact or computed thanks to finite differences.}
\end{center}
\end{table}
Figure \ref{fig:valid-errorUY} shows the relative reconstruction error between the data $u$ and the reconstructed $y$ for both experiments:
\be
\label{mneq:14}
\ba{rcl}
\textrm{Relative reconstruction error SP}_{n} &=& \displaystyle\frac{|U_{0,n}-Y^{sp}_{0,n}|}{\|U_{0}+U_{1}\|}\\
\textrm{Relative reconstruction error FD}_{n} &=& \displaystyle\frac{|U_{0,n}-Y^{fd}_{0,n}|}{\|U_{0}+U_{1}\|}
\ea
\ee
and similarly for $u_{1}$ and $y_{1}$, where $U_{0,n}$ is the $n$-th spectral coefficient of the data $u_{0}$, $Y^{sp}_{0,n}$ is the $n$-th spectral coefficient (in the basis $(\phi_{j})$ defined by formulas (\ref{mneq:base}) in section \ref{sec2.4}) of the reconstructed $y_{0}$ when the control $w$ is obtained thanks to exact eigenvalues and $Y^{fd}_{0,n}$ is the $n$-th spectral coefficient of $y_{0}$ when the control $w$ is obtained thanks to finite differences eigenvalues. The norm $\|U_{0}+U_{1}\|$ in our basis $(\phi_{j})$ is given by:
\be
\label{mneq:13}
\|U_{0}+U_{1}\|^2 = \sum_{n=1}^{N_{e}} U_{0,n}^2 + U_{1,n}^2  
\ee
For exact eigenvalues, we can see that the errors are negligible on the first $100$-th spectral coefficients, and quite small on the next ones. We have similar results for finite differences eigenvalues, except that we have an error on the 50-th coefficient. This error does not occur when the reconstruction is done with the same finite differences eigenvectors basis, and it can probably be explained as follows: to compute the reconstructed $y$ from the finite difference control $w$, we first compute an approximation of $w$ as a function of $(x,y)$ (i.e. on the grid) from its spectral coefficients (on the finite differences eigenvectors basis), then we compute the coefficients of this function on the exact basis (thanks to a very simple integration formula). We thus introduce two sources of errors, projection on the grid and projection on the exact basis, which do not have anything to do with our spectral Galerkin method. Therefore we will not discuss the matter in further detail here. 

\subsubsection{Impact of the number of eigenvalues}\label{sec3.5-2}

In this paragraph, we still use the same data (\ref{mneq:11}), but the number of eigenvalues and eigenvectors $N_{e}$ used to compute the $M_{T}$ matrices is varying. Table \ref{tab:valid:2} shows $M_{T}$ condition numbers and reconstruction errors for various $N_{e}$ with exact or finite-differences-computed eigenvalues. The reconstruction is still performed with 2000 exact eigenvalues. We can see that the finite differences eigenvalues lead to almost as good results as exact eigenvalues. We also observe in both cases the decrease of the reconstruction error with an increasing number of eigenvalues, as predicted in lemma \ref{lem3}. A 5\% error is obtained with 70 eigenvalues (the input data being the 50-th eigenvalue), and 100 eigenvalues lead to less than 2\%.
\begin{table}
\begin{center}
\begin{tabular}{|c|cc|cc|}
\hline
         &  \multicolumn{2}{|c|}{Condition Number} & \multicolumn{2}{|c|}{Reconstruction error}\\ 
$N_{e}$  &     Exact & Finite differences &  Exact & Finite differences\\
\hline 52  & 6.5 & 6.4 & 29.2\% & 29.0\%\\
\hline 55  & 6.6 & 6.6 & 17.7\% & 17.6\%\\
\hline 60  & 6.6 & 6.6 & 17.3\% & 17.2\%\\
\hline 70  & 7.1 & 7.0 & 4.9 \% & 5.0 \%\\
\hline 80  & 7.3 & 7.2 & 3.2 \% & 3.3 \%\\
\hline 100 & 7.5 & 7.4 & 1.6 \% & 1.8 \% \\
\hline 200 & 8.3 & 8.3 & 0.5 \% & 1.0 \% \\
\hline 300 & 8.8 &     & 0.3 \% &       \\
\hline 500 & 9.5 &     & 0.2 \% &     \\
\hline
\end{tabular}
\caption{\label{tab:valid:2} 
Validation experiments in the square: Condition numbers and validation errors for various numbers $N_{e}$ of eigenvalues used to compute the control function, where the eigenvalues are exact or computed thanks to finite differences. The input data is the $50$-th eigenvalue, and the reconstruction is performed with 2000 eigenvalues.}
\end{center}
\end{table}

%%%%%%%%%%%%%%%%%%%%%%%%%%%%
\section{Numerical experiments}
\label{sec4}

\subsection{Frequency localization}\label{sec4.1}
\label{sec:onemode}

In this subsection, the geometry (square) as well as the number of eigenvalues used (200 for HUM, 2000 for verification) are fixed. Note also that in this paragraph we use only exact eigenvalues for HUM and verification.\\
The data is also fixed to a given eigenmode, that is:
\be
u_0 = e_{50} \quad u_1 = 0
\ee
where $e_{n}$ is the n-th eigenvector of $-\Delta$ on the square.\\

\begin{figure}
\begin{center}
\includegraphics[angle=-90,width=\textwidth]{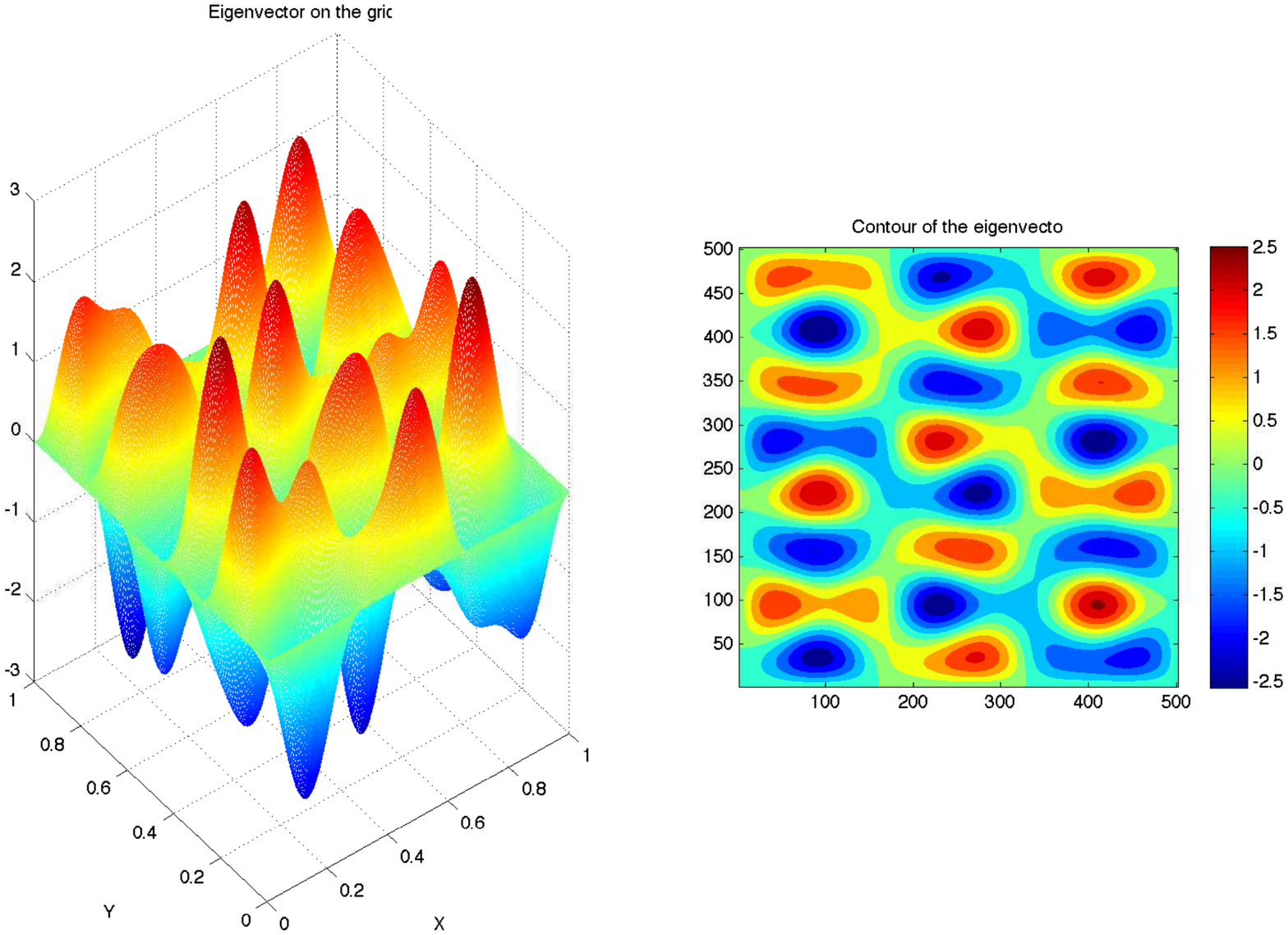}
\caption{\label{fig:1mode50} 
Representation on the grid in 3D (left panel) or contour plot (right panel) of the 50-th eigenvector in the square.}
\end{center}
\end{figure}

\begin{figure}
\begin{center}
\includegraphics[width=\textwidth,height=.6\textwidth]{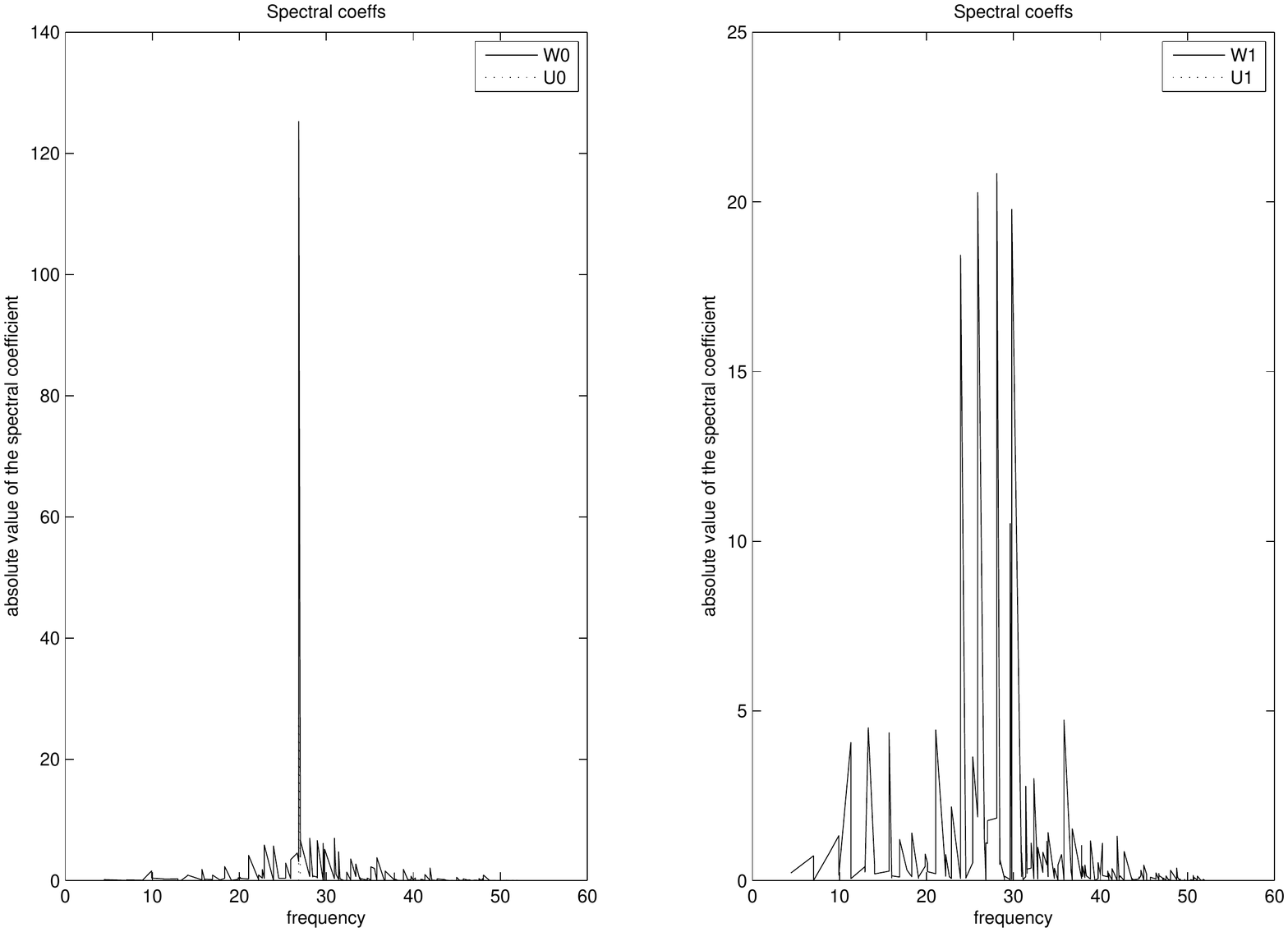}
\includegraphics[width=\textwidth,height=.6\textwidth]{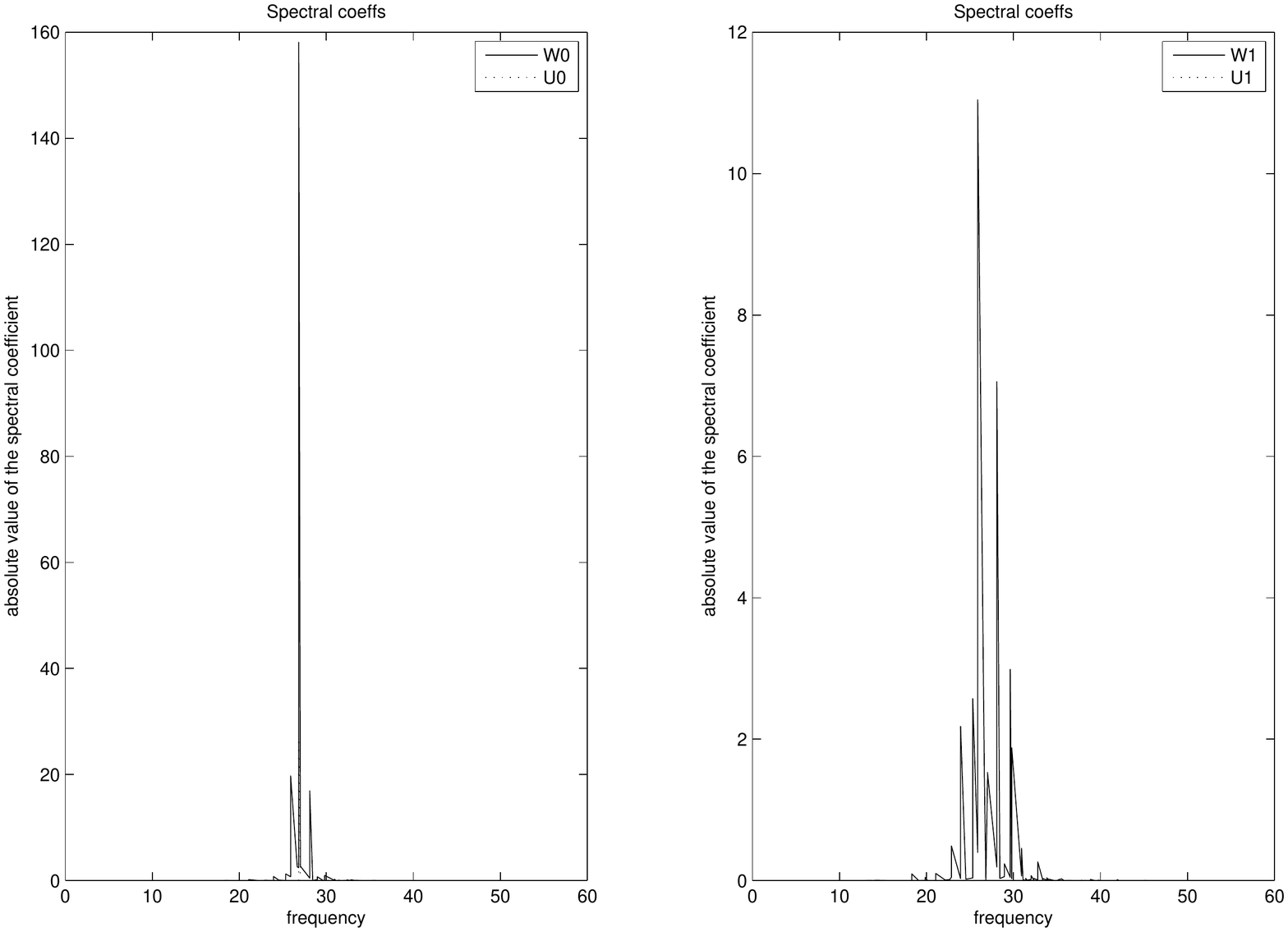}
\caption{\label{fig:1mode:sq:freq1} 
One-mode experiment in the square: localization of the Fourier frequences of $(u_0,u_1)$ (dashed line) and $(w_0,w_1)$ (solid line) for a given time $T$ and a given domain $U $ without smoothing (top) and with time- and space-smoothing (bottom). The x-coordinate represents the eigenvalues. The input data $u_{0}$ is equal to the 50-th eigenvector, equal to an eigenvalue of about 26.8, and $u_{1}=0$.}
\end{center}
\end{figure}

The first output of interest is the spreading of $w$ spectral coefficients, compared to $u$. Figure \ref{fig:1mode:sq:freq1} shows the spectral coefficients of the input $(u_{0},u_{1})$ and the control $(w_{0},w_{1})$ with and without smoothing. As predicted by theorem \ref{thm2} and lemma \ref{lem3} we can see that the main coefficient of $(w_{0},w_{1})$ is the 50-th of $w_{0}$, and also that the smoothing noticeably improves the localization of $w$.\\

\begin{figure}
\begin{center}
\includegraphics[width=\textwidth,height=.6\textwidth]{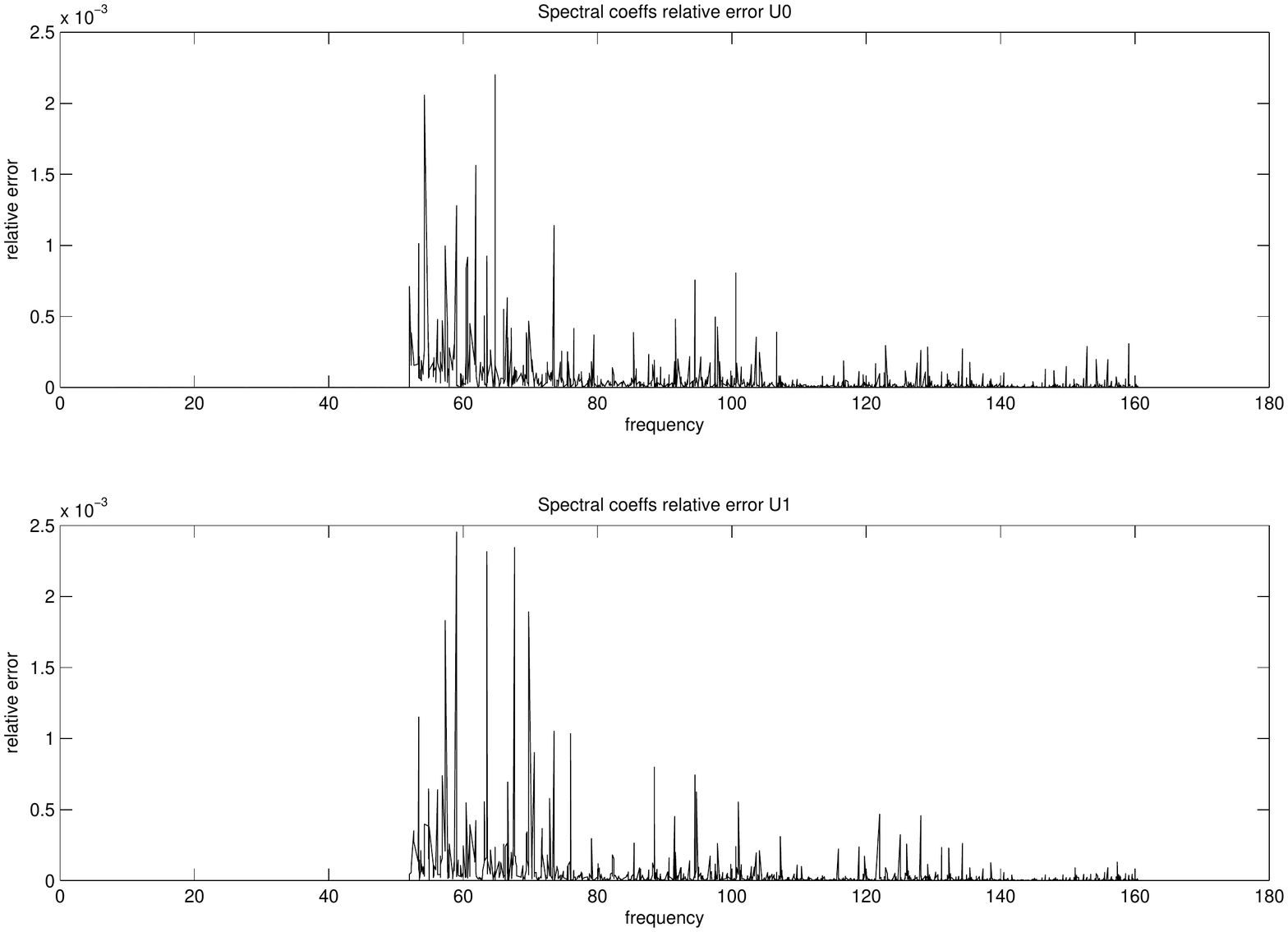}
\includegraphics[width=\textwidth,height=.6\textwidth]{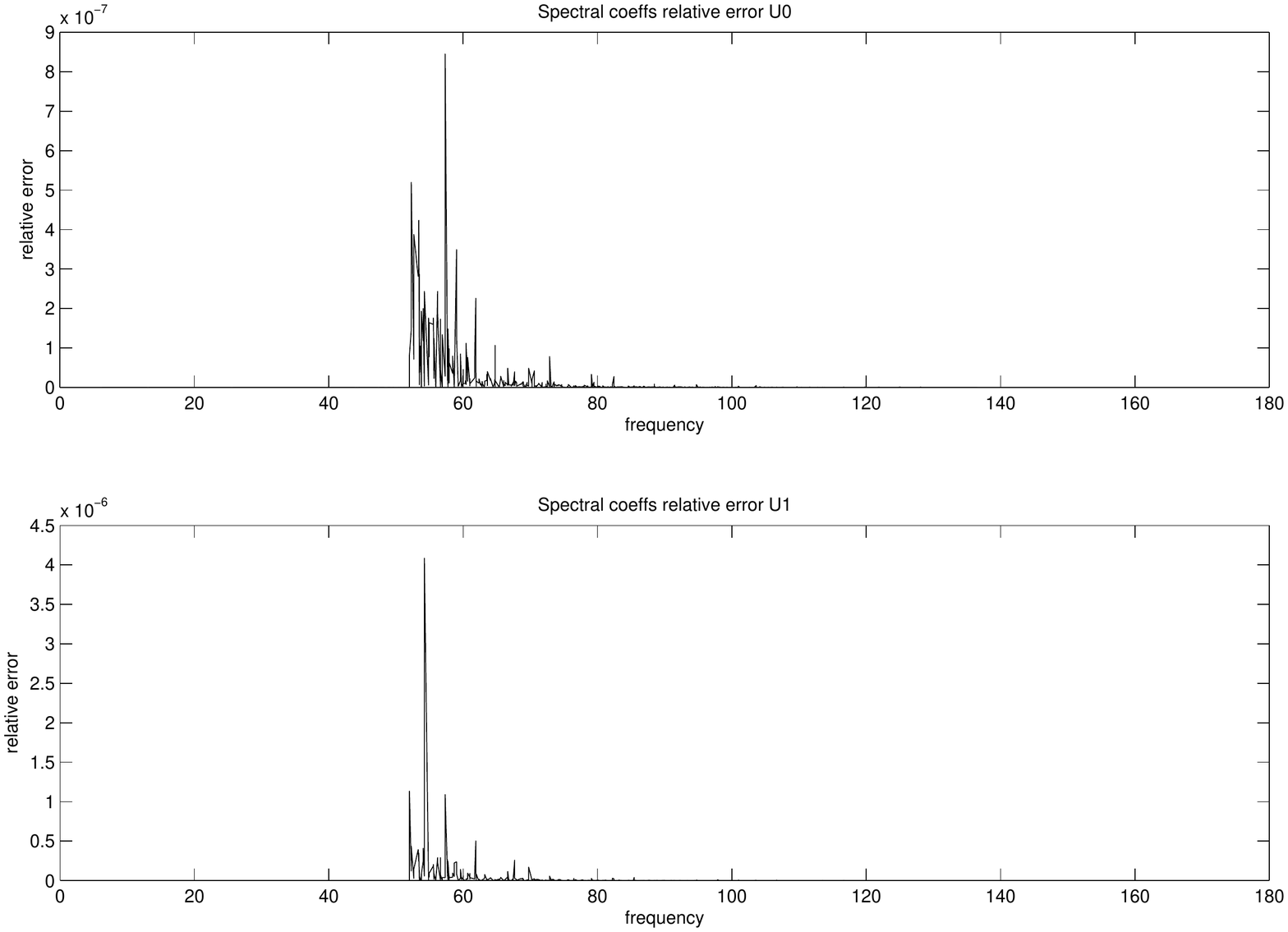}
\caption{\label{fig:1mode:sq:freq2} 
One-mode experiment in the square: localization of the Fourier coefficients of $(u_0-y_0,u_1-y_1)$, where $u$ is the data and $y$ is the reconstructed function obtained from the control function $w$, for a given time $T$ and a given domain $U$ without smoothing (top panels) and with time- and space-smoothing (bottom panels).}
\end{center}
\end{figure}

Similarly we can look at the spectral coefficients of the reconstruction error. Figure \ref{fig:1mode:sq:freq2} presents the reconstruction error (see paragraph \ref{sec3.4recerror} for a definition) with or without smoothing. We notice that the errors occur mostly above the cutoff frequency (used for $M_{T,\omega}$ computation, and thus for the control computation). Another important remark should be made here: the smoothing has a spectacular impact on the frequency localization of the error, as well as on the absolute value of the error (maximum of $2.10^{-3}$ without smoothing, and $8.10^{-7}$ with smoothing), as announced in theorem \ref{thm2} and lemma \ref{lem3}.\\

\begin{rmq}
For other domains, such as the disc and trapezoid, as well as other one-mode input data, we obtain similar results. The results also remain the same if we permute $u0$ and $u1$, i.e. if we choose $u0=0$ and $u1$ equal to one fixed mode.
\end{rmq}

%%%%%%%%%%%%%%%%%%%%%%%%%%%%%%%%%%%%%%

\subsection{Space localization}\label{sec4.2}

\subsubsection{Dirac experiments}\label{sec4.2-1}

\begin{figure}
\begin{center}
\includegraphics[angle=-90,width=\textwidth]{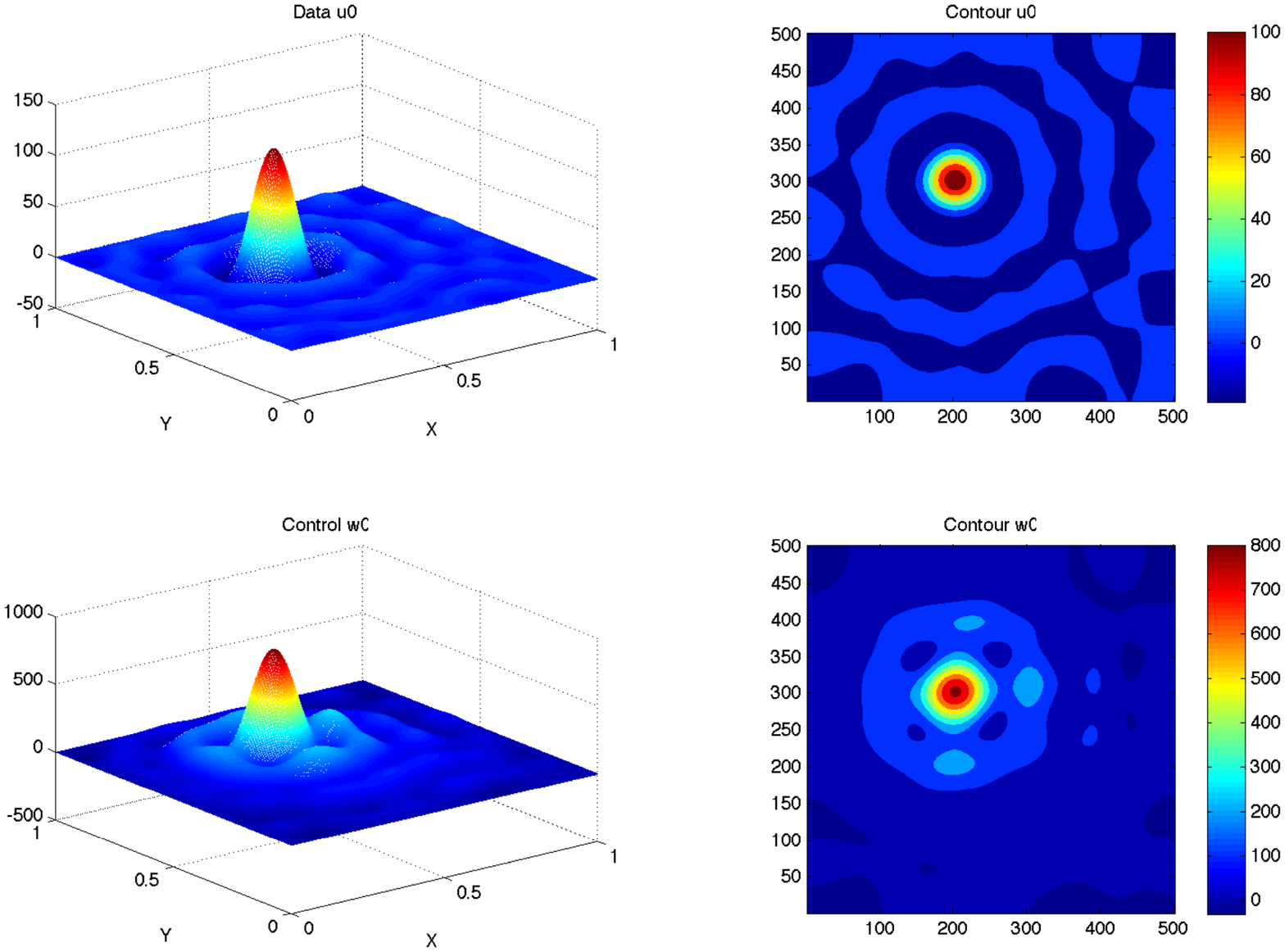}
\caption{\label{fig:spaceloc:sq:1} 
Space  localization of the data $u_{0}$ (top panels) and the control $w_{0}$ (bottom panels), for a Dirac experiment in the square, with exact eigenvalues. These plots correspond to an experiment without smoothing, but it is similar with smoothing. Left panels represent 3D view, and right panels show contour plots.}
\end{center}
\end{figure}

\begin{figure}
\begin{center}
\includegraphics[angle=-90,width=\textwidth]{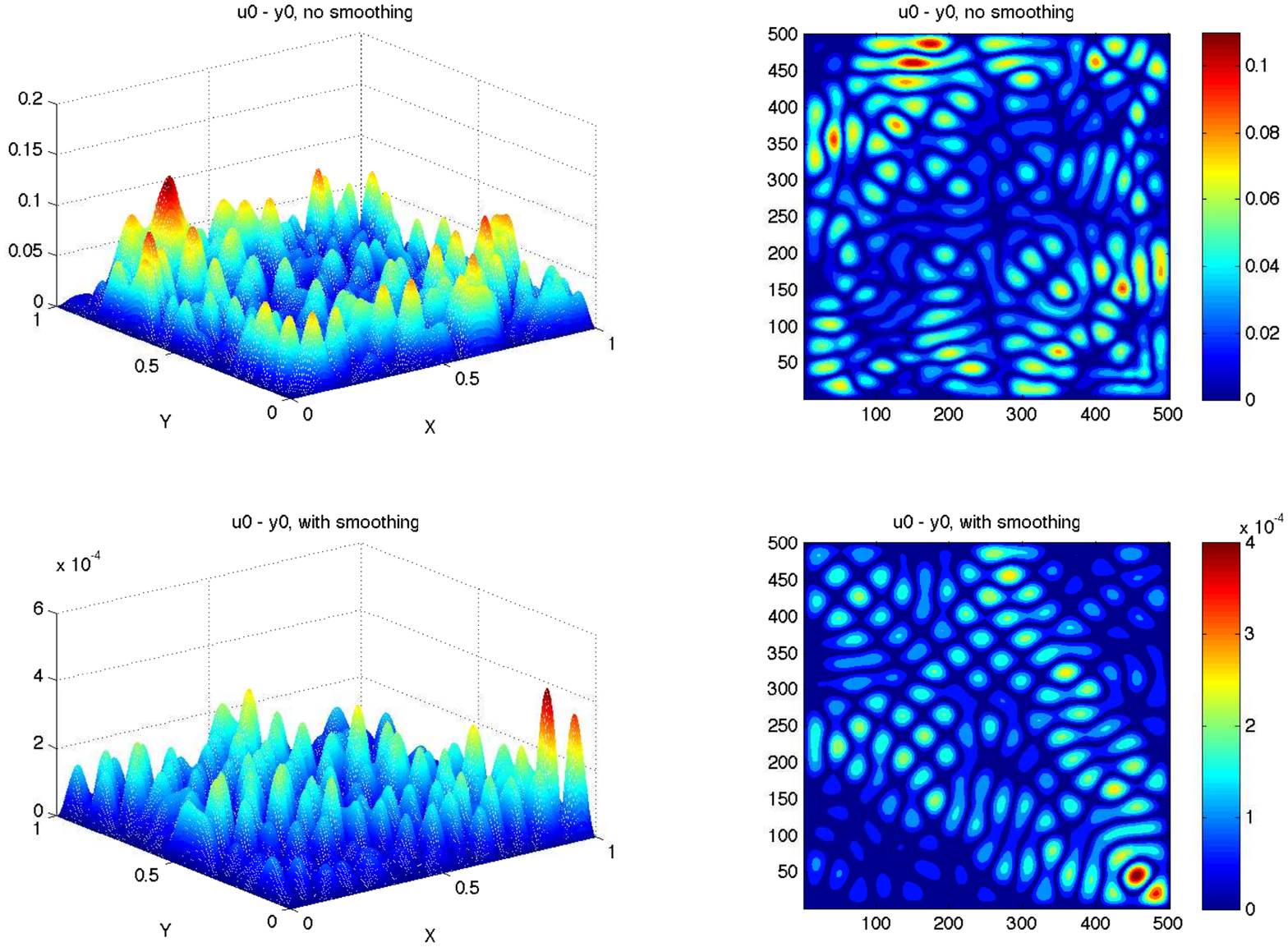}
\caption{\label{fig:spaceloc:sq:3} 
Difference between the data $u_{0}$ and the reconstructed function $y_{0}$ without smoothing (top panels) and with smoothing (bottom panels) for a dirac experiment in the square, with exact eigenvalues. Left panels represent 3D view, and right panels show contour plots.}
\end{center}
\end{figure}

In this section we investigate the localization in space. To do so, we use ``Dirac" functions $\delta_{(x,y) = (x_{0},y_{0})}$ as data, or more precisely truncations to a given cutoff frequency of Dirac functions:
\be
\label{mneq:11bis}
\baa{rcl}
u_{0} &=& \sum_{i=1}^{N_{i}} e_{n}(x_{0},y_{0})\, e_{n}\\
u_{1} &=& 0
\eaa
\ee
where $N_{i}$ is the index corresponding to the chosen cutoff frequency, with $N_{i}=100$ or $120$ in the sequel. Figure \ref{fig:spaceloc:sq:1} shows the data $u_{0}$ and the control $w_{0}$ in the square with exact eigenvalues, without smoothing, the results being similar with smoothing. We can see that the support of $w_{0}$ is very similar to $u_{0}$'s. Figure \ref{fig:spaceloc:sq:3} presents the reconstruction error associated to this experiment. We can see as before that the smoothing produces highly reduced errors. 

\begin{figure}
\begin{center}
\includegraphics[angle=-90,width=\textwidth]{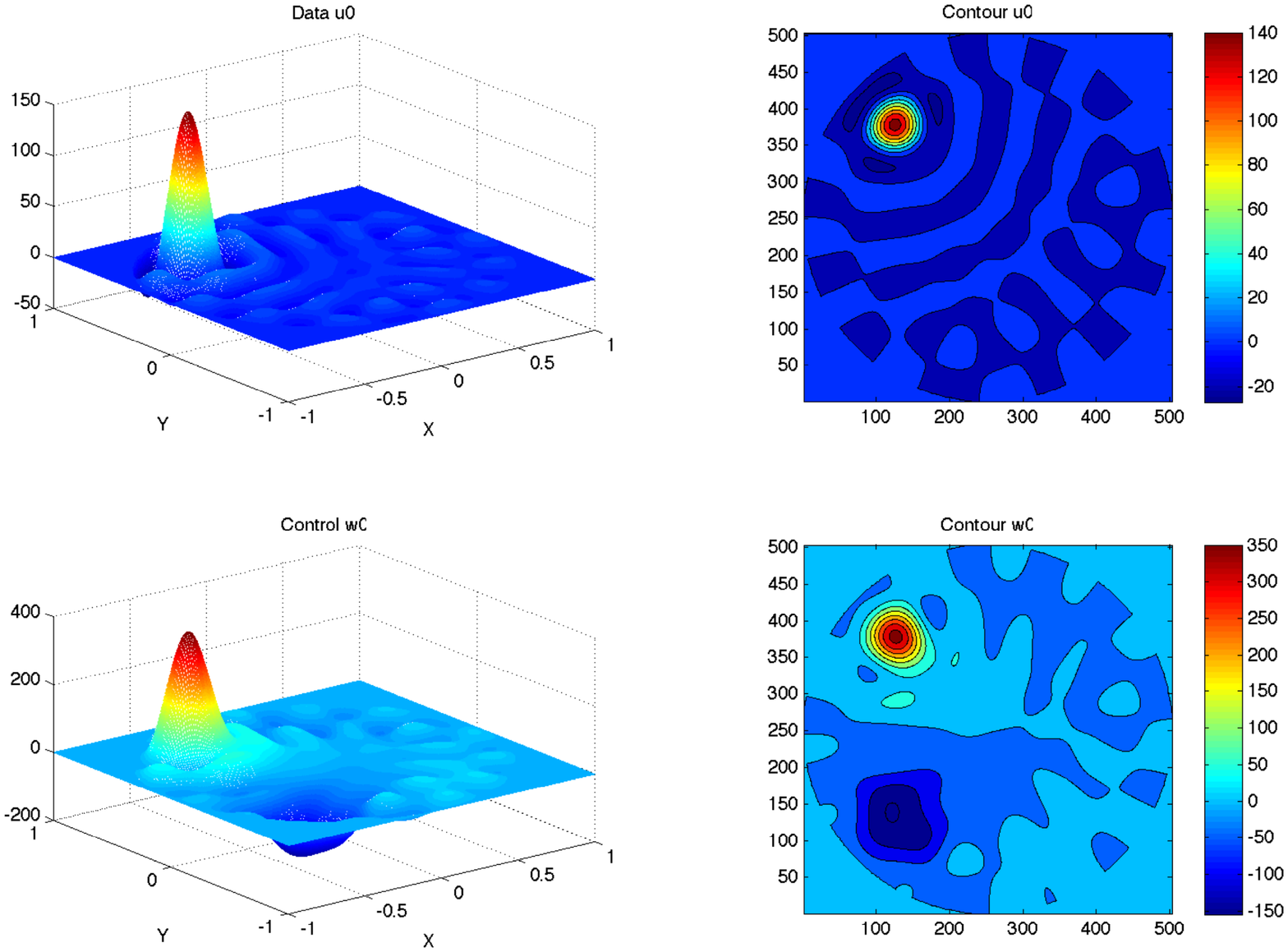}
\caption{\label{fig:spaceloc:disc:1} 
Space localization of the data $u_{0}$ (top panels) and the control $w_{0}$ (bottom panels), for a dirac experiment  in the disc. These plots correspond to an experiment without smoothing, but it is similar with smoothing. Left panels represent 3D view, and right panels show contour plots. In this experiment, the input data is defined with $N_{i}=100$ eigenvectors. }
\end{center}
\end{figure}

\begin{figure}
\begin{center}
\includegraphics[angle=-90,width=\textwidth]{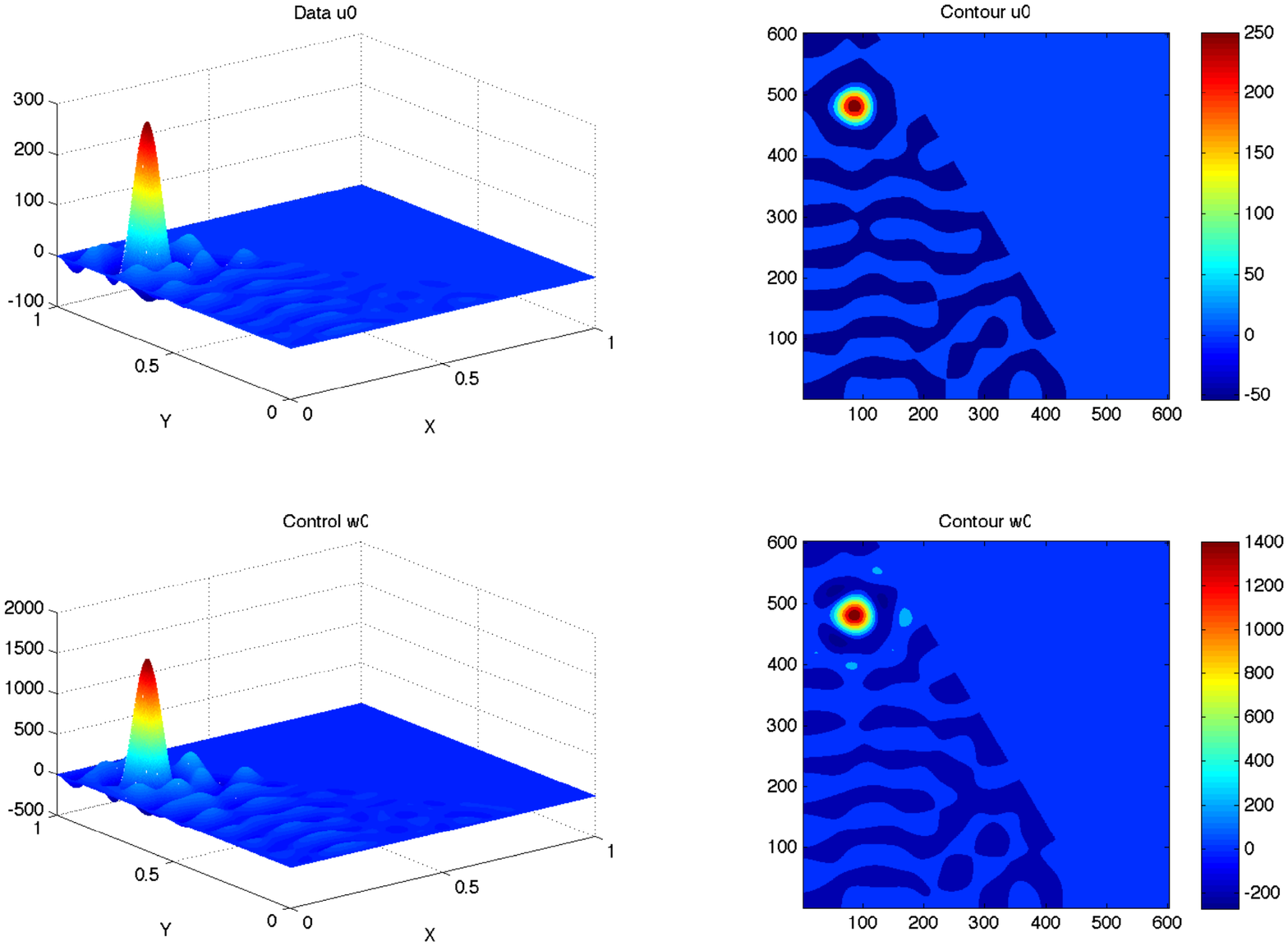}
\caption{\label{fig:spaceloc:trapez:1} 
Space localization of the data $u_{0}$ (top panels) and the control $w_{0}$ (bottom panels), for a dirac experiment  in the trapezoid. These plots correspond to an experiment without smoothing, but it is similar with smoothing. Left panels represent 3D view, and right panels show contour plots. In this experiment, the input data is defined with $N_{i}=120$ eigenvectors.}
\end{center}
\end{figure}

\begin{figure}
\begin{center}
\includegraphics[angle=-90,width=\textwidth]{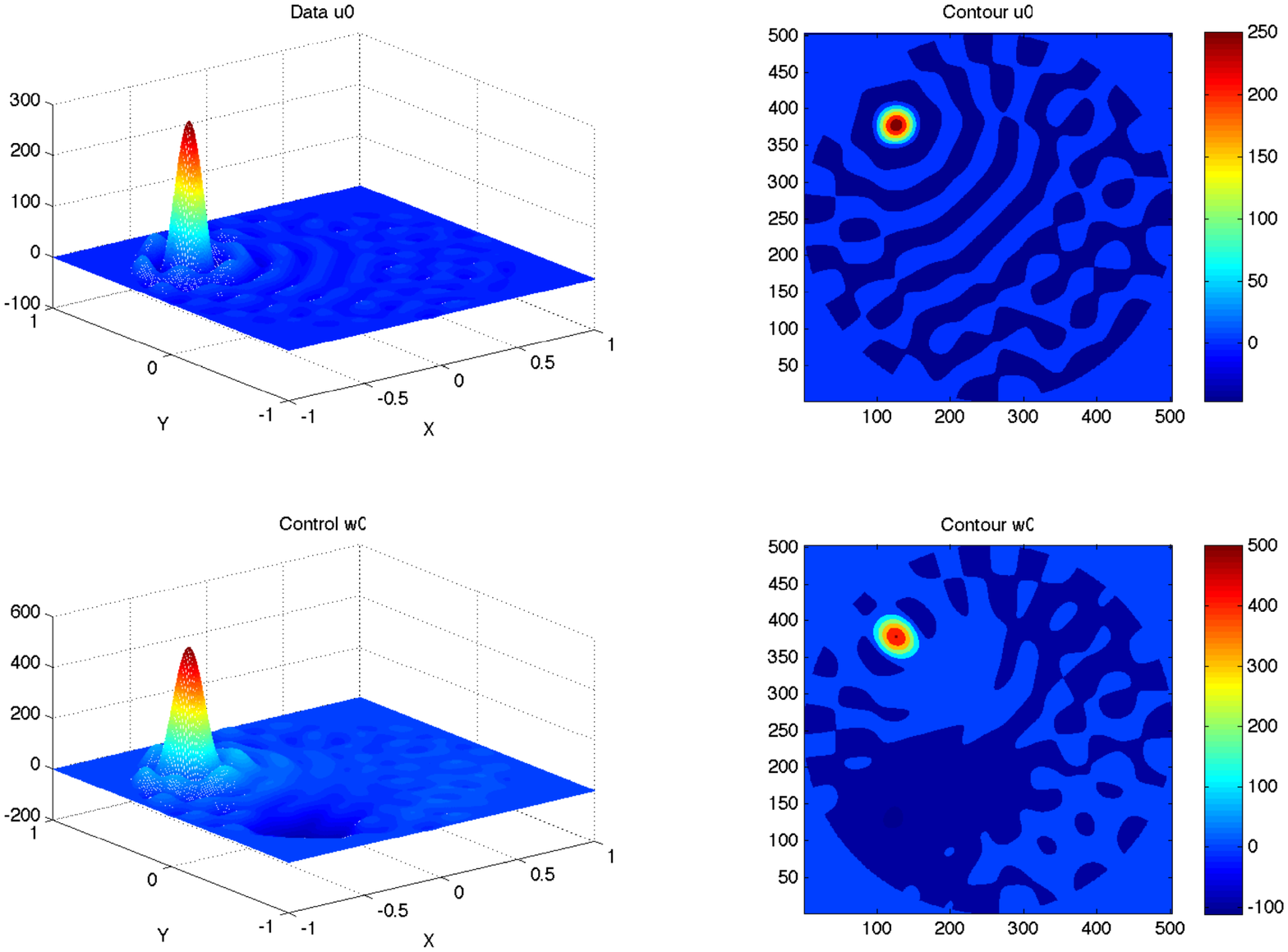}
\caption{\label{fig:spaceloc:disc:1bis} 
Space localization of the data $u_{0}$ (top panels) and the control $w_{0}$ (bottom panels), for a dirac experiment  in the disc. These plots correspond to an experiment with smoothing, and it is similar without smoothing. Left panels represent 3D view, and right panels show contour plots. In this experiment, the input data is defined with $N_{i}=200$ eigenvectors.}
\end{center}
\end{figure}

\begin{figure}
\begin{center}
\includegraphics[angle=-90,width=\textwidth]{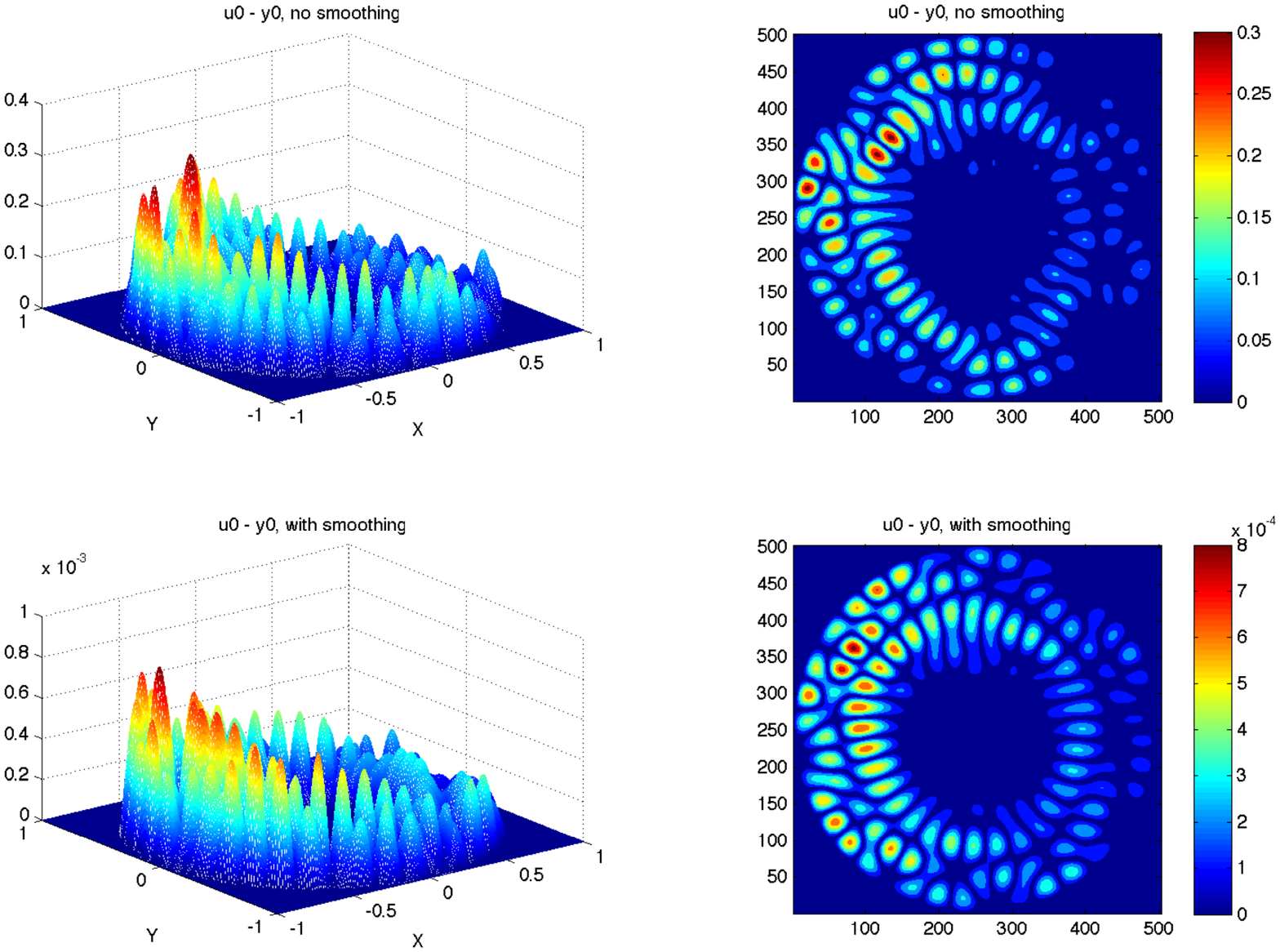}
\caption{\label{fig:spaceloc:disc:2} 
Difference between the data $u_{0}$ and the reconstructed function $y_{0}$  for a dirac experiment in the disc without smoothing (top panels) and with time- and space-smoothing (bottom panels). Left panels represent 3D view, and right panels show contour plots. In this experiment, the input data is defined with $N_{i}=100$ eigenvectors.}
\end{center}
\end{figure}

Similarly, we performed experiments with numerical approximation of a Dirac function as input data in the disc and in a trapezoid. 
Figures \ref{fig:spaceloc:disc:1} and \ref{fig:spaceloc:trapez:1} present the space-localization of $u_{0}$ and $w_{0}$ without smoothing (we get similar results with smoothing). As previously, the control $w_{0}$ is supported by roughly the same area than the input $u_{0}$. In the disc we can see a small disturbance, located in the symmetric area of the support of $u_{0}$ with respect to the control domain $U$. However, this error does not increase with $N_{i}$, as we can see in figure \ref{fig:spaceloc:disc:1bis} (case $N_{i}=200$) so it remains compatible with conjecture \ref{conj1}.\\
Figure \ref{fig:spaceloc:disc:2} shows the reconstruction errors for these experiments, with or without smoothing. As before we notice the high improvement produced by the smoothing. We get similar errors in the trapezoid.

\subsubsection{Box experiments in the square}\label{sec4.2-2}
\begin{figure}
\begin{center}
\includegraphics[angle=-90,width=\textwidth]{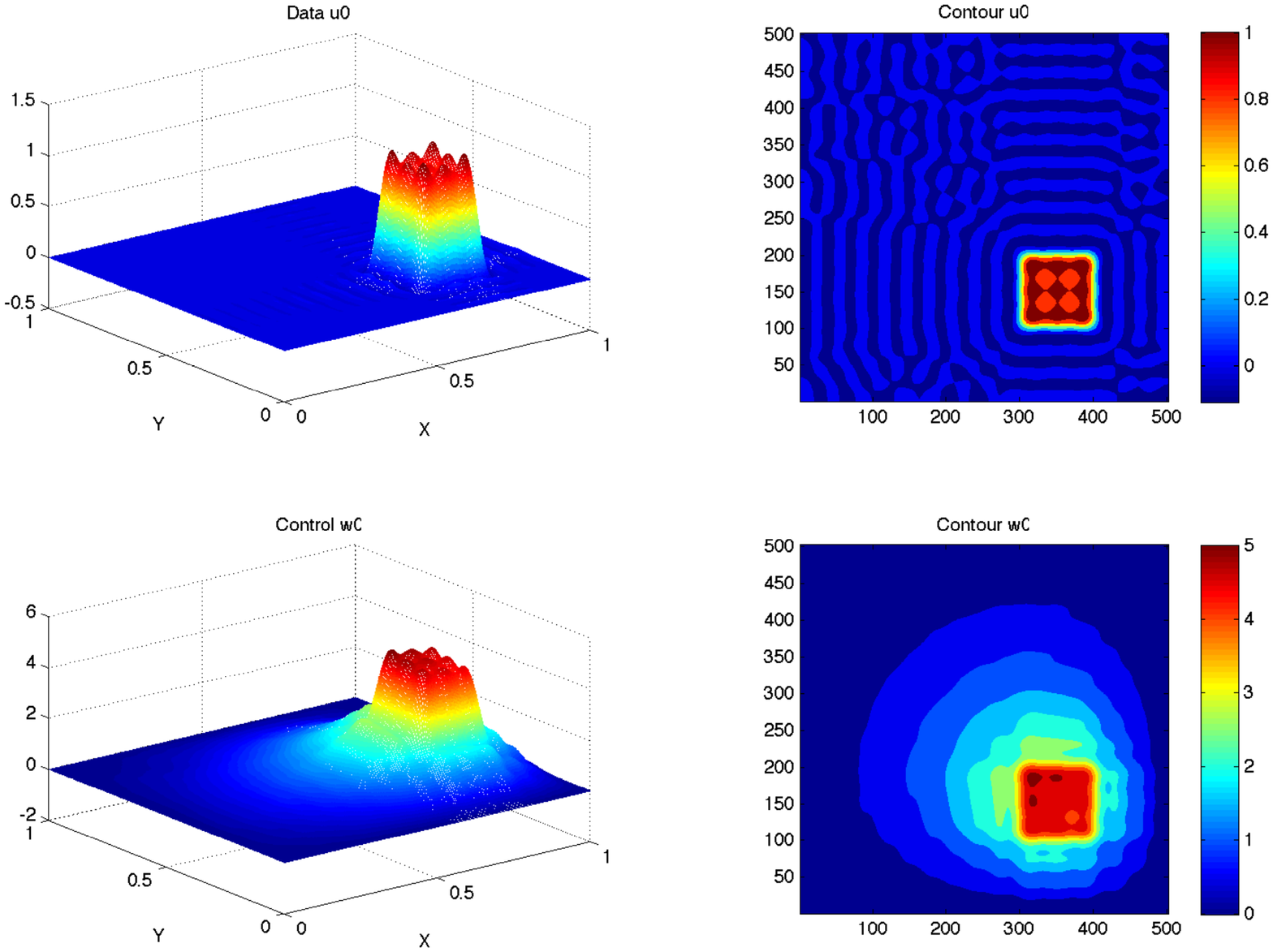}
\caption{\label{fig:spaceloc:sqbox:1} 
Space localization of the control function $w_{0}$ (bottom panels) with respect to the data $u_{0}$ (top panels), in the square, without smoothing: 3D plots on the left, and contour plots on the right.}
\end{center}
\end{figure}

\begin{figure}
\begin{center}
\includegraphics[angle=-90,width=\textwidth]{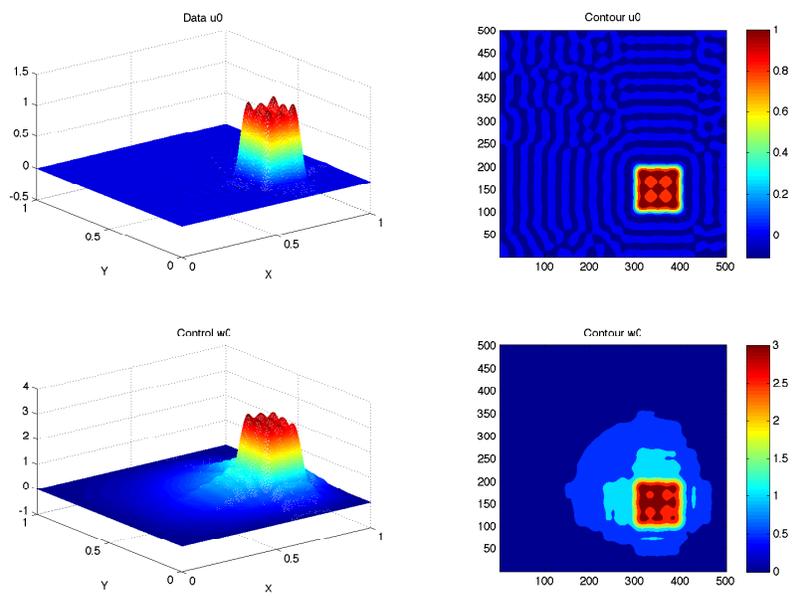}
\caption{\label{fig:spaceloc:sqbox:2} 
Space localization of the control function $w_{0}$ (bottom panels) with respect to the data $u_{0}$ (top panels), in the square, with smoothing. Left panels represent 3D view, and right panels show contour plots.}
\end{center}
\end{figure}

\begin{figure}
\begin{center}
\includegraphics[angle=-90,width=\textwidth]{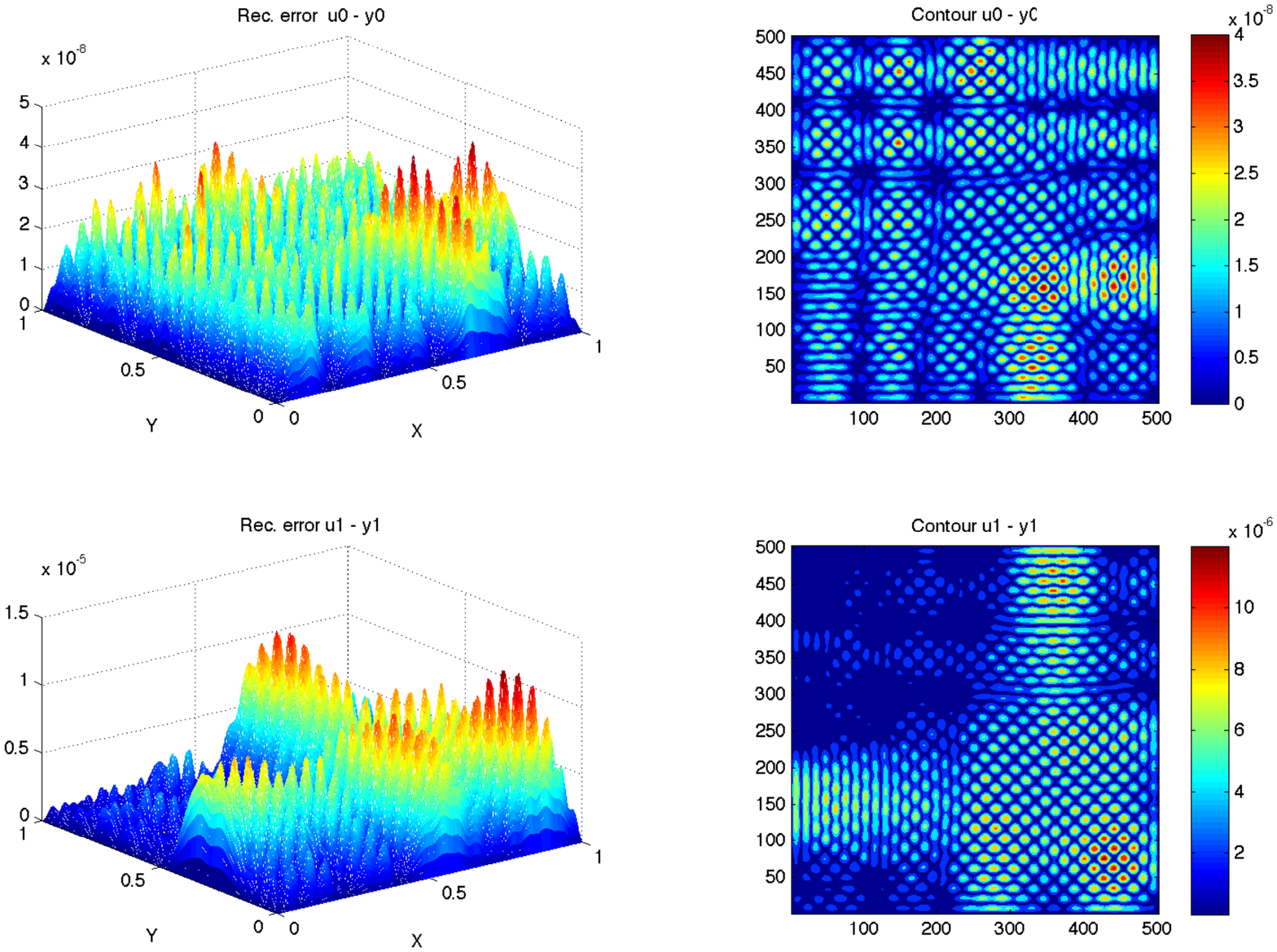}
\caption{\label{fig:spaceloc:sqbox:3} 
Difference between the data $u_{0}$ and the reconstructed function $y_{0}$ (top panels) and $u_{1}$ and $y_{1}$ (bottom panels) with smoothing in the square. The data is the identity function of a square whose edges are parallel to the $x$ and $y$ axes. Left panels represent 3D view, and right panels show contour plots.}
\end{center}
\end{figure}

\begin{figure}
\begin{center}
\includegraphics[angle=-90,width=\textwidth]{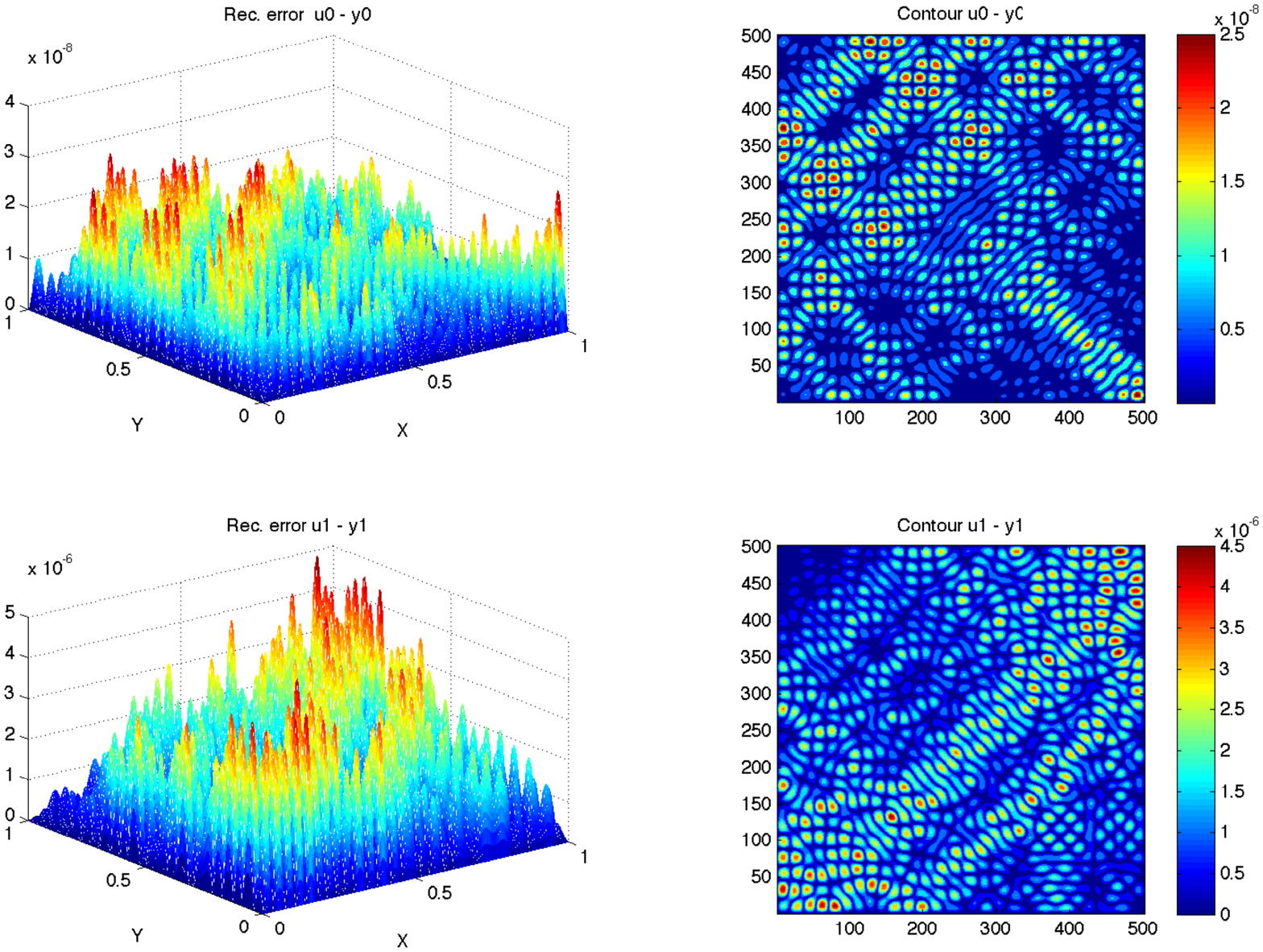}
\caption{\label{fig:spaceloc:sqbox:4} 
Difference between the data $u_{0}$ and the reconstructed function $y_{0}$ (top panels) and $u_{1}$ and $y_{1}$ (bottom panels) with smoothing in the square. The data is the identity function of a square whose edges are parallel to the diagonals of the square. Left panels represent 3D view, and right panels show contour plots.}
\end{center}
\end{figure}

In this paragraph we consider the case $u_{0} = {\bf 1}_{\textrm{box}}$, where $\textrm{box} = [0.6,0.8]\times[0.2,0.4]$ is a box in the square. The control domain $U$ is $0.1$ wide: $U=\{x < 0.1 \textrm{ and } y > 0.9 \}$. 
These experiments were performed in the square with 1000 exact eigenvalues for the $M_{T}$ matrix computation, the input data $u_{0}$ being defined thanks to 800 eigenvalues.
Figures \ref{fig:spaceloc:sqbox:1} and \ref{fig:spaceloc:sqbox:2} show the space localization of the data $u_{0}$ and the control $w_{0}$ without and with smoothing. As before we can notice that the space localization is preserved, and that with smoothing the support of $w_{0}$ is more sharply defined.
Figures \ref{fig:spaceloc:sqbox:3} and \ref{fig:spaceloc:sqbox:4} show the reconstruction errors for two different data, the first being the same as in figure \ref{fig:spaceloc:sqbox:1}, and the second being similar but rotated by $\pi/4$. We show here only the case with smoothing, the errors being larger but similarly shaped without. We can notice that the errors lows and highs are located on a lattice whose axes are parallel to the box sides. This is compatible with the structure of the wave-front set associated to both input data.

%%%%%%%%%%%%%%%%%%%%%%%%%%%%%%%%%%%%%%

\subsection{Reconstruction error}\label{sec4.3}

\begin{figure}
\begin{center}
\includegraphics[width=1\textwidth]{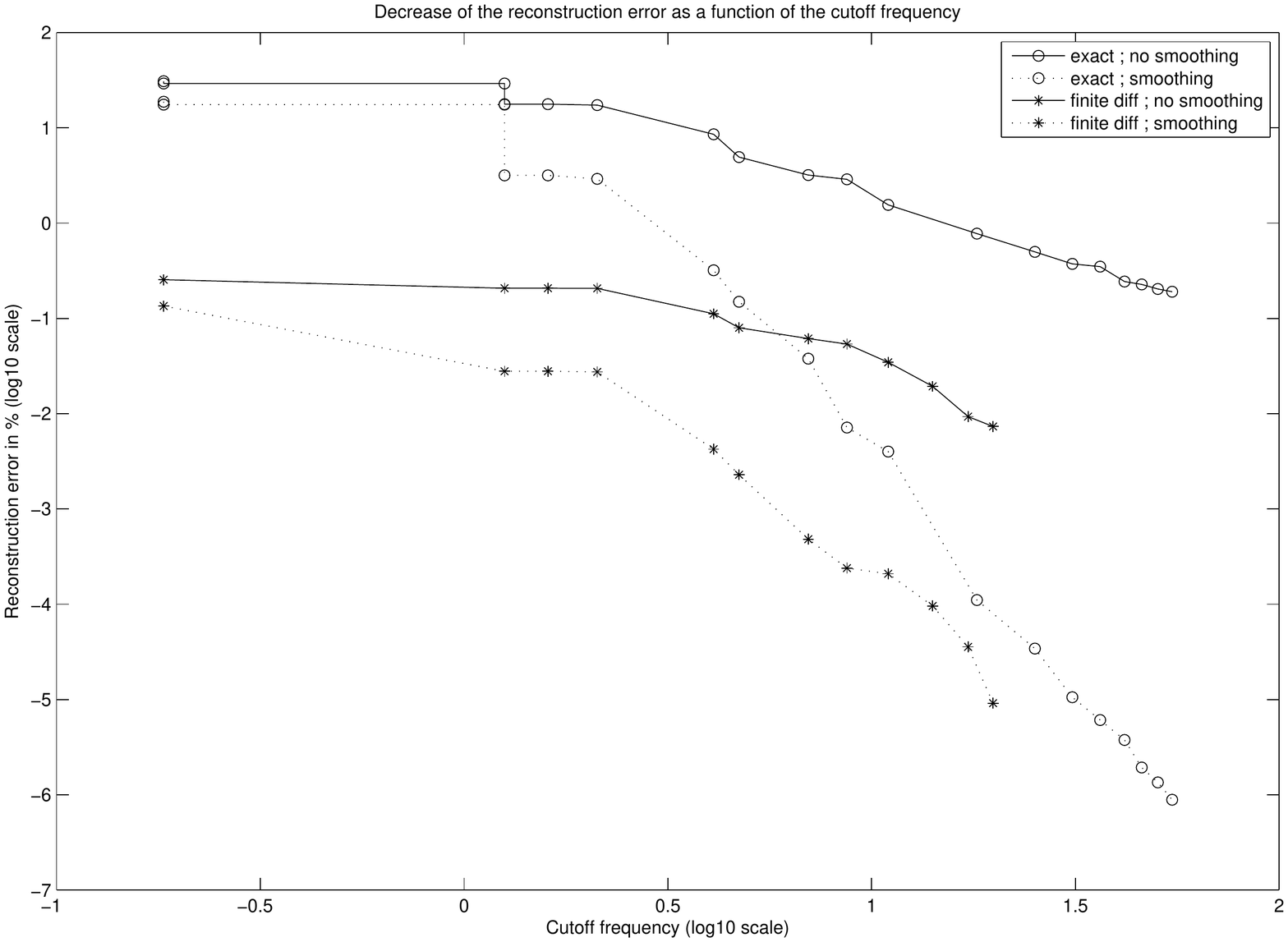}
\caption{\label{fig:rec-error4} 
Reconstruction errors for the finite differences and exact methods, as a function of the cutoff frequency (i.e., the largest eigenvalue used for the control computation), with or without time- and space-smoothing. }
\end{center}
\end{figure}

\begin{figure}
\begin{center}
\includegraphics[width=1\textwidth]{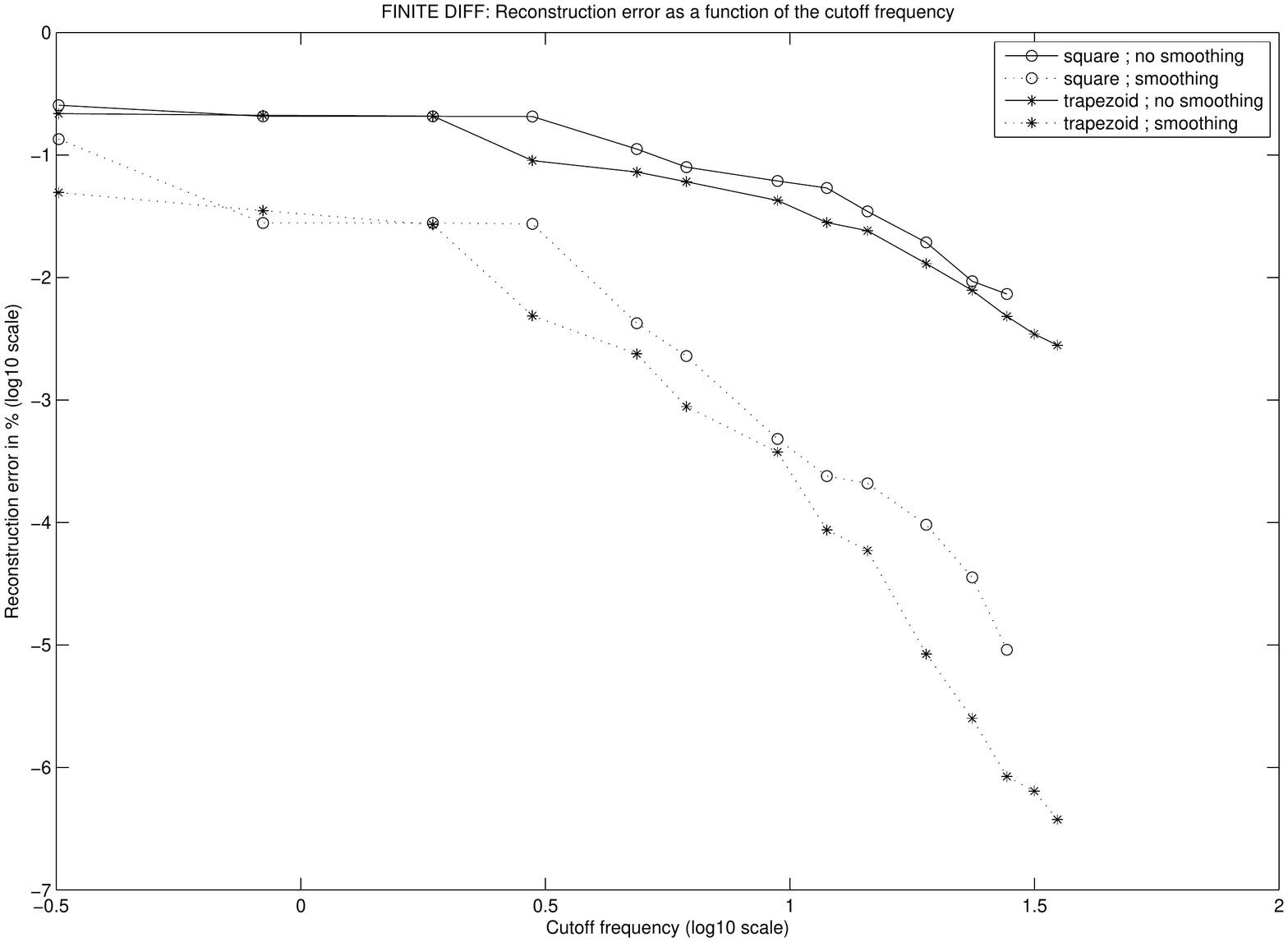}
\caption{\label{fig:rec-error5} 
Reconstruction errors for the finite differences method, in the square and in the trapezoid, as a function of the cutoff frequency (i.e., the largest eigenvalue used for the control computation), with or without time- and space-smoothing. }
\end{center}
\end{figure}

In this section we investigate lemma \ref{lem3} or more precisely the subsequent remark \ref{rmklem3}. This remark states that the reconstruction error should decrease as the inverse of the cutoff frequency without smoothing, and as the inverse of the fifth power of the cutoff frequency with smoothing. To investigate this, we perform a ``one-mode" experiment (see paragraph \ref{sec:onemode}) using the 50-th mode as input data. We then compute the control with an increasing cutoff frequency, up to 47 (finite differences case) or 82 (exact case), and we compute the reconstruction error, thanks to a larger cutoff frequency (52 in the finite differences case, or 160 in the exact case). \\
Figure \ref{fig:rec-error4}  represents the reconstruction error (with exact or finite differences eigenvalues) as a function of the cutoff frequency (i.e., the largest eigenvalue used for the control function computation). Figure \ref{fig:rec-error5} presents the same results (with finite differences eigenvalues only) for two different geometries: the square, and the trapezoid (general domain). The log scale allows us to see that the error actually decreases as the inverse of the cutoff frequency without smoothing, and as the inverse of the fifth power of the cutoff frequency with smoothing, according to remark \ref{rmklem3}.\\

\subsection{Energy of the control function}\label{sec4.4}
\begin{figure}
\begin{center}
\includegraphics[width=.85\textwidth]{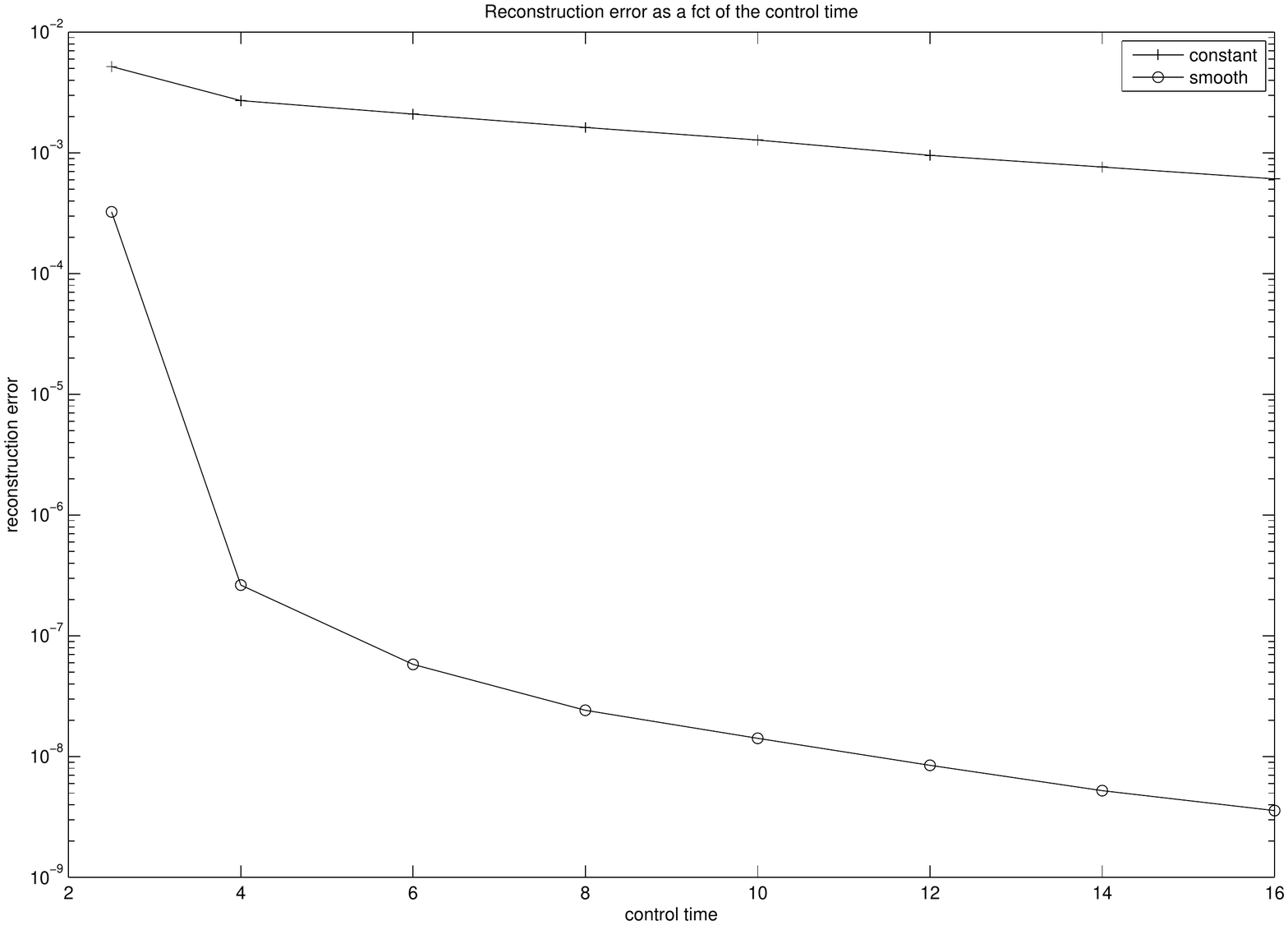}\\
\includegraphics[width=.85\textwidth]{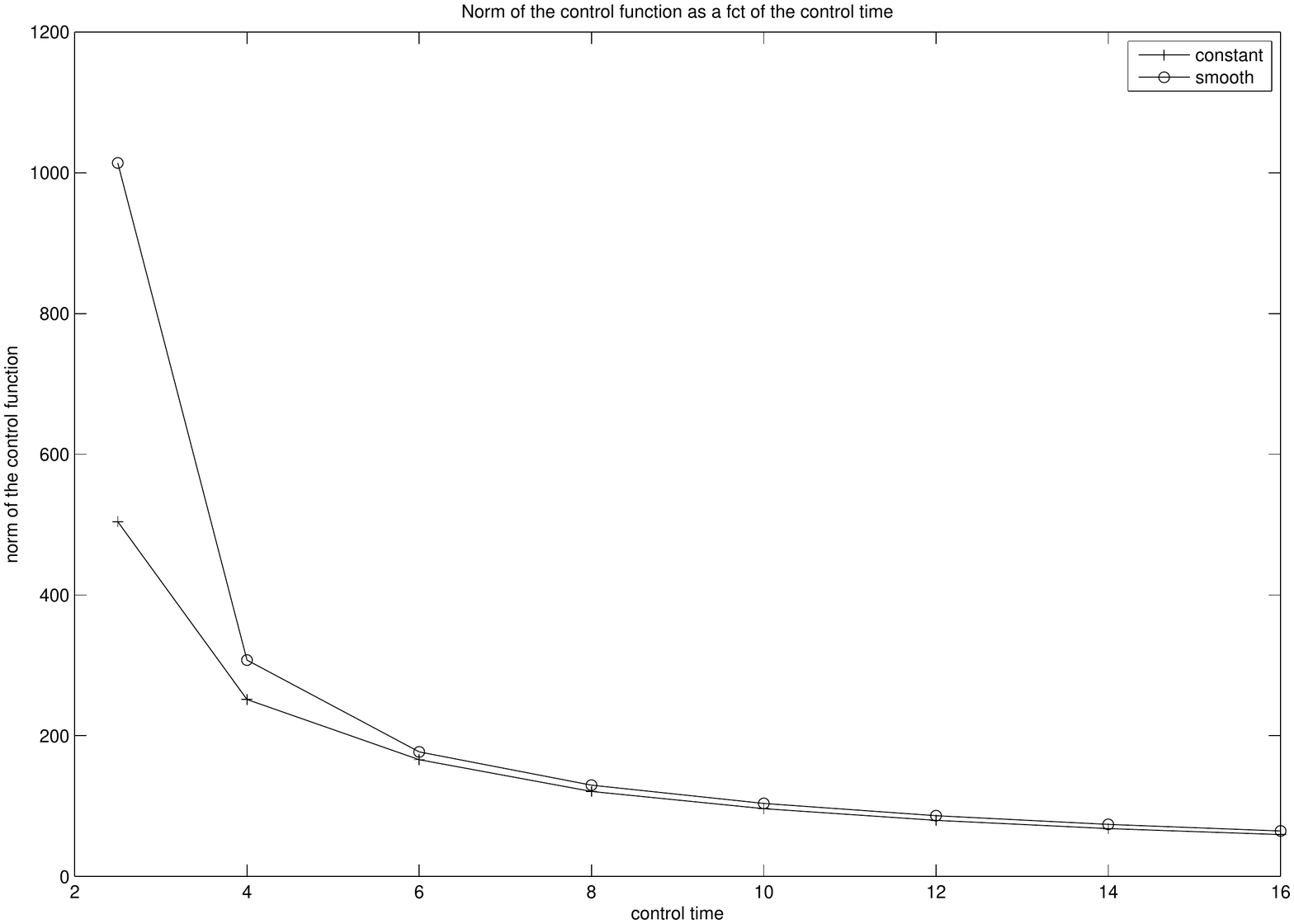}
\caption{\label{fig:examplestime1} 
Experiments in the square, with exact eigenvalues: impact of the smoothing on the reconstruction error (top) and on the norm of the control function (bottom), as a function of the control time. }
\end{center}
\end{figure}

\begin{figure}
\begin{center}
\includegraphics[width=\textwidth]{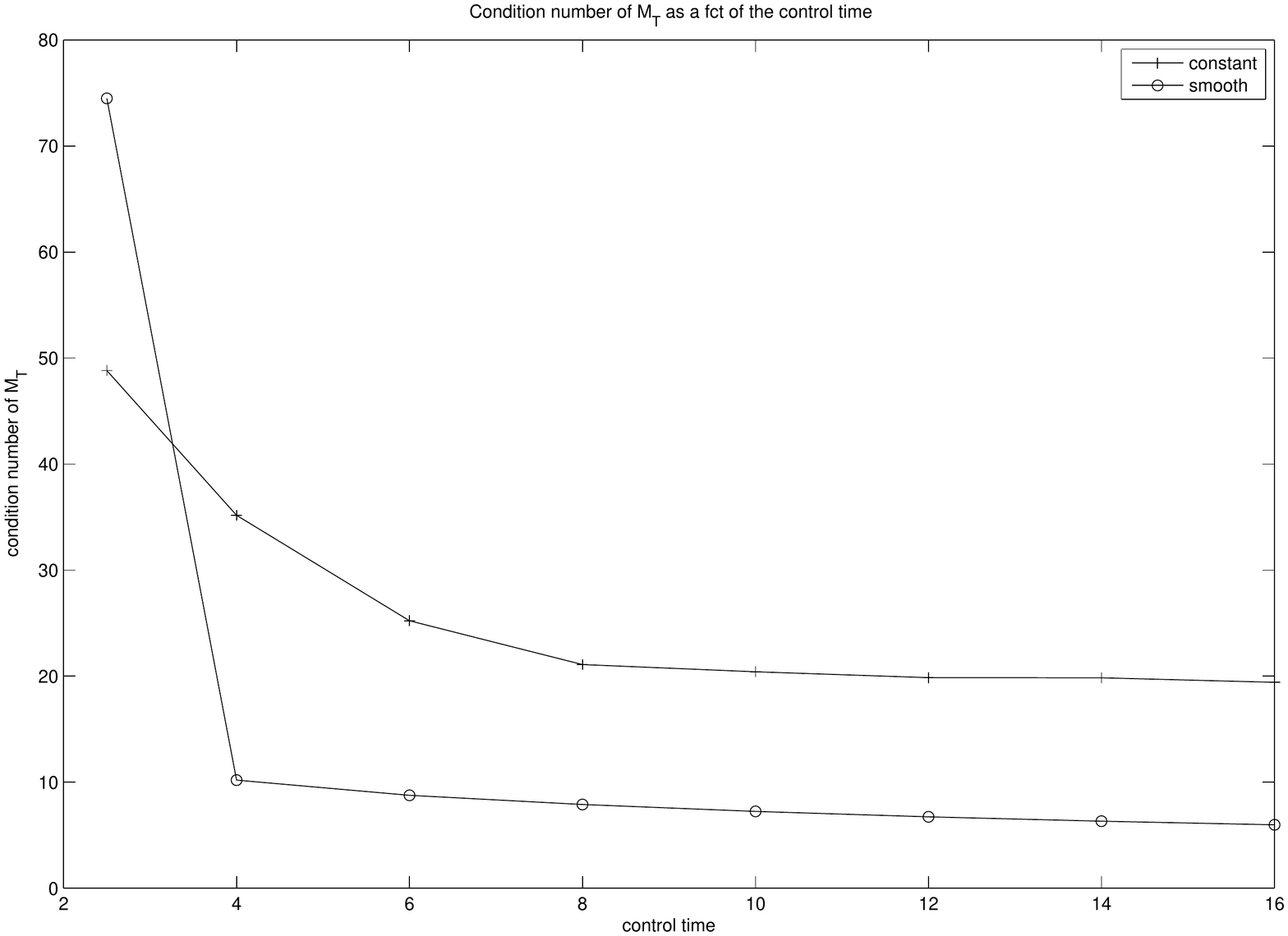}
\caption{\label{fig:examplestime2} 
Experiments in the square, with exact eigenvalues: impact of the smoothing on the condition number of the $M_{T}$ matrix, as a function of the control time. }
\end{center}
\end{figure}

In this paragraph we investigate the impact of the smoothing, the width of the control domain $U$ and the control time $T$ on various outputs such as the condition number of $M_{T}$, the reconstruction error $\|u-y\|$, and the norm of the control function $\|w\|$.\\
To do so we performed several one-mode experiments (see paragraph \ref{sec:onemode}, mode 500) in the square, with exact eigenvalues, 1000 eigenvalues used for computation of $M_{T}$, 2000 eigenvalues used for reconstruction and verification. We chose various times: 2.5 and 8, plus their ``smoothed" counterparts, according to the empirical formula $T_{\textrm{smooth}} = 15/8 * T$. This increase of $T_{\textrm{smooth}}$ is justified on the theoretical level by formulas (\ref{t5}) and (\ref{t4}) which show that the efficiency of the control is related to a mean value of $\chi(t,x)$ on the trajectories.  Similarly, we chose various width of $U $: 1/10 and 3/10, plus their ``smoothed" counterpart, which are double. Table \ref{tab:examples} presents the numerical results for these experiments. This table draws several remarks. First, the condition number of $M_{T}$, the reconstruction error and the norm of the control $w$ decrease with increasing time and $U$. Second, if we compare each non-smooth experiment with its ``smoothed" counterpart (the comparison is of course approximate, since the ``smoothed" time and width formulas are only reasonable approximations), the condition number seems similar, as well as the norm of the control function $w$, whereas the reconstruction error is far smaller with smoothing than without. \\
\begin{table}
\begin{center}
\begin{tabular}{|c|c|c|c|c|c|}
\hline Smooth & Width of $U $ & Time & Condition number & Rec. error & $\|w\|$\\
\hline no      &   1/10      & 2.5  &   48.8303   &   0.00518843    &    504.287      \\
	   yes     &   2/10      & 2.5  &   204.048   &   0.00140275    &    1837.12     \\
	   yes     &   2/10      & 4.7  &   21.7869   &   4.65892E-07     &    364.013       \\
\hline no      &   1/10      & 8  &   21.1003   &  0.00162583      &    120.744      \\
	   yes     &   2/10      & 8 &   16.5497   &   8.51442E-08     &    189.361     \\
	   yes     &   2/10      & 15  &  12.017    &  8.39923E-09     &     100.616     \\
\hline no      &   3/10      & 2.5  &  4.20741    &   0.0014823     &     147.009      \\
	   yes     &   6/10      & 2.5  &  9.05136     &  1.94519E-06      &   336.704       \\
	   yes     &   6/10      & 4.7  &  3.09927    &   2.99855E-08     &     125.481     \\
\hline no      &   3/10      & 8  &  3.20921    &  0.000488423     &    39.9988       \\
	   yes     &   6/10      & 8 &   2.74172    &  6.0204E-09     &    69.8206     \\
	   yes     &   6/10      & 15  &  2.4113     &    8.55119E-10     &    37.1463      \\
\hline
\end{tabular}
\caption{\label{tab:examples} 
Impact of the control time, the width of $U $ and the smoothing, on the condition number of $M_{T}$, on the reconstruction error and on the norm of the control function. These results come from one-mode experiments in the square, with exact eigenvalues.}
\end{center}
\end{table}
Figures \ref{fig:examplestime1} and \ref{fig:examplestime2} emphasize the impact of the control time, they present the reconstruction error, the norm of the control, and the condition number of $M_{T}$, as a function of the control time (varying between 2.5 and 16), with or without smoothing. Conclusions are similar to the table conclusions.

%%%%%%%%%%%%%%%%%%%%%%%%%%%%%%%%%%%%%%

\subsection{Condition number}\label{sec4.5}

\begin{figure}
\begin{center}
\includegraphics[width=1\textwidth]{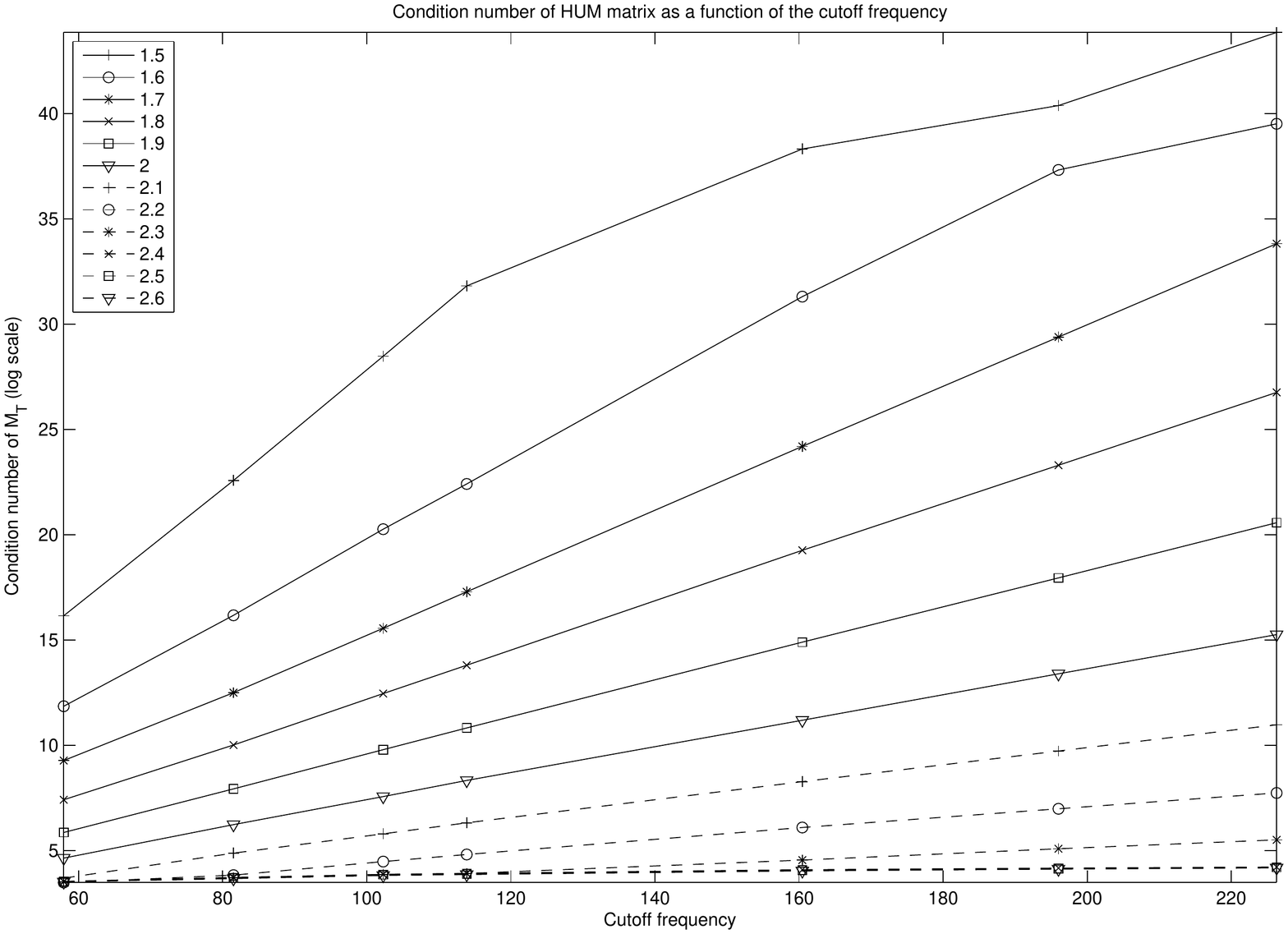}
\caption{\label{fig:cond-eigenvalues-constant} 
Condition number of the $M_{T,\omega}$ matrix as a function of the cutoff frequency $\omega$ for various control times. }
\end{center}
\end{figure}

\begin{figure}
\begin{center}
\includegraphics[width=1\textwidth]{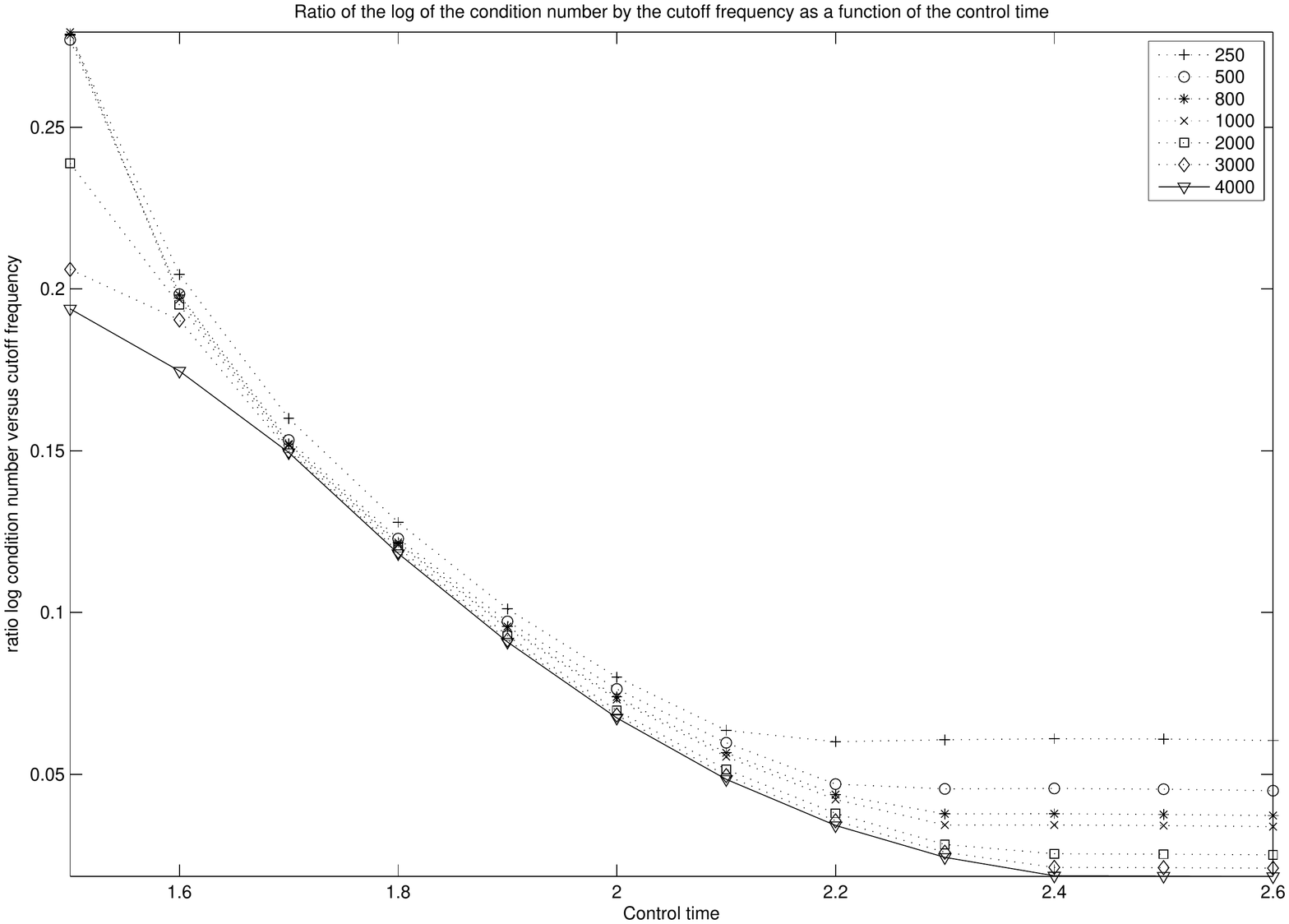}
\caption{\label{fig:cond-time-constant} 
Ratio of the log of the condition number of the $M_{T,\omega}$ matrix and the cutoff frequency $\omega$, as a function of the control time for various eigenvalues numbers. }
\end{center}
\end{figure}

In this section, we investigate conjecture \ref{conj2}. To do so, we compute the condition number of $M_{T,\omega}$, as we have:
\be
\textrm{cond}(M_{T,\omega}) = \|M_{T,\omega}\| .  \|M_{T,\omega}^{-1}\| \simeq \|M_{T}\| .  \|M_{T,\omega}^{-1}\|
\ee
Figure \ref{fig:cond-eigenvalues-constant} shows the condition number of the $M_{T}$ matrix as a function of the control time or of the last eigenvalue used for the control function computation. According to conjecture \ref{conj2}, we obtain lines of the type
\be
\log \left( \textrm{cond}(M_{T,\omega})\right) = \omega . C(T,U)
\ee
Figure \ref{fig:cond-time-constant} shows for various eigenvalues numbers the following curves:
\be
T \mapsto \frac{\log \left( \textrm{cond}(M_{T,\omega}) \right)}{\omega}
\ee
Similarly, we can draw conclusions compatible with conjecture \ref{conj2}, as these curves seems to converge when the number of eigenvalues grows to infinity.

%%%%%%%%%%%%%%%%%%%%%%%%%%%%%%%%%%%%%%

\subsection{Non-controlling domains}\label{sec4.6}

\begin{figure}
\begin{center}
\includegraphics[angle=-90,width=1\textwidth]{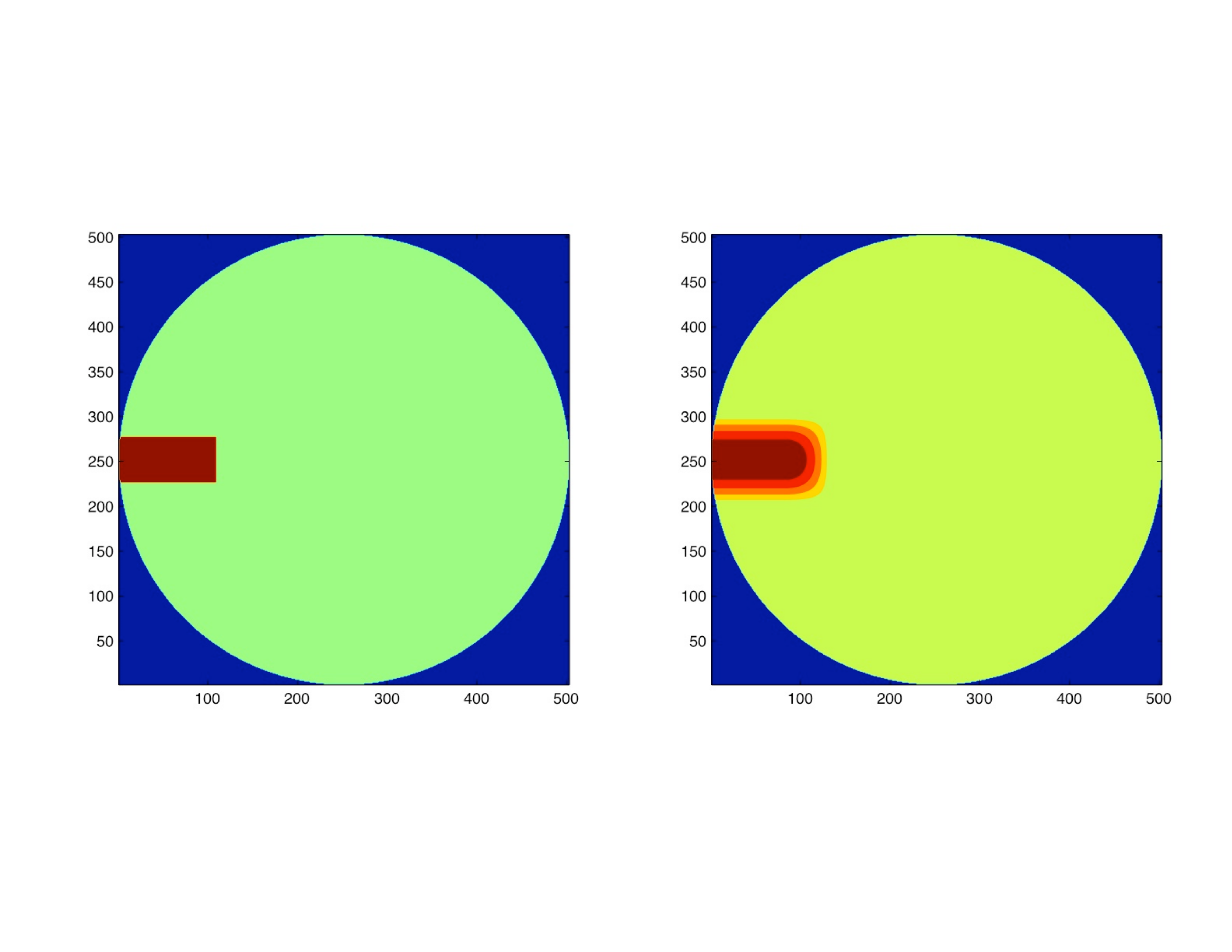}
\caption{\label{fig:omeg-spec} 
Non-controlling domain $U $ without (left) or with (right) smoothing. This domain consists of the neighborhood of a radius which is truncated around the disc boundary.}
\end{center}
\end{figure}

\begin{figure}
\begin{center}
\includegraphics[angle=-90,width=1\textwidth]{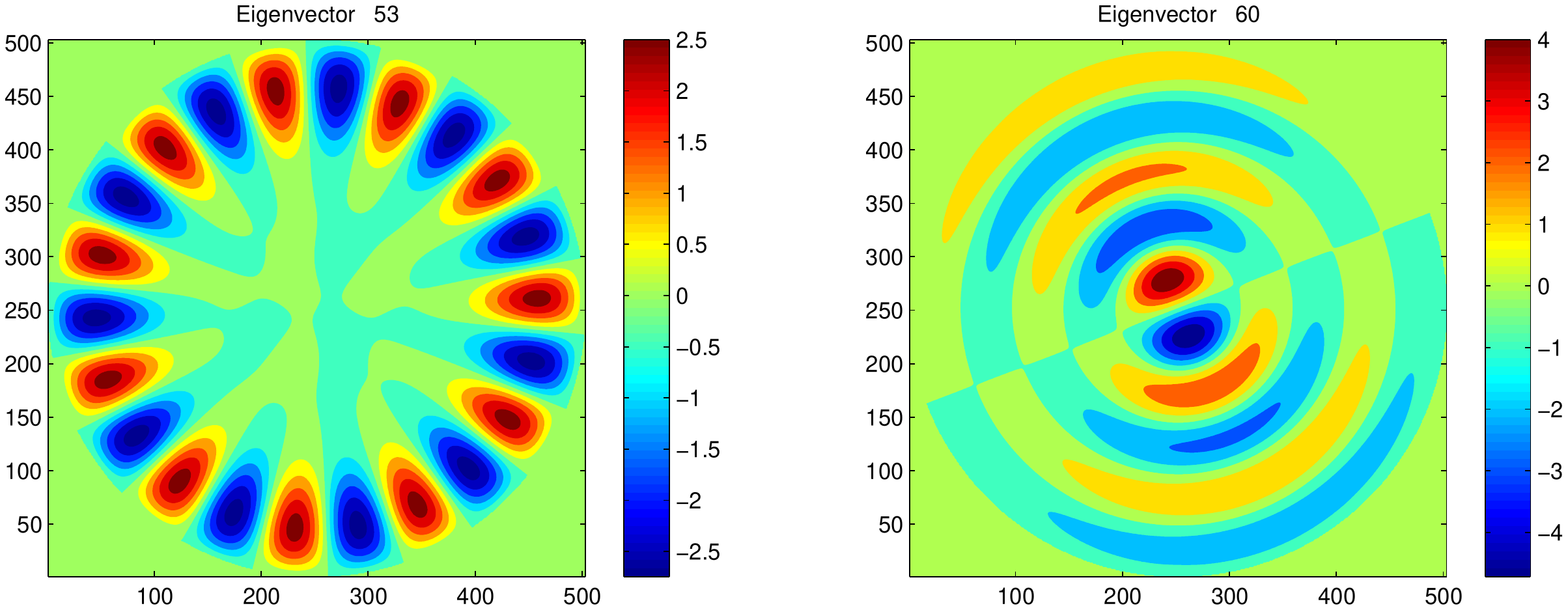}
\caption{\label{fig:modes-spec} 
Special modes chosen for experiment with non-controlling domains, corresponding to the 53rd and 60th eigenvalues. }
\end{center}
\end{figure}

In this section we investigate two special experiments with non-controlling domains, i.e. such that the geometric control condition is not satisfied whatever the control time.\\
First we consider the domain presented in Figure \ref{fig:omeg-spec}. 
For this domain the condition number of the $M_{T}$ matrix is large, and subsequently we should be experiencing difficulties to reconstruct the data $u$. We perform one-mode experiments with two different eigenvectors, one being localized in the center of the disc (eigenvalue 60), the other being localized around the boundary (eigenvalue 53)  as can be seen on Figure \ref{fig:modes-spec}. 
The various outputs are presented in Table \ref{tab:disc-special}, and we can see that the inversion is fairly accurate for the 53rd eigenmode, while it is logically poor for the 60th eigenmode. Moreover, the energy needed for the control process, i.e. the norm of the control $w$, is small for the 53rd eigenvector, while it is large for the 60-th. We can also notice that the smoothing has the noticeable effect to decrease the reconstruction error, the norm of the control function $w$ being similar.\\
\begin{table}
\begin{center}
\begin{tabular}{|c|c|cc|cc|}
\hline
         &  Condition & \multicolumn{2}{|c|}{53rd Eigenmode} & \multicolumn{2}{|c|}{60th Eigenmode}\\ 
smooth?  &     number & Rec. error &  $\|w\|$ &      Rec. error &  $\|w\|$\\
\hline no & 4960 & 0.63\% & 15  & 26\% & 4583\\
\hline yes & 6042 & 6.5 $10^{-4}$\% & 24  & 0.11\% & 8872\\
\hline
\end{tabular}
\caption{\label{tab:disc-special} 
Influence of the shape of the data on the reconstruction error and the norm of the control, with a non-controlling domain $U $. On the left, the 53rd eigenmode is localized around the boundary of the circle, as is the control domain. On the right, the 60th eigenmode is localized around the center of the circle.}
\end{center}
\end{table}

\begin{figure}
\begin{center}
\includegraphics[angle=-90,width=1\textwidth]{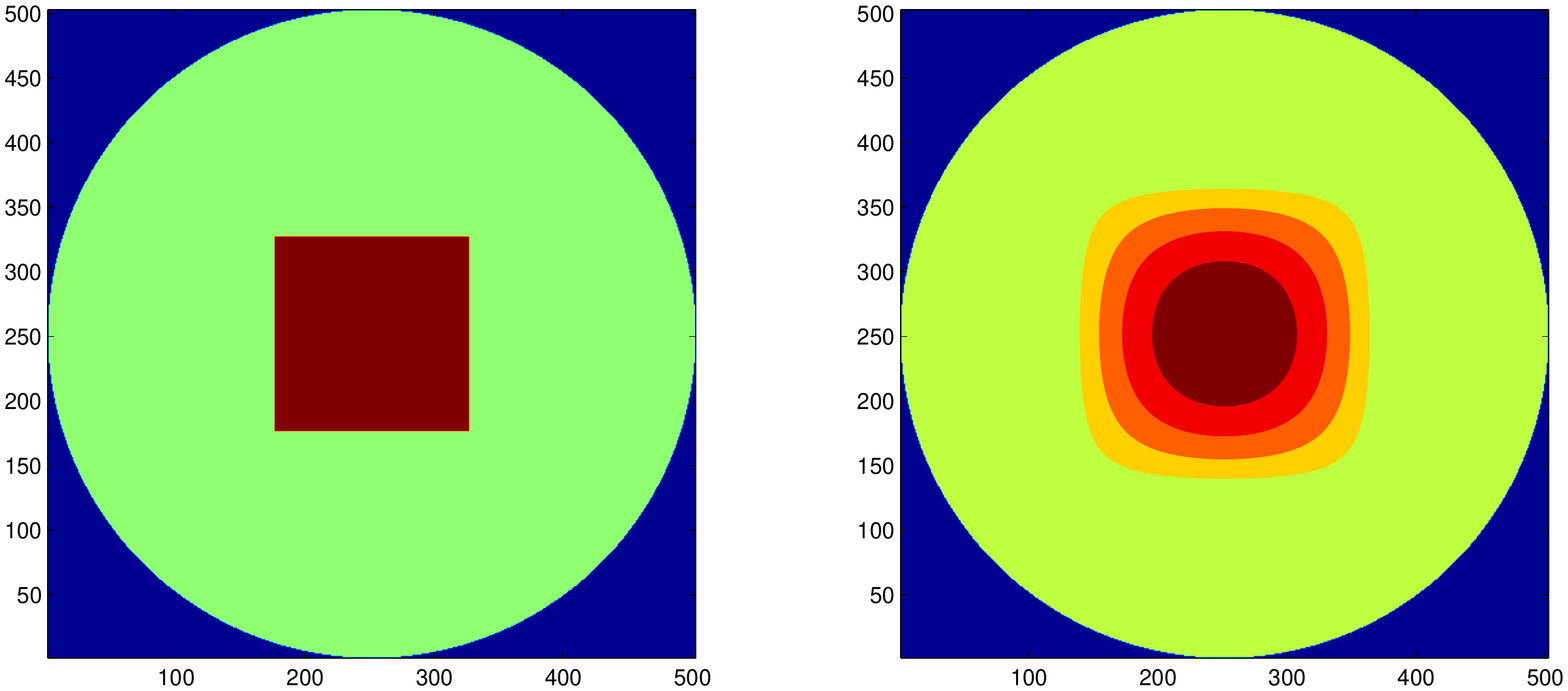}
\caption{\label{fig:omeg-spec2} 
Non-controlling domain $U $ without (left) or with (right) smoothing.  This domain consists of the neighborhood of a radius which is truncated around the disc center. }
\end{center}
\end{figure}

In the second experiment we change the point of view: instead of considering one given domain and two different data, we consider one given data, and two different non-controlling domains. The data is again $u_{53}$ (see Figure \ref{fig:modes-spec}), which is localized at the boundary of the disc. The first domain is the previous one (see Figure \ref{fig:omeg-spec}), the second domain is presented in Figure \ref{fig:omeg-spec2}, it is localized at the center of the disc. 
In either case, the condition number of the $M_{T}$ matrix is large, and the data should prove difficult to reconstruct. Table \ref{tab:disc-special2} present the outputs we get for the two domains. As previously, we observe that the control process works fairly well for the appropriate control domain, with a small error as well as a small energy for the control. Conversely, when the control domain does not ``see" the input data, the results are poorer: the energy needed is large with or without smoothing, the error is also large without smoothing, it is however small with smoothing. 
\begin{table}
\begin{center}
\begin{tabular}{|c|ccc|ccc|}
\hline
         &  \multicolumn{3}{|c|}{First domain} & \multicolumn{3}{|c|}{Second domain}\\ 
smooth?  &     Cond. nb. & Rec. error &  $\|w\|$ &     Cond. nb. & Rec. error &  $\|w\|$\\
\hline no & 4960 & 0.63\% & 15 & 3.6 $10^6$ & 68\% & 1.9 $10^4$\\
\hline yes & 6042 & 6.5 $10^{-4}$\% & 24 & 3.3 $10^5$ & 9.4 $10^{-3}$\% & 6.5 $10^5$\\
\hline
\end{tabular}
\caption{\label{tab:disc-special2} 
Reconstruction error and the norm of the control with a data $u_{53}$ localized at the boundary of the disc and two different non-controlling domains $U $, the first one being localized around the boundary, the second one around the center.}
\end{center}
\end{table}

%%%%%%%%%%%%%%%%%%%%%%%%%%%%%%%%%%%%%%

\section*{Acknowledgement}
\addcontentsline{toc}{section}{Acknowledgement}

The experiments have been realized with Matlab\footnote{The Mathworks, Inc. \url{http://www.mathworks.fr}} software on Laboratoire Jean-Alexandre Dieudonn\'e (Nice) and Laboratoire Jean Kuntzmann (Grenoble) computing machines. The INRIA Gforge\footnote{http://gforge.inria.fr} has also been used. The authors thank J.-M. Lacroix (Laboratoire J.-A. Dieudonn\'e) for his managing of Nice computing machine. \\
This work has been partially supported by Institut Universitaire de France.

%%%%%%%%%%%%%%%%%%%%%%%%%%%%
%%%%%%%%%%%%%%%%%%%%%%%%%%%%
\addcontentsline{toc}{section}{References}
\bibliography{Bibli_wavecontrol}

\end{document}